\documentclass[12pt]{article}
\usepackage{graphicx} 
\usepackage{mathrsfs}
\usepackage{amsfonts,amscd,amsthm}
\usepackage{amsmath}
\usepackage{amssymb}
\usepackage{bbm}
\usepackage{graphicx}
\usepackage{enumerate}
\usepackage{geometry}
\usepackage{xcolor}
\usepackage{txfonts}
\usepackage{cite}
\usepackage{bbding}
\newtheorem{definition}{Definition}[section]
\newtheorem{theorem}{Theorem}[section]
\newtheorem{corollary}{Corollary}[section]
\newtheorem{lemma}{Lemma}[section]

\newtheorem{proposition}{Proposition}[section]
\newtheorem{remark}{Remark}[section]
\newtheorem{example}{Example}[section]

\newcommand{\B}{{\mathbb{B}}}

\newcommand{\R}{\mathbbm{R}}

\DeclareMathOperator*\bdry{bdry}
\DeclareMathOperator*\argmax{arg\,max}

\DeclareMathOperator*\co{conv}

\DeclareMathOperator*\cl{cl}

\DeclareMathOperator*\inte{int}
\DeclareMathOperator*\ri{ri}

\DeclareMathOperator*\dom{dom}
\DeclareMathOperator*\rge{rge}

\DeclareMathOperator*\extr{extr}
\DeclareMathOperator*\gph{gph}

\DeclareMathOperator*\rank{rank}
\DeclareMathOperator*\cone{cone}
\DeclareMathOperator*\aff{aff}
\DeclareMathOperator*\proj{proj}
\usepackage{comment}

\begin{document}

\begin{center}
{\Large\Large\sc  {\bf    Lipschitz continuity of solution multifunctions of extended $\ell_1$ regularization problems}}
\end{center}

\begin{center}

{\sc Kaiwen Meng}
\\ {\small School of  Mathematics, Southwestern University of Finance and Economics, Chengdu 611130, China}\\
Email: mengkw@swufe.edu.cn\\[0.5cm]

{\sc  Pengcheng Wu}
\\ {\small Department of Applied Mathematics, The Hong Kong Polytechnic University, Hong Kong }\\
Email: pcwu0725@163.com \\[0.5cm]

{\sc Xiaoqi Yang}\\ {\small  Department of Applied Mathematics, The Hong Kong Polytechnic University,  Hong Kong}\\
Email: mayangxq@polyu.edu.hk

\end{center}

{\bf\noindent  Abstract:}  The Lasso and the basis pursuit in compressed sensing and machine learning are convex optimization problems with three parameters: the regularization scalar, the observation vector and the data matrix. Relative to the first two parameters, we obtain the Lipschitz continuity of the solution multifunction on its convex domain. When the data matrix of the Lasso also perturbs, where non-polyhedral structure may display, we obtain full characterizations for the Lipschitz continuity of the solution multifunction on the product of a compact and convex set in the space of first two parameters and a neighborhood of the fixed data matrix. Moreover for the solution multifunction of the Lasso, we show that the Lipschitz continuity implies its single-valuedness. Our analysis is based on polyhedron theory, a sufficient condition that ensures the Lipschitz continuity of a polyhedral multifunction with a convex domain, and an explicit representation of the solution multifunction, where the latter is a consequence of the Lipschitz continuity of the solution multifunction relative to the first two parameters.

 \vskip0.1cm



\vskip0.1cm


\section{Introduction}

In this paper, we consider the Lipschitz continuity of the optimal solution multifunction of the following extended $\ell_1$ regularization problem:
\begin{equation}\label{ext-lasso-problem-lambda}
\min_{x \in \R^n} \; \frac{1}{2\lambda}\|Ax-\proj{_{\rge A}}b\|^2 + \|x\|_1,
\end{equation}
where $(\lambda,b,A) \in \mathbb{R_+}\times \mathbb{R}^n\times \mathbb{R}^{m\times n}$ are parameters and $\proj_{\rge A}b$ denotes the projection of vector $b$ onto the range space  of matrix $A$ ($\rge A$ in short). 

When $\lambda > 0$, the extended $\ell_1$ regularization problem \eqref{ext-lasso-problem-lambda} and the following Lasso 
\begin{equation}\label{lasso-problem-lambda-b1}
\min_{x \in \R^n} \;  \frac12\|Ax-b\|^2+\lambda \|x\|_1,
\end{equation}
have the same set of optimal solutions. 
When $\lambda = 0$, the extended $\ell_1$ regularization problem 
\eqref{ext-lasso-problem-lambda} can be seen from the viewpoint of penalization as the following modified Basis Pursuit (BP in short) problem 
\begin{equation}\label{ext-BP-problem}
\min_{x \in \mathbb{R}^n}\;\|x\|_1\quad\mbox{subject to}\quad Ax-\mbox{proj}_{\rge A}b=0,
\end{equation}
and in this case  the optimal solution set of problem (\ref{ext-lasso-problem-lambda}) is defined  to be that of (\ref{ext-BP-problem}). 
The modified BP is seen as an extension of the classical BP by replacing observation data vector $b$ with its projection $\proj_{\rge A}b$. This variation is to account for the non-feasibility of the classical BP when $A$ is not of full row rank.

The Lasso problem and Basis Pursuit problem  have attracted extensive attention from optimization, statistics and computer science communities and have been accepted as two of the most useful approximate models for finding a sparse solution (see Tibshirani \cite{Tibshirani13}, Candes et al. \cite{CandesRT06} and Donoho \cite{Donoho06} and the references therein). Many effective algorithms have been designed for solving these problems. However the study of stability properties of the optimal solution multifunction of (\ref{lasso-problem-lambda-b1}) and (\ref{ext-BP-problem}) is very limited. The recovery bound or so-called relative upper Lipschitz continuity property of optimal solutions of problem (\ref{lasso-problem-lambda-b1}) with respect to parameter $\lambda$ only was obtained under various assumptions such as restricted isometric property and restricted eigenvalue condition on data matrix $A$, see Hu et al. \cite{Hu17} and the references therein. 

Stability properties of multifunctions have extensively been investigated, see {\cite{Adly2014,AubANA1984,bonnans2013perturbation,dontchev2009implicit,Rock2009VaAn,outrata1998KZ,GfrOut2016,Meng2020} } and the references therein. Robinson \cite{Robinson1981} obtained the local upper Lipschitz continuity of a polyhedral multifunction (i.e., its graph is expressed as the union of finitely many polyhedral sets). Mangasarian and Shiau \cite{Mangasarian1987} established the Lipschitz continuity for solutions of linear inequalities, linear programs and certain linear complementarity problems.
Rockafellar and Wets \cite[Example 9.35]{Rock2009VaAn} (see also Walkup and Wets \cite{walkup1969w}) obtained the Lipschitz continuity for a multifunction with graph being one polyhedron.
In \cite{Han2024P} Han and Pang derived the (Lipschitz) continuous (single-valued) solution
function of parametric variational inequalities under functional and constraint perturbations.
Mordukhovich et al. \cite{Mordukhovich2014OS} established characterizations of the single-valuedness and local Lipschitz continuity of the local optimal solutions for parametric second-order cone programs. Beer et al. \cite{Beer2021CLP} analyzed the Lipschitz modulus of the feasible set of linear and convex inequality systems. Under a strong regularity assumption Guo and Xu \cite{Guo2024X} obtained the Lipschitz continuity of isolated solutions of perturbed stochastic generalized equations around the given point.
It is worth noting that for a polyhedral multifunction with convex domain, Gowda and Sznajder \cite{Gowda1996SLip} showed that it is Lipschitz continuous on the domain if and only if the inverse of the multifunction is open. This seemly complete characterization of Lipschitz continuity may not be easily applied, as the verification of the openness of the inverse of the multifunction can be as difficult as verifying the Lipschitz continuity of the multifunction itself. 

In the recent years, the uniqueness and the local  Lipschitz continuity of the optimal solution multifunction of \eqref{lasso-problem-lambda-b1} and its extensions have been studied by employing some variational conditions. Define the multifunction $S:\R\times \R^m \times \R^{m\times n}\rightrightarrows \R^n$ as follows,  
for any triple parameters $(\lambda, b, A)\in \R\times \R^m \times \R^{m\times n}$,  
  \[
S(\lambda,b, A):=\left\{
  \begin{array}{cl}
\mbox{the optimal solution set of}\; (\ref{ext-lasso-problem-lambda}) &\mbox{if}\;\lambda\geq 0,\\[0.1cm]
  \emptyset   &\mbox{if}\;\lambda<0.  
  \end{array}
      \right.
 \] 
It is clear that $\dom S$ (i.e., the domain  of $S$) is   equal to $\R_+ \times \R^m\times \R^{m\times n}$ and thus closed, which is in some sense an appealing requirement for its Lipschitz behavior. In view of the equivalence 
$$
Ax-\proj{_{\rge A}}b=0\ \Longleftrightarrow \  A^\top (Ax- b)=0, 
$$
we get from first-order optimality conditions of  \eqref{lasso-problem-lambda-b1} and \eqref{ext-BP-problem} that for any $(\lambda, b, A)\in  \dom S$,
\begin{equation}\label{zyx-lambda}
 S(\lambda,b, A)=\left\{
  \begin{array}{ll}
\left\{x\in \R^n\mid   0 \in   A^\top (Ax- b)+ \lambda \partial  \|\cdot\|_1(x) \right\}  &\mbox{if}\; \lambda>0,  \\[0.2cm]
 \left \{x\in \R^n\mid A^\top (Ax- b)=0, \;0\in \rge A^\top  + \partial \|\cdot\|_1(x) \right\} &\mbox{if}\;\lambda=0.
  \end{array}
  \right.
\end{equation}
Letting $I:=\{1,\cdots, n\}$, and for every vector  $x\in \R^n$, we define the following index sets: 
 \begin{equation}\label{def-index-partition-signs-x}
     I^+(x):=\{i\in I\mid     x_i>0\},\; I^0(x):=\{i\in I\mid     x_i=0\},\; I^-(x):=\{i\in I\mid    x_i<0\}.
 \end{equation}
Let us recall the conditions that are used in \cite{Berk2023bh,Cui2024,Nghia2024} for the study of the uniqueness and local Lipschitz continuity of $S$. For every $( \lambda, b, A)\in \dom S$, we define the following conditions.
\begin{description}
    \item[Condition 1.1] There are some ${x}\in S( \lambda,   b,A)$ and  $y\in \R^m$  such that the columns $A_i$ of $A$ with $i\in I^+( {x})\cup I^-( {x})$ are linearly independent  and
    \[
              A_i^\top y=1\;\forall i\in I^+( {x}),\;-1<A_i^\top y<1\;\forall i\in I^0( {x}),\;A_i^\top y=-1\;\forall i\in I^-( {x}).
        \]
        \item[Condition 1.2] There is some $ {x}\in S( \lambda,   b,A)$ such that the columns $A_i$ of  $A$ with $i\in I$ such that $|A_i^\top ( {b}-A {x})|=\lambda$ are linearly independent.
    \item[Condition 1.3] There is some $ {x}\in S( \lambda,   b,A)$ such that the columns $A_i$ of  $A$ with $i\in I^+( {x})\cup I^-( {x})$ are linearly independent  and
     \[
    |A_i^\top ( {b}-A {x})|< {\lambda}\quad \forall i\in I^0( {x}).
   \]
\end{description}

While Condition  1.1 was  proposed in \cite{zhang2015yc} in a later time, {Conditions 1.2 and 1.3 have  appeared earlier  in  \cite{Tibshirani13,Mairal2012}  and  \cite{fuchs04},} respectively. { Sensitivity analysis of solution mapping for general regularization problems under Condition 1.3 was studied in \cite{Vaiter2017}.}
It has been shown in \cite[Lemma 4.6]{Berk2023bh} that the following implications hold for the case of $\lambda>0$:
\begin{equation} \label{impls}
    \mbox{Condition 1.3}\Longrightarrow \mbox{Condition 1.2}\Longrightarrow  \mbox{Condition 1.1}.
\end{equation}


Condition 1.1 guarantees the uniqueness of the solution of \eqref{lasso-problem-lambda-b1}, see \cite{zhang2015yc}. Under Condition 1.2 (resp. Condition 1.3), Berk et al. \cite{Berk2023bh} established local Lipschitz continuity (resp. smoothness) property of the unique optimal solution of the Lasso problem (\ref{lasso-problem-lambda-b1}) with $(\lambda,b) \in \R_{++} \times \R^m$ and local Lipschitz continuity of that with $(\lambda,b,A) \in \R_{++} \times \R^m \times \R^{m\times n}$ respectively, where the coderivative and variational analysis tools have been employed in their study. 
{Nghia \cite{Nghia2024} obtained full characterizations for local (single-valued) Lipschitz continuity not only for a group Lasso (including the Lasso) but also for the nuclear norm minimization problem (with non-polyhedral structures) respectively. Cui et al. \cite{Cui2024} considered a broader class of models where they replaced the $\ell_1$-norm in the Lasso with a general convex function having (Fenchel) conjugate that is $\mathcal{C}^2$-cone reducible by employing Robinson's strong regularity on the dual problem.} The equivalences among Condition 1.2, the single-valuedness and Lipschitz continuity of $S$ around $(\bar\lambda, \bar b,\bar A)$, and the single-valuedness of $S$ around $(\bar\lambda, \bar b,\bar A)$ are obtained in \cite{Cui2024}, {see also \cite{Berk2023bh,Nghia2024}}. In these works \cite{Berk2023bh,Cui2024,Nghia2024}, the Robinson's strong regularity in \cite{robinson1980strongly} was employed to obtain the local Lipschitz continuity. {After our initial submission of this paper, we notice that Hu et al. \cite{zhang-arxiv2024} also obtained the Lipschitz continuity for the Lasso.} 



In this paper we will employ the polyhedral theory to establish the Lipschitz continuity of the optimal solution multifunction $S$ of problem (\ref{ext-lasso-problem-lambda}). We will investigate the following two cases 
\begin{enumerate}
    \item[(i)] the matrix $A$ is fixed, and
    \item[(ii)] the matrix $A$  varies  locally around a fixed matrix $\bar A$.
\end{enumerate}

For the case (i), we obtain the Lipschitz continuity of $S(\cdot,\cdot,A)$ on its domain $\R_+ \times \R^m$. This is done by observing a verifiable sufficient condition for a polyhedral multifunction to be Lipschitz continuous on its domain. In this sufficient condition, we impose two assumptions: the first one is that its domain is convex and the second one is that the graph of the restriction of the multifunction on the projection of each of the  polyhedral sets (in the polyhedral multifunction decomposition) onto the variable space is convex. The convexity assumption of the domain allows one to connect any two points in the domain by a straight line, while this straight line can be divided into a finite number of line segments and each line segment lies in one of the polyhedral sets on which the multifunction is Lipschitz continuous, as the result of the second assumption. The Lipschitz continuity of the multifunction on its domain is acquired by that of these arbitrary line segments with a Lipschitz constant being the maximum of the Lipschitz constants corresponding to the finitely many polyhedral sets. The remaining task is to verify the two assumptions in the sufficient condition. In our approach, we employ    a collection   of index partitions  of $I$, which correspond  one-to-one to the non-empty faces of the polyhedral feasible set of the following dual problem of the extended $\ell_1$ regularization problem (\ref{ext-lasso-problem-lambda}):   
\begin{equation}\label{lasso-dual0001}
\max \;\langle \proj{_{\rge A}} b, y\rangle-\frac12 \lambda \|y\|^2\quad\mbox{subject to}\quad    \|A^\top y\|_\infty\leq 1.
\end{equation}
We show that the graph of $S(\cdot,\cdot,A)$ is expressible  as the union of some polyhedral cones which are constructed by these index partitions of $I$. Thus $S(\cdot,\cdot,A)$ is verified to be a polyhedral multifunction. Moreover, we establish that the domain of   $S(\cdot,\cdot,A)$ is partitioned as the union of the projections of the above polyhedral cones onto the parameters space. We then prove that the graph of the restriction of $S(\cdot,\cdot,A)$ on each of these projections is convex, where we employ the single-valuedness property of the mapping $(\lambda,b) \to Ax$ on the solution multifunction  $S(\lambda,b,A)$ as established in \cite{zhang2015yc}.   

For the case (ii), as a consequence of the newly established Lipschitz continuity of $S(\cdot,\cdot,A)$ on its domain, we begin by presenting an explicit representation of $S$. 
We then apply this explicit representation for $S$ 
to investigate full characterizations of Conditions 1.2 and 1.3 at some $(\bar\lambda, \bar b,\bar A)\in \dom S$ with $\bar\lambda>0$. In particular, we show the equivalences among Condition 1.2 at $(\bar\lambda, \bar b,\bar A)$, the Lipschitz continuity of $S$ around $(\bar\lambda, \bar b,\bar A)$ and the piecewise linearity of $S(\cdot,\cdot,\bar A)$ around $(\bar\lambda, \bar b)$. As for Condition 1.3, we demonstrate that it holds at $( \bar\lambda, \bar b, \bar A)$ if and only if $S$ is $C^\infty$ around $(\bar\lambda, \bar b,\bar A)$ or equivalently $S(\cdot, \cdot, \bar A)$ is linear around $(\bar \lambda, \bar b)$. These equivalences improve the ones in \cite{Berk2023bh,Cui2024} in a sense that the Lipschitz continuity of $S$ around $(\bar \lambda,\bar b,\bar A)$ can actually imply the single-valuedness of $S$ around $(\bar \lambda,\bar b,\bar A)$. While a counterexample shows that the global Lipschitz continuity of $S$ on its domain is not possible when $A$ varies, we establish that $S$ is Lipschitz continuous on the product of a nonempty, compact, and convex subset $\mathcal{U}$ of $\mathbb{R}_{++}\times \mathbb{R}^m$ and a neighborhood of $\bar A$ if and only if Condition 1.2 holds at $(\lambda, b, \bar A)$ for all  $(\lambda, b)\in \mathcal{U}$.  Furthermore, when $\bar{A}$ is of full rank, we are able to quantify the radius of the ball neighborhood of $\bar A$. 

The rest of the paper is organized as follows. Section 2.1 shows a sufficient condition for a polyhedral multifunction to be Lipschitz continuous on its convex domain. Section 2.2 includes several technical lemmas. When $A$ is fixed, Section 3 establishes the Lispchitz continuity of the solution multifunction $S$ on its convex domain and derives, as a corollary, an explicit representation of $S$. Section 4 obtains characterizations of Conditions 1.2 and 1.3 and establishes the Lipschitz continuity of $S$ when $A$ varies locally. Section 5 presents conclusions of the paper.

\section{Preliminaries and technical lemmas}\label{sec-preli-tech-lem}

 Throughout the paper we use the standard notations of variational analysis; see the seminal book \cite{Rock2009VaAn} by Rockafellar and Wets.
The Euclidean $\ell_2$-norm, $\ell_1$-norm and $\ell_\infty$-norm of a vector $x \in \R^n$ are denoted by $||x||$, $\|x\|_1$ and $\|x\|_\infty$ respectively. Let $\B(x,\alpha)$ denote the closed  ball centered at $x$ with radius $\alpha$, $\B$ denote the closed unit Euclidean ball,    $\mathbb R_+:=[0,\infty)$ and $\mathbb R_{++}:=(0,\infty)$.  For $x,y\in\mathbb R^n$, we define $[x,y]:=\{x+\theta(y-x)\mid \theta\in[0,1]\}$ and $(x,y):=\{x+\theta(y-x)\mid \theta\in(0,1)\}$.

Let $C\subset \R^n$ be a set.   We denote the interior, the relative interior, the boundary, the closure and the affine hull of $C$ respectively by $\inte C$, $\ri C$, $\bdry C$, $\cl C$ and $\aff C$.  Moreover, we denote by $|C|$ the cardinality of the set $C$, i.e., the number of elements in $C$.

A polyhedral set $C$ is defined as the intersection of a finite family of
closed half-spaces in $\R^n$, i.e., as the set of solutions to some finite system of inequalities of the form
\begin{equation}\label{polydef}
\alpha_j^\top x\leq\beta_j,\quad j=1,\cdots,k,
\end{equation}
for $\alpha_1,\cdots,\alpha_k\in \mathbb R^n$ and $\beta_1,\cdots,\beta_k\in \mathbb R$.
For a nonempty subset $F$ of polyhedral set $C$, $F$ is said to be a face of $C$ if there exists a vector $v\in \R^n$ (possibly $v=0$) such that 
$
F=\argmax\{\;v^\top x \mid x\in C\}.
$
The zero-dimensional
faces of $C$ are called the extreme points of $C$ and the collection of all such points is denoted as $\extr (C)$. In general, if $F$ is a half-line face
of $C$, we shall call the direction of $F$ an extreme direction of $C$.
Moreover, we deduce from \cite[Page 101]{S86TheoryOL} that $F$ is a nonempty  face of $C$ represented by (\ref{polydef}) if and only if there exists $J\subset\{1,\cdots,k\}$ such that 
$
F=\{x\in C\mid \alpha_j^\top x=\beta_j \quad\forall j\in J\} \not= \emptyset.
$

Let $A \in \mathbb{R}^{m\times n}$ be an $m\times n$ matrix, $x\in\mathbb R^n$, $i\in I:=\{1,\cdots,n\}$ and $K\subset I$. We denote that $A_i$ is the $i$th column of $A$ and $A_K\in\mathbb R^{m\times|K|}$ is the submatrix of $A$ consisting of all $A_i$'s with $i\in K$. We also denote that $x_K\in\mathbb R^{|K|}$ is the subvector of $x$ consisting of all $x_i$'s with $i\in K$. Denote the minimum singular value and spectral  norm of $A$ as $\sigma_{\min}(A)$ and $\|A\|,$ respectively.

For a multifunction $G:\R^n\rightrightarrows \R^m$, we denote by
$$
\gph G:=\{(x, y)\mid y\in G(x)\},\ \dom G:=\{x\mid G(x)\not=\emptyset\}\ \mbox{and}\ \rge G:=\{y\mid \exists x\in\dom G~\mbox{s.t.}~ y\in G(x)\}
$$
the graph, the domain and the range of $G$, respectively. For a set $\Omega\subset \R^n$, we denote by
  \[
  G|_\Omega(x):=\left\{
  \begin{array}{ll}
    G(x) & \mbox{if}\;x\in \Omega, \\[0.2cm]
    \emptyset & \mbox{if}\;x\not\in \Omega,
  \end{array}
    \right.
  \]
  the multifunction of $G$ restricted on $\Omega$. It is clear to see that
  \[
  \gph G|_\Omega=\gph G\cap (\Omega\times \R^m)\quad\mbox{and}\quad \dom  G|_\Omega=\Omega\cap \dom G.
  \]
A multifunction $G:\R^n\rightrightarrows \R^m$ is said to be a polyhedral one if its graph is the union of finitely many polyhedra in $\mathbb R^n\times \mathbb R^m$. Moreover, for a polyhedral multifunction $G$, if it is further single-valued, we called it a piecewise linear function (cf. \cite[Definition 2.47]{Rock2009VaAn}). Let $g:\R^n\to\R^m$ be single-valued. We say that $g$ is linear  if there exists some matrix $B\in\R^{m\times n}$ such that $g(x)=Bx$ and that $g$ is linear on some nonempty subset $D$ of $\R^n$ if there exists a linear function $\tilde g$ such that $g|_D=\tilde g|_D$.

The definition below is borrowed from \cite[Definition  9.26]{Rock2009VaAn}.
 \begin{definition}[Lipschitz continuity of multifunctions]\label{def-lip-con-setvalued}
 A multifunction $G:\R^n\rightrightarrows \R^m$ is Lipschitz continuous on $X$, a subset of $\R^n$, if it is nonempty-closed-valued on $X$ and there exists $\kappa\in \R_+$, called a Lipschitz constant, such that
 \[
 G(x')\subset G(x)+\kappa \|x'-x\|\B \quad \forall x',x\in X.
 \]
Given a point  $x\in \dom G$, we say that  $G$ is Lipschitz continuous  around   $x$ (resp. relative to $X$) if there is a neighborhood $V$ of $x$ such that $G$ is  Lipschitz continuous on $V$ (resp. $V\cap X$). 
\end{definition}

Let $g:\R^n\to \R$ be a convex function and let $\bar{x} \in \R^n$ be  a point.
The vector $v\in \R^n$ is a   subgradient of $g$ at $\bar{x}$, written $v\in \partial g(\bar{x})$, if $g(x)\geq g(\bar{x})+ v^\top (x-\bar{x})$ for all $x\in\mathbb R^n$.

\subsection{A sufficient condition for Lipschitz continuity of a polyhedral multifunction}

In this subsection, we establish a verifiable sufficient condition for the Lipschitz continuity of a polyhedral multifunction $G$ over its domain. Our approach relies on two key aspects. First, we utilize the decomposition of the graph of the polyhedral multifunction $\gph G$ into the union of certain polyhedral sets. Second, we consider the convexity of the graph of the restricted function of $G$, where the restriction is defined on the projection of each polyhedral set from the decomposition of $\gph G$ onto the decision variable space. 

{The following Lemma \ref{lem-Lip-con-poly-theo-ori} is an extension of Walkup and Wets’s theorem \cite{walkup1969w} (see also \cite[Example 9.35]{Rock2009VaAn}) showing that any polyhedral convex mapping is Lipschitz continuous relative to its domain.} Furthermore, it was demonstrated in \cite{facchinei03,lu2008variational,scholtes2012introduction} that when a polyhedral multifunction $G$ is single-valued on a convex subset $C$ of $\dom G$, then $G$ exhibits Lipschitz continuity (and is piecewise affine) on $C$. 

\begin{lemma} [Lipschitz continuity of polyhedral multifunctions] \label{lem-Lip-con-poly-theo-ori}
Consider a polyhedral multifunction $G$, whose graph   can be expressed  as  the union of finitely many polyhedral sets $G_i$ with $1\leq i\leq l$.
Assume that the following conditions hold:
\begin{description}
  \item[(i)] $\dom G$ is convex;
  \item[(ii)] 
  $\gph G|_{Q_i}$  is convex for all $1\leq i\leq l$, where $Q_i$ is the image set of $G_i$ under the projection mapping $(x, y)\mapsto x$.
\end{description}
Then $G$ is Lipschitz continuous
not only on each  $Q_i$ with some constant $\kappa_i>0$, but also
on its domain with constant $\max_{1\leq i\leq l}\kappa_i$.
\end{lemma}
\begin{proof}
To begin with, we notice from assumption (ii) that $\dom G=\bigcup_{i=1}^l Q_i$ with
  $Q_i$ being  a polyhedral set  and
$$
\gph G|_{Q_i}=\{(x,y)\in \gph G\mid x\in Q_i\}=\left\{(x,y)\in \bigcup_{j=1}^l G_j\mid x\in Q_i\right\}=\bigcup_{j=1}^l \left\{(x,y)\in G_j\mid x\in Q_i\right\}.
$$
That is,  $\gph  G|_{Q_i}$ can be expressed as a union of finitely many polyhedral sets. It then follows from assumption  (ii) again and \cite[Theorem 19.6]{CvxAn} that    $\gph  G|_{Q_i}$ is a polyhedral set. Applying \cite[Example 9.35]{Rock2009VaAn},   $G$ is Lipschitz continuous on each $Q_i$ with some constant $\kappa_i>0$.
That is,
\begin{equation}\label{Thm-Lip1}
G(x')\subset G(x)+\kappa_i\|x'-x\|\mathbb B,\quad\forall x',x\in Q_i.
\end{equation}
Next, we will prove the Lipschitz continuity of $G$ on   $\dom G$. Indeed, for any $x',x\in\dom G$, we denote
$
x_\theta:=x'+\theta(x-x')$, $\theta\in[0,1].$
It follows from the convexity of $\dom G$ (see assumption {\bf (i)}) that we get that $x_\theta\in \dom G$ for all $\theta\in[0,1]$. This, combined with the fact that $\dom G=\bigcup_{i=1}^lQ_i$ and $Q_i$ is closed and convex for each $i\in\{1,\cdots,l\}$, yields that there exist $i_j\in\{1,\cdots,l\}$ and $\theta_j\in[0,1]$ ($j=0,1,\cdots,k$), with 
$
0=\theta_0<\theta_1<\cdots<\theta_{k-1}<\theta_k=1,
$
such that
$
x_\theta\in Q_{i_j}$, for all $\theta\in[\theta_{j-1},\theta_j].
$
Then we deduce from \eqref{Thm-Lip1} that
\begin{equation*}
G(x_{\theta_{j-1}})\subset G(x_{\theta_j})+\kappa_{i_j}\|x_{\theta_{j-1}}-x_{\theta_j}\|\mathbb B
=G(x_{\theta_j})+\kappa_{i_j}(\theta_j-\theta_{j-1})\|x'-x\|\mathbb B,\quad j\in\{1,\cdots,k\}.
\end{equation*}
Thus,
\begin{equation*}
\begin{array}{rlc}
G(x')=G(x_{\theta_0})\subset&G(x_{\theta_k})+\sum\limits_{j=1}^k\kappa_{i_j}(\theta_j-\theta_{j-1})\|x'-x\|\mathbb B\\
\subset&G(x)+\max\limits_{j\in\{1,\cdots,k\}}\kappa_{i_j}\|x'-x\|\mathbb B\\
\subset&G(x)+\max\limits_{1\leq i\leq l}\kappa_{i}\|x'-x\|\mathbb B.
\end{array}
\end{equation*}
This completes the proof.
\end{proof}

\begin{remark}
In \cite{Mangasarian1987}, Mangasarian and Shiau studied the Lipschitz continuity of solution multifunctions to system of linear inequalities, linear program and linear complementarity problem. The solution  multifunctions of these problems have a polyhedral multifunction graph. Lemma \ref{lem-Lip-con-poly-theo-ori} can be applied to study the Lipschitz continuity of these multifunctions. However,  the details are not presented in here.

\end{remark}

Next, we give an example to show that in Lemma \ref{lem-Lip-con-poly-theo-ori}, when $G$ is Lipschitz continuous, both assumptions (i) and (ii) are not necessarily held.
\begin{example}\label{exa1}
Let $\Omega=\R\times\{1\}$, $\Omega_2=\R\times\{-1\}$, $\Omega_3=\{0\}\times[-1,1]$ be polyhedral sets and $\Omega=\Omega\cup\Omega_2\cup\Omega_3$, and define $G:\Omega\rightrightarrows\R$ as $G(x,y)=\{0,1\}$ for all $(x,y)\in\Omega$. Clearly, $G(x,y)$ is a constant, and thus $G$ is a polyhedral multifunction and Lipschitz continuous on $\Omega$. However, $\dom G = \Omega$ is not convex and, for any $i=1,2,3$, $\gph  G|_{\Omega_i}$ is the union of two parallel lines/line segments, which is not convex.
\end{example}

While both assumptions (i) and (ii) in Lemma \ref{lem-Lip-con-poly-theo-ori} are not necessary, we next provide two examples to show that the Lipschitz continuity may not be valid whenever either assumption (i) or (ii) in Lemma \ref{lem-Lip-con-poly-theo-ori}
does not hold.
\begin{example}
Take $\Omega$ as in Example \ref{exa1} 
and define $G:\Omega\to\R$ as $G(x,y)=xy$, i.e., $G(x,y)=x$ if $(x,y)\in\Omega$, $G(x,y)=-x$ if $(x,y)\in\Omega_2$ and $G(x,y)=0$ if $(x,y)\in\Omega_3$. Then $\gph  G|_{\Omega}=\{(x,1,x)~|~x\in\R\}$, $\gph  G|_{\Omega_2}=\{(x,-1,-x)~|~x\in\R\}$ and $\gph  G|_{\Omega_3}=\{(0,y,0)~|~y\in[-1,1]\}$. Therefore, all $\gph  G|_{\Omega_i}$, $i=1,2,3$, are polyhedral and $G$ is thus a polyhedral multifunction. However, $\dom G=\Omega$ is not convex. Thus assumption (i) does not hold. And $G$ is not Lipschitz continuous on $\Omega$. Indeed, taking $(x_1^{\nu},y_1^{\nu})=(\nu,1)$ and $(x_2^{\nu},y_2^{\nu})=(\nu,-1)$ for all $\nu\in\mathbb N$, we have
$$
|G(x_1^\nu,y_1^{\nu})-G(x_2^{\nu},y_2^{\nu})|=2\nu = \nu\|(x_1^{\nu},y_1^{\nu})-(x_2^{\nu},y_2^{\nu})\|,~~\forall \nu\in\mathbb N.
$$
\end{example}

\begin{example} This example is borrowed from \cite{Mangasarian1987}. 
We consider the linear complementarity problem
$$
Mx+q\geq0,\quad x\geq0,\quad x^\top(Mx+q)=0,
$$
with
$
M=\left[~
\begin{aligned}
0~~1\\
-1~~0
\end{aligned}
~\right].
$
For any $q\in\mathbb R^2$, the solution set of the linear complementarity problem is denoted by $G(q)$. We find that $\dom G=\R\times\R_+$, which is convex. Let
$$
\begin{aligned}
&Q_1:=(-\infty,0]\times[0,+\infty),~Q_2:=\{0\}\times[0,+\infty),\\
&Q_3:=(0,+\infty)\times[0,+\infty),~Q_4:=\R\times\{0\}.
\end{aligned}
$$
Then $\dom G = Q_1 \cup Q_2 \cup Q_3 \cup Q_4$ and
$$
G(q)=\left\{
\begin{aligned}
&(q_2,-q_1)^\top~~~&&\mbox{if}~q\in \ri Q_1,\\
&([0,q_2]\times\{0\})^\top~~~&&\mbox{if}~q\in \ri Q_2,\\
&(0,0)^\top~~~&&\mbox{if}~q\in \ri Q_3,\\
&(0,\max\{-q_1,0\})^\top~~~&&\mbox{if}~q\in Q_4.
\end{aligned}
\right.
$$
It is easy to show that $G$ is a polyhedral
multifunction, which is not Lipschitz continuous on $\dom G$. Indeed $\gph  G|_{Q_1}$ and $\gph  G|_{Q_3}$ are not convex. Thus assumption (ii) does not hold.
\end{example}

\subsection{Technical lemmas on index partitions}

In this subsection, we present some preliminary lemmas to be used in Sections 3 and 4. In particular, Lemmas 2.2 - 2.5 will be used in Section 3 and all the lemmas in this subsection will be used in Section 4 respectively. 

Let $I:= \{1,\dots, n\}$  and the collection of all index partitions $(I^+, I^0, I^-)$ of $I$ be denoted by  $\mathcal{I}$, where $(I^+, I^0, I^-)$ being a partition of $I$ means that
 $I^+$, $I^0$ and $I^-$  are pairwise  disjoint and their union is $I$.  
Clearly, the index partition   $(I^+(x), I^0(x), I^-(x))$   defined by \eqref{def-index-partition-signs-x} belongs to $\mathcal{I}$ for all $x\in \R^n$, and  
\begin{equation}\label{norm-1-subdiff}
\partial \|\cdot\|_1(x)=\{v\in \R^n\mid v_i=1\;\forall i\in I^+(x),\;v_i\in [-1,1]\;\forall i\in I^0(x),\;v_i=-1\;\forall i\in I^-(x)\}. 
\end{equation}
For every $(I^+, I^0, I^-)\in \mathcal{I}$, define
\begin{equation}\label{w-partition}
    W(I^+,I^0,I^-):=\{x\in \R^n\mid x_i>0\;\forall i\in I^+, x_i=0\;\forall i\in I^0, \;x_i<0\;\forall i\in I^-\},
\end{equation}
which is a nonempty and relatively open subset of $\R^n$.
Note that the collection of
$\{W(I^+,I^0, I^-)\mid (I^+,I^0,I^-)\in \mathcal{I}\}$ forms a partition of $\R^n$.

Define a  multifunction  $\mathbb{Y}: \R^{m\times n}\rightrightarrows \R^m$ by 
\[
A\mapsto \{y\in \R^m\mid \|A^\top y\|_\infty\leq 1\}. 
\]
It is clear that for every $A\in \R^{m\times n}$, $\mathbb{Y}(A)$ is  a nonempty   polyhedron containing $0$ in its interior, which  represents the feasible set of the dual problem \eqref{lasso-dual0001} of the extended $\ell_1$ regularization problem (\ref{ext-lasso-problem-lambda}). We denote the finite  collection of all nonempty faces of $\mathbb{Y}(A)$ by $\mathcal{F}(A)$. 
For every  $(I^+,I^0, I^-)\in \mathcal{I}$, define a multifunction $\mathbb{F}(\cdot, (I^+,I^0, I^-)): \R^{m\times n}\rightrightarrows \R^m$  by 
\begin{equation}\label{def-FFF}
    A\mapsto \{y\in \R^m\mid   A_i^\top y=1\;\forall i\in I^+,\; |A_i^\top y|\leq 1\;\forall i\in I^0,\;  A_i^\top y=-1\;\forall i\in I^-\}.
\end{equation}
Clearly, each $\mathbb{F}(A, (I^+,I^0, I^-))$ is a (possibly empty) face of $\mathbb{Y}(A)$.
For every nonempty subset $F$ of $\mathbb{Y}(A)$,  we  define the following index sets: 
\begin{equation}\label{active-index-sets}
\left\{
    \begin{array}{lll}
\mathcal{A}^+(F, A)&:=& \{i\in I\mid A_i^\top y=1\quad \forall y\in F\}, \\
\mathcal{A}^0(F, A)&:=&  \{i\in I\mid  |A_i^\top y| <1\quad \exists y \in F\},\\
\mathcal{A}^-(F, A)&:=& \{i\in I\mid  A_i^\top y=-1\quad \forall y\in F\}.
    \end{array}
    \right.
\end{equation}
It is evident that the index partition  $(\mathcal{A}^+(F, A), \mathcal{A}^0(F, A),\mathcal{A}^-(F, A))$ belongs to  $\mathcal{I}$ for all nonempty and convex subset $F$ of $\mathbb{Y}(A)$.

Define a multifunction $\mathcal{J}: \R^{m\times n}\rightrightarrows \mathcal{I}$ by  
 \begin{equation}\label{def-partition-faces}
   A\mapsto \left\{\left(I^+, I^0, I^-\right)\in \mathcal{I}\;\left|\;\exists y\in \R^m:  A_i^\top y=1\,\forall i\in I^+,\, |A_i^\top y|< 1\,\forall i\in I^0,\,  A_i^\top y=-1\,\forall i\in I^-\right.
\right\}.
\end{equation}
 Clearly, for every $A\in \R^{m\times n}$,  the index partition $(\emptyset, I, \emptyset)\in \mathcal{J}(A)$ due to $|A_i^\top 0|=0<1$ for all $i\in I$, and  
 each $\mathbb{F}(A, (I^+,I^0,I^-))$ with $(I^+,I^0,I^-)\in \mathcal{J}(A)$ is by definition  a nonempty face of $\mathbb{Y}(A)$. 
Some basic properties of the nonempty faces $F$ in $\mathcal{F}(A)$ and the index partitions $(I^+, I^-,I^-)$ in $\mathcal{J}(A)$ are summarized as follows. 
\begin{lemma}\label{lem-xxys}
Let $A\in \R^{m\times n}$.
The following   hold:
\begin{description}
  \item[(i)] For every  $F\in \mathcal{F}(A)$, we have $(\mathcal{A}^+(F,A), \mathcal{A}^0(F,A),\mathcal{A}^-(F,A)) \in \mathcal{J}(A)$ and 
   \[
(\mathcal{A}^+(F,A), \mathcal{A}^0(F,A),\mathcal{A}^-(F,A))=(\mathcal{A}^+(\bar{y}, A), \mathcal{A}^0(\bar{y},A),\mathcal{A}^-(\bar{y},A)) \quad \forall \bar{y}\in \ri F.
\]
\item[(ii)] For every   $(I^+, I^0,I^-)\in \mathcal{J}(A)$, we have $\mathbb{F}(A, (I^+, I^0,I^-))\in \mathcal{F}(A)$
and 
     \begin{equation}\label{xiangduineibu}
           \ri \mathbb{F}(A, (I^+, I^0,I^-))=\{y\in \R^m\mid A_i^\top y=1\forall i\in I^+,\;  |A_i^\top y|< 1 \forall i\in I^0, \;A_i^\top y=-1 \forall i\in I^-\}.
     \end{equation}
\end{description}
\end{lemma}
\begin{proof} To show (i), let $F\in \mathcal{F}(A)$ and $\bar{y}\in \ri F$. By definition we have $\mathcal{A}^+(F,A)\subset \mathcal{A}^+(\bar{y},A)$. We claim that $\mathcal{A}^+(\bar{y},A)\subset \mathcal{A}^+(F,A)$. Suppose by contradiction there is some $i_0\in \mathcal{A}^+(\bar{y},A)$ but $i_0\not\in \mathcal{A}^+(F,A)$, i.e., $A_{i_0}^\top y<1$ for some $y\in F$. As $F$ is convex and $\bar{y}\in \ri F$, we have $[\bar{y}, y)\subset \ri F$. Thus, there is some $\tilde{y}$ close enough to $y$ in $[\bar{y}, y)$ such that $A_{i_0}^\top \tilde{y}<1$.
Applying \cite[Proposition 7]{ioffe2018variational} to the system of linear inequalities $-1\leq A_i^\top y\leq 1$ for all $i\in I$,  we have $(\mathcal{A}^+(\bar{y},A), \mathcal{A}^0(\bar{y},A),\mathcal{A}^-(\bar{y},A))=(\mathcal{A}^+(\tilde{y},A), \mathcal{A}^0(\tilde{y},A),\mathcal{A}^-(\tilde{y},A))$. This implies that $i_0\in \mathcal{A}^+(\tilde{y},A)$ or equivalently $A_{i_0}^\top \tilde{y}=1$, a contradiction. Therefore, we have $\mathcal{A}^+(F,A)=\mathcal{A}^+(\bar{y},A)$. Similarly, we have $\mathcal{A}^-(F,A)=\mathcal{A}^-(\bar{y},A)$. As a consequence, we get   $\mathcal{A}^0(F,A)=\mathcal{A}^0(\bar{y},A)$. This verifies the equality in  (i).  By definition we have $(\mathcal{A}^+(y, A), \mathcal{A}^0(y,A),\mathcal{A}^-(y,A))\in \mathcal{J}(A)$ for any $y\in \mathbb{Y}(A)$. Thus, we get (i).  

To show (ii), let $(I^+, I^-,I^-)\in \mathcal{J}(A)$. Note that by definition the right-hand-side of \eqref{xiangduineibu} is nonempty. Then  $\mathbb{F}(A, (I^+, I^0,I^-))$ is also nonempty, and we get the equality   \eqref{xiangduineibu} from  
   \cite[Theorem 6.5]{CvxAn}  or  \cite[Corollary 2.1]{zheng2021}. As $\mathbb{F}(A, (I^+, I^0,I^-))$ is   a face of $\mathbb{Y}(A)$, we have $\mathbb{F}(A, (I^+, I^0,I^-))\in \mathcal{F}(A)$.    This completes the proof. \end{proof}

As a natural consequence of the above lemma,  there is a one-to-one correspondence between index  partitions in  $\mathcal{J}(A)$ and nonempty faces in $\mathcal{F}(A)$, entailing that $|\mathcal{J}(A)|=|\mathcal{F}(A)|$ (i.e., the cardinality of $\mathcal{J}(A)$ is equal to that of $\mathcal{F}(A)$). We present such a result in the following lemma without giving a proof as it is rather straightforward.   

\begin{lemma}\label{lem-one-to-one}
     Let $A\in \R^{m\times n}$. 
There is a one-to-one correspondence between the index partitions $(I^+, I^0, I^-)$ in the collection $\mathcal{J}(A)$ and the nonempty faces $F$ in $\mathcal{F}(A)$:
\[
(I^+, I^0, I^-)\quad \longleftrightarrow \quad F
\]
with $F=\mathbb{F}(A, (I^+, I^0, I^-))$ and $(I^+, I^0, I^-)=(\mathcal{A}^+(F, A), \mathcal{A}^0(F, A),\mathcal{A}^-(F, A))$. 
\end{lemma}

Define the projection mappings $P_{1,2}:\R\times\R^m\times  \R^n\to \R\times \R^m$ and $P_3:\R\times \R^m\times  \R^n\to \R^n$ by  $(\lambda,b, x)\mapsto (\lambda, b)$ and  $(\lambda,b, x)\mapsto x$, respectively.  
For every $(I^+,I^0, I^-)\in \mathcal{I}$, define a multifunction $\mathbb{G}(\cdot, (I^+,I^0, I^-)): \R^{m\times n}\rightrightarrows \R_+\times \R^m\times \R^n$   by
\begin{equation}\label{def-sf-new}
  A\mapsto   
  \begin{array}{l}
   \left\{
    (\lambda, b, x) \in \R_+ \times \R^m \times \R^n
    \left|
    \begin{array}{lcl}
    x_i\geq 0, &\hskip-0.2cm A_i^\top (b-Ax)=\lambda &\hskip-0.2cm \forall i\in I^+\\
     x_i=0, &\hskip-0.2cm  | A_i^\top (b-Ax)|\leq \lambda&\hskip-0.2cm \forall i\in I^0\\
     x_i\leq 0, &\hskip-0.2cm A_i^\top (b-Ax)=-\lambda &\hskip-0.2cm \forall i\in I^-
     \end{array}
     \right.
     \right\}.
     \end{array}
\end{equation}
For all $A\in \R^{m\times n}$ and $(I^+, I^0, I^-)\in \mathcal{I}$, it is evident that    $\mathbb{G}(A, (I^+,I^0, I^-))$,  $P_{1,2}(\mathbb{G}( A, (I^+,I^0, I^-)))$ and $P_{3}(\mathbb{G}( A, (I^+,I^0, I^-)))$ are all nonempty and polyhedral cones. 

For further developments in the sequel,   we list out in the following two lemmas a few  basic properties of polyhedral cones $\mathbb{G}(A, (I^+,I^0, I^-))$, $P_{1,2}(\mathbb{G}( A, (I^+,I^0, I^-)))$ and $P_{3}(\mathbb{G}( A, (I^+,I^0, I^-)))$, where $A\in \R^{m\times n}$  and $(I^+,I^0,I^-)\in \mathcal{J}(A)$.
\begin{lemma}\label{lem-basic-proSF}
Let $A\in \R^{m\times n}$ and  $(I^+,I^0,I^-)\in \mathcal{J}(A)$.  In terms of $F:=\mathbb{F}(A,(I^+,I^0,I^-))$, we have 
   \begin{equation}\label{clD-primal-dual}
  \mathbb{G}(A, (I^+,I^0, I^-))= \{(\lambda,Ax+\lambda y+v, x)\mid\lambda\geq 0,\;x\in \cl W(I^+,I^0,I^-),y\in F\cap \rge A,\;v\in \ker A^\top       \},
      \end{equation}
\begin{equation}\label{ri-S_F}
\begin{array}{ll}
 \ri \mathbb{G}(A, (I^+,I^0, I^-))&=\left\{
    (\lambda, b, x)  
    \left|
    \begin{array}{lcl}
     \lambda> 0&&\\
    x_i> 0, &\hskip-0.2cm A_i^\top (b-Ax)=\lambda &\hskip-0.2cm \forall i\in I^+\\
     x_i=0, &\hskip-0.2cm  |A_i^\top (b-Ax)|< \lambda&\hskip-0.2cm \forall i\in I^0\\
     x_i< 0, &\hskip-0.2cm A_i^\top (b-Ax)=-\lambda &\hskip-0.2cm \forall i\in I^-
     \end{array}
     \right.
     \right\}\\[0.5cm]
     &=\{ (\lambda,Ax+\lambda y, x) \mid  \lambda> 0,x\in W(I^+,I^0,I^-),\;y\in \ri F\},
\end{array}
\end{equation}
and
\begin{equation}\label{def-sf-aff}
\aff \mathbb{G}(A, (I^+,I^0, I^-))=\left\{
(\lambda, b, x)
\left|
\begin{array}{rl}
 A_i^\top (b-Ax)=\lambda & \forall i\in I^+\\
 x_i=0 & \forall i\in  I^0\\
  A_i^\top (b-Ax)=-\lambda &\forall i\in I^-
 \end{array}
 \right.
 \right\}.
\end{equation}
Moreover, for any  distinct index partitions $(I^+,I^0, I^-)$ and $(\bar I^+, \bar I^0, \bar I^-)$ in $\mathcal J (A)$,  we have  
\begin{equation}\label{huashui}
    \ri   \mathbb{G}(A, (I^+,I^0, I^-)) \cap  \mathbb{G}(A, (\bar I^+, \bar I^0, \bar I^-)) = \emptyset.
\end{equation}
\end{lemma}
\begin{proof} 
To show \eqref{clD-primal-dual}, denote its right-hand-side  by $G^*$.  We first show $G^*\subset  \mathbb{G}(A, (I^+,I^0, I^-))$. Let $(\lambda, b, x)\in G^*$. Then by definition  $\lambda\geq 0$, $x\in \cl W(I^+,I^0,I^-)$ and $b=Ax+\lambda y+v$ for some $y\in F\cap \rge A$ and $v\in \ker A^\top$. Then given that $A_i^\top (b-Ax)=\lambda A_i^\top y$ for all $i$, we get from  the inclusion $y\in \mathbb{F}(A,(I^+,I^0,I^-))$ that 
\[
     A_i^\top (b-Ax) =\lambda \;\forall i\in I^+, \;\;
      |A_i^\top (b-Ax)| \leq \lambda \;\forall i\in I^0, \;\;
    A_i^\top (b-Ax)  =-\lambda \;\forall i\in I^-.
\]
This gives $(\lambda, b, x)\in\mathbb{G}(A, (I^+,I^0, I^-))$ by definition. So we get  $G^*\subset  \mathbb{G}(A, (I^+,I^0, I^-))$. 

Next we show $\mathbb{G}(A, (I^+,I^0, I^-))\subset G^*$. Let $(\lambda, b, x)\in\mathbb{G}(A, (I^+,I^0, I^-))$.  Then by definition  $\lambda\geq 0$, $x\in \cl W(I^+,I^0,I^-)$ and 
\begin{equation}\label{zhilaohu}
A_i^\top (b-Ax)=\lambda\; \forall i\in I^+, \quad 
       | A_i^\top (b-Ax)|\leq \lambda \; \forall i\in I^0,\quad 
    A_i^\top (b-Ax) =-\lambda \; \forall i\in I^-.
\end{equation}
If $\lambda=0$, we get from (\ref{zhilaohu}) that $A^\top(b-Ax)=0$. Then there is some $v\in \ker A^\top$ such that  $b=Ax+0y+v$ for any  $y\in F\cap\rge A$. So in this case, we get $(\lambda, b, x) \in G^*$ by definition.  Alternatively if $\lambda>0$, in terms of 
$y:=(\proj_{\rge A}b-Ax)/\lambda\in \rge A,$
we get $y\in F$ from (\ref{zhilaohu}) by noting that (\ref{zhilaohu}) still  holds with $\proj_{\rge A}b$ in place of $b$. Let $v:=\proj_{\ker A^\top}(b)$. Then we have 
$b=\proj_{\rge A}b+\proj_{\ker A^\top}(b)=Ax+\lambda y+v$ with $y\in F\cap \rge A$ and $v\in \ker A^\top$. So in this case, we also get  $(\lambda, b, x) \in G^*$ by definition. To sum up, we get $\mathbb{G}(A, (I^+,I^0, I^-))=G^*$, i.e., \eqref{clD-primal-dual} holds. 

By the formula \eqref{xiangduineibu} for  $\ri\mathbb{F}(A,(I^+,I^0, I^-))$, we note that the right-hand-side of  the first equality in   (\ref{ri-S_F})  is  nonempty as it contains any triple $(\lambda, Ax+\lambda y, x)$ with  $\lambda>0$, $x\in W(I^+, I^0, I^-)$  and $y\in \ri \mathbb{F}(A,(I^+,I^0, I^-))$.  Then we get  from  \cite[Theorem 6.5]{CvxAn}  the first equality in   (\ref{ri-S_F}), which, together with  \eqref{xiangduineibu}, implies the second equality in (\ref{ri-S_F}). The  equality  (\ref{def-sf-aff}) follows readily from  the first equality in   (\ref{ri-S_F})  and  the definition of affine hull.

It remains to  show \eqref{huashui} for any  distinct index partitions $(I^+,I^0, I^-)$ and $(\bar I^+, \bar I^0, \bar I^-)$ in $\mathcal J (A)$. Suppose by contradiction that $\ri \mathbb{G}(A, (I^+,I^0, I^-)) \cap \mathbb{G}(A, (\bar I^+, \bar I^0, \bar I^-))$  contains  some $(\lambda,b,x)$.  We claim that $I^+=\bar I^+$. Indeed, for any $i\in I^+$, we deduce from \eqref{ri-S_F} and  the inclusion  $(\lambda,b,x)\in \ri \mathbb{G}(A, (I^+,I^0, I^-))$ that $x_i>0$. Then by the inclusion  $(\lambda,b,x)\in \mathbb{G}(A, (\bar I^+, \bar I^0, \bar I^-))$, we get $i\in \bar I^+$ and hence $I^+\subset \bar I^+$.  To prove the opposite inclusion, let $i\in \bar I^+$. It follows from   the inclusion  $(\lambda,b,x)\in \mathbb{G}(A, (\bar I^+, \bar I^0, \bar I^-))$ that $A_i^\top (b-Ax)=\lambda$.  As  $(\lambda,b,x)\in \ri \mathbb{G}(A, (I^+,I^0, I^-))$, we get  from   \eqref{ri-S_F} that    $i\in I^+$ and hence that $\bar I^+\subset  I^+$. Therefore, we get $I^+= \bar I^+$.  Similarly, we can get  $I^-= \bar I^-$. As a result, we get  $(I^+,I^0, I^-)=(\bar I^+, \bar I^0, \bar I^-)$, a contradiction. This completes the proof.  \end{proof}

\begin{lemma}\label{lem-basic-pro}
Let $A\in \R^{m\times n}$ and $(I^+,I^0,I^-)\in \mathcal{J}(A)$.  In terms of $F:=\mathbb{F}(A,(I^+,I^0,I^-))$,   the following  hold:
\begin{description}
         \item[(i)]   
        $P_{1,2}(\mathbb{G}( A, (I^+,I^0, I^-)))=\{(\lambda,Ax+\lambda y+v)
   \mid  \lambda\geq 0,\;x\in \cl W(I^+,I^0, I^-),y\in F\cap \rge A,\;v\in \ker A^\top
       \}$   is a polyhedral cone  having nonempty interior with
  \[
  \begin{array}{ll}
      \inte P_{1,2}(\mathbb{G}( A, (I^+,I^0, I^-)))&=P_{1,2}(\ri \mathbb{G}(A, (I^+,I^0, I^-)))\\[0.2cm]
      &=\{(\lambda, Ax+\lambda y)\mid \lambda>0,\;x\in W(I^+,I^0, I^-),\;y\in \ri F\}.
  \end{array}        
  \]
 \item[(ii)]  $P_{1,2}(\mathbb{G}( A, (I^+,I^0, I^-)))$ can be decomposed   as follows:
 \[
 P_{1,2}(\mathbb{G}( A, (I^+,I^0, I^-)))=C_1+C_2+\{0\}\times \ker A^\top,
 \]
 where 
 $$
C_1:=\{(0, Ax)\mid x\in \cl W (I^+,I^0, I^-)\}=\cone\left(\left\{ (0, A_i  ) \mid i\in I^+\right\}\cup \{ (0, -A_j )  \mid j\in I^-  \} \right)
 $$
 and 
 $$
C_2:=\{(\lambda, \lambda y)\mid \lambda\geq 0, y\in F\cap \rge A\}=\cone\left\{ (1, V  ) \mid V\in  {\rm extr} (F\cap \rge A) \right\}
 $$
 are polyhedral cones in $\R\times \rge A$. 
\item[(iii)] The $\ell_1$-norm function  $\|\cdot\|_1$ is linear on $P_{3}(\mathbb{G}( A, (I^+,I^0, I^-)))$  with
\[
P_{3}(\mathbb{G}( A, (I^+,I^0, I^-)))=\cl W (I^+,I^0, I^-)\;\mbox{and} \;\ri P_{3}(\mathbb{G}( A, (I^+,I^0, I^-)))= W (I^+,I^0, I^-).
\]
          \end{description}
\end{lemma}
 \begin{proof}   In view of (\ref{clD-primal-dual}), we get the   formula for  $P_{1,2}(\mathbb{G}( A, (I^+,I^0, I^-)))$ in (i) immediately.   From this formula, one can  see that 
 $P_{1,2}(\mathbb{G}( A, (I^+,I^0, I^-)))=C_1+C_2+\{0\}\times \ker A^\top$. It is clear that 
 $C_1$ and $C_2$ are two polyhedral cones in $\R\times \rge A$.    
As $C_1$ and $C_2$ can be generated from   vectors of the forms  $(0, A_i ) $ and $(0, -A_j )$ with $i\in I^+$ and  $j\in I^-$, and  vectors of the form $(1, V) $ with $V \in {\rm extr} (F\cap \rge A)$ respectively,    (ii) is verified.  Moreover, $C_1+C_2$ can be generated from  these vectors, which are called the generators of $C_1+C_2$. 

We now show that $P_{1,2}(\mathbb{G}( A, (I^+,I^0, I^-)))$ has a nonempty interior. As $P_{1,2}(\mathbb{G}( A, (I^+,I^0, I^-)))=C_1+C_2+\{0\}\times \ker A^\top$ and  $C_1+C_2$ is included in $\R\times \rge A$ (i.e., the linear orthogonal complementary space of  $\{0\}\times \ker A^\top$), 
to show that $P_{1,2}(\mathbb{G}( A, (I^+,I^0, I^-)))$ has a nonempty interior  is to illustrate that the dimension of $C_1+C_2$ is $1+\dim(\rge A)$, i.e.,    the rank of the generators of $C_1+C_2$ mentioned above is $1+\dim(\rge A)$. 

When   $I^+\cup I^-=\emptyset$, we have  $W(I^+,I^0, I^-)=\{0\}$ by definition, and then $C_1+C_2=\{(\lambda, \lambda y)\mid \lambda\geq 0, y\in \mathbb{Y}(A)\cap \rge A\}$, whose dimension is $1+\dim(\rge A)$ due to   $0\in \inte \mathbb{Y}(A)\cap \ri (\rge A)$.
Now consider the case that  $I^+\cup I^-\not=\emptyset$.  In this case,
it is trivial to verify that
$
{\rm aff} (F\cap \rge A)=\{y\in \rge A\mid A_i^\top y=1\;\forall i\in I^+,\;A_i^\top y=-1\;\forall i\in I^-\}
$, 
whose dimension is  $\dim(\rge A)-\rank(\{A_i ,i\in I^+\cup I^-\})$ (see \cite[Theorem 2.3.3]{Luc2016}). It is obvious that $0\not \in \aff (F\cap \rge A)$. As $ F\cap \rge A$  is the convex hull of ${\rm extr} (F\cap \rge A)$,  the rank of ${\rm extr} (F\cap \rge A)$ is thus  equal to the dimension of $\aff (F\cap \rge A)$ plus 1, i.e., 
\begin{equation}\label{dim}
\rank ({\rm extr} (F\cap \rge A))=1+\dim(\rge A)-\rank(\{A_i ,i\in  I^+\cup I^-\}).
\end{equation}
Noting from the definition of $F$ that $A_i^\top V=1$ and $A_j^\top V=-1$ for all  $i\in  I^+$ and  $j\in I^-$ and $V\in {\rm extr}(F\cap \rge A)$,  the rank of the generators of $C_1+C_2$ mentioned above   is thus  equal to 
$
\rank(\{A_i ,i\in   I^+\cup I^-\})+\rank ({\rm extr} (F\cap \rge A))
$, 
which is $1+\dim(\rge A)$ by \eqref{dim}. 
 Therefore,  $P_{1,2}(\mathbb{G}( A, (I^+,I^0, I^-)))$ must have  a nonempty interior.  According to \cite[Theorem 6.6]{CvxAn},  we have $\inte P_{1,2}(\mathbb{G}( A, (I^+,I^0, I^-)))=P_{1,2}(\ri \mathbb{G}(A, (I^+,I^0, I^-)))$ and then   from the second equality in (\ref{ri-S_F})  we get the second equality for $\inte P_{1,2}(\mathbb{G}( A, (I^+,I^0, I^-)))$. All the results in (i) are verified. 
 

  To show (iii), we note that the formula for $P_3(\mathbb{G}( A, (I^+,I^0, I^-)))$  follows directly from \eqref{clD-primal-dual} and that  the formula for $\ri P_3(\mathbb{G}( A, (I^+,I^0, I^-)))$ is trivial.  By the formula for $P_3(\mathbb{G}( A, (I^+,I^0, I^-)))$, $\|\cdot\|_1$ is clearly  linear on  $P_3(\mathbb{G}( A, (I^+,I^0, I^-)))$.  This completes the proof. \end{proof}

In the remainder of the section, we introduce two  multifunctions which take values in $\mathcal{I}$ and have close relations with the multifunctions $\mathcal{J}$ and $P_{1,2}(\mathbb{G}(\cdot,(I^+,I^0,I^-)))$. These two multifunctions will be used when we study the Lipschitz continuity of $S$ as  $A$ varies. The first one is    $\mathcal{J}^0:\R^{m\times n}\rightrightarrows \mathcal{I}$ defined  by 
\begin{equation}\label{def-J0}
A\mapsto \left\{(I^+,I^0,I^-)\in \mathcal{J}(A)\;\mid\; \mbox{either}\;I^+\cup I^-=\emptyset\;\mbox{or}\;\rank A_{I^+\cup I^-}=|I^+\cup I^-|\right\},
\end{equation}
and the second one is   $ \mathcal{K}: \R\times \R^m\times \R^{m\times n}\rightrightarrows\mathcal{I}$ defined by 
\begin{equation}\label{def-active}
 (\lambda, b, A)\mapsto\left\{(I^+,I^0,I^-)\in \mathcal{J}(A)\;\mid\;  (\lambda,  b)\in P_{1,2}(\mathbb{G}(A,(I^+,I^0,I^-)))\right\}.
 \end{equation}
 It is evident that $\dom \mathcal{J}^0=\R^{m\times n}$ due to $(\emptyset, I, \emptyset)\in \mathcal{J}^0(A)$ and that $\mathcal{K}(\cdot, \cdot, A)\subset\mathcal{J}(A)$ for all $A\in \R^{m\times n}$.

 The following lemma establishes the inner semicontinuity of $\mathcal{J}^0$. Additionally, it demonstrates a method for generating a new index partition from an existing one within $\mathcal{J}^0(\bar{A})$, where  $\bar{A} \in \mathbb{R}^{m\times n}$. 
   When the matrix $\bar{A}$ has full column rank or full row rank, we also show that, under certain  conditions, both $\mathcal{J}(A)$ and $\mathcal{J}^0(A)$ coincide and remain constant for all $A$ within a quantified neighborhood of $\bar{A}$. 

\begin{lemma}\label{lem-J0A} 
Let  $\bar A\in \R^{m\times n}$. The following hold:
\begin{description}
    \item[(i)] $\mathcal{J}^0$ is inner semicontinuous at $\bar A$, i.e.     $\mathcal{J}^0(\bar A)\subset \mathcal{J}^0(A)$ for all $A$ near $\bar A$. 
    \item[(ii)] For any given $(\bar I^+, \bar I^0, \bar I^-)\in \mathcal{J}^0(\bar A)$,  we have  $(I^+,I^0,I^-)\in\mathcal{J}^0(\bar A)$     
    for all  $(I^+,I^0,I^-)\in \mathcal{I}$ such that $I^+\subset \bar I^+$ and  $I^-\subset\bar I^-$.  
    \item[(iii)] If $\rank(\bar A)=n$, then $\mathcal{J}^0(A)=\mathcal{J}( A)=\mathcal{I}$  for all $A\in \inte\B(\bar A, \sigma_{\min}(\bar A))$.
    \item[(iv)]  If $\rank (\bar A)=m$ and  $\mathcal{J}^0(\bar A)=\mathcal{J}(\bar A)=\mathcal{\bar I}$ for some subcollection $\mathcal{\bar I}$ of  $\mathcal{I}$, then $\mathcal{J}^0(A)=\mathcal{J}(A)=\mathcal{\bar I}$   
    for all $A$ near $\bar A$, which entails that   $\rank(A)=m$ for all $A$ such that $\|A-\bar A\|<\kappa(\bar A)$, where   
    \begin{equation}\label{hmz-banjing}
    \kappa(\bar A):=\sup\{\gamma>0\mid \mathcal{J}^0(A)=\mathcal{J}(A)=\mathcal{\bar I}\;\;\forall A\in \B(\bar A, \gamma)\}>0.
    \end{equation}     
\end{description}
\end{lemma}
\begin{proof} 
To show (i), fix any  $(\bar I^+, \bar I^0, \bar I^-)\in \mathcal{J}^0(\bar A)$.
Given that  $(\emptyset, I, \emptyset)\in \mathcal{J}^0(A)$ for all $A\in \R^{m\times n}$,  we only need to consider the case that $\bar I^+\cup \bar I^-\not=\emptyset$.   Let $\bar y\in \mathbb{F}(\bar A,(\bar I^+,\bar I^0,\bar I^-))$ and  $\bar s:=\bar A_K^\top \bar y$ with $K:=\bar I^+\cup \bar I^-$.  Obviously,  $\bar s$ is independent of a particular choice of $\bar y$. Given that $(\bar I^+, \bar I^0, \bar I^-)\in \mathcal{J}^0(\bar A)$,   there is  a neighborhood $\mathcal{V}$ of $\bar A$ such that $A_K$ has full column rank for all $A\in \mathcal{V}$.   This gives 
\begin{equation}\label{equation-AK}
\{y\in \R^m\mid   A_K^\top y=\bar s\}=\ker   A_K^\top+\{  A_K(  A_K^\top   A_K)^{-1} \bar s\}\quad \forall A\in \mathcal{V}.
\end{equation}
Fix any $y^*\in \ri \mathbb{F}(\bar A,(\bar I^+,\bar I^0,\bar I^-))$. We get $\bar A_K^\top y^*=\bar s$ from \eqref{xiangduineibu}. It follows from \eqref{xiangduineibu} and \eqref{equation-AK} with $A=\bar A$ that there is some $w^*\in \ker {\bar A_K^\top}$ such that
$
y^*=  w^*+\bar A_K(\bar A_K^\top  \bar A_K)^{-1} \bar s
$  and
\begin{equation}\label{niaoj}
 |\bar A_i^\top y^*|<1\quad \forall i\in \bar I^0.
\end{equation}
Then for all $A$ near $\bar A$, we claim that  $(\bar I^+, \bar I^0, \bar I^-)\in \mathcal{J}(A)$   or equivalently that 
\[
  \{y\in \R^m \mid   A_i^\top y=1\;\forall i\in \bar  I^+,\;    |A_i^\top y|<1\;\forall i\in \bar I^0,\;   A_i^\top y=-1\;\forall i\in \bar I^-\}\not=\emptyset.
\]
 Otherwise,  there are some sequences  $A^\nu\to \bar A$ in $\mathcal{V}$  and  
$w^\nu\to w^*$ with $w^\nu\in \ker {(A_K^\nu)^\top}$ for all $\nu$ such that $y^\nu:=w^\nu+A^\nu_K[(A^\nu_K)^\top A^\nu_K]^{-1}\bar s\to y^*$  and $\max_{i\in \bar  I^0}|(A^\nu_i)^\top y^\nu|\geq 1$ for all $\nu$. The existence of such a sequence $\{w^\nu\}$ is possible because  $\bar A_K$ is of full column rank and hence $\ker {(A_K^\nu)^\top}\to \ker {\bar A_K^\top}$, see \cite[Theorem 4.32 (b)]{Rock2009VaAn}. Therefore, $\max_{i\in \bar  I^0} |\bar A_i^\top y^*|\geq 1$, a contradiction to \eqref{niaoj}.  Thus, we  have  $(\bar I^+, \bar I^0, \bar I^-)\in \mathcal{J}(A)$  for all $A$ near $\bar A$. Noting that $A_{K}$ is of full column rank for all $A\in \mathcal{V}$, we further have  $(\bar I^+, \bar I^0, \bar I^-)\in \mathcal{J}^0(A)$ by definition  for all $A$ near $\bar A$. 
This verifies (i)  as $(\bar I^+, \bar I^0, \bar I^-)\in \mathcal{J}^0(\bar A)$ is given arbitrarily and $\mathcal{J}^0(\bar A)$ is a finite collection of index partitions. 

To show (ii), let $(\bar I^+, \bar I^0, \bar I^-)\in \mathcal{J}^0(\bar A)$. It also suffices to consider the case that $\bar I^+\cup \bar I^-\not=\emptyset$, and let $K$,  $y^*$ and $w^*$ be given as above. Let $(I^+,I^0,I^-)\in \mathcal{I}$ be such that $I^+\subset \bar I^+$ and  $I^-\subset\bar I^-$. Since $(I^+,I^0,I^-)\in \mathcal{I}$, we get   
\begin{equation}\label{rengrangxiang000}
I^0=(\bar I^+\setminus I^+)\cup(\bar I^-\setminus I^-)\cup \bar I^0.
\end{equation}
Let $z\in \R^n$ be any vector such that
\begin{equation}\label{xsge}
    z_i=1~\forall i\in I^+,~z_i=1-\alpha~\forall i\in\bar I^+\setminus I^+,~z_i=-1+\alpha~\forall i\in\bar I^-\setminus I^-,~z_i=-1~\forall i\in I^-, 
\end{equation}
where $\alpha\in (0, 1)$.  
Let  $y:=w^*+\bar A_K(\bar A_K^\top \bar A_K)^{-1}z_K$. Then $\bar A_K^\top y=z_K$. So we have 
\begin{equation}\label{rengrangxiang}
    \bar A_i^\top y=1~\forall i\in I^+,~|\bar A_i^\top y|<1~\forall i\in (\bar I^+\setminus I^+)\cup(\bar I^-\setminus I^-),~\bar A_i^\top y=-1~\forall i\in I^-.
\end{equation}
If $\bar I^0=\emptyset$, we get from \eqref{rengrangxiang000} and \eqref{rengrangxiang} that  $(I^+,I^0,I^-)\in\mathcal J(\bar A)$ and hence $(I^+,I^0,I^-)\in\mathcal J^0(\bar A)$ by definition. It remains to consider the case that $\bar I^0\not=\emptyset$, for which, we have 
\begin{equation}\label{xsge000}
|\bar A_i^\top y|\leq  |\bar A_i^\top y^*|+|\bar A_i^\top (y-y^*)|=|\bar A_i^\top y^*|+|\bar A_i^\top \bar A_K(\bar A_K^\top \bar A_K)^{-1}(z_K-\bar s)|\leq \rho+\sqrt n\beta\alpha\quad \forall i\in \bar I^0,
\end{equation}
where $\rho:=\max_{i\in\bar I^0}|\bar A_i^\top y^*|\in [0, 1)$ (due to \eqref{niaoj}) and  $\beta:=\max_{i\in\bar I^0}\|\bar A_i^\top \bar A_K(\bar A_K^\top \bar A_K)^{-1}\|$. By letting  $\alpha$ in \eqref{xsge}  be
$
 (1-\rho)/\sqrt n(1+\beta), 
$
which is clearly  a number in (0, 1), we can still have \eqref{rengrangxiang}, and moreover, we can get from (\ref{xsge000}) that
$
|\bar A_i^\top y|\leq  \rho+\sqrt n\beta\alpha<1\quad \forall i\in \bar I^0.
$
This, together with \eqref{rengrangxiang000}, also  implies that  $(I^+,I^0,I^-)\in\mathcal J(\bar A)$ and hence   $(I^+,I^0,I^-)\in\mathcal J^0(\bar A)$ by definition. 

 To show (iii), fix any $(I^+, I^0, I^-)\in \mathcal{I}$ and any $A\in \R^{m\times n}$ such that $\|A-\bar A\|<\sigma_{\min}(\bar A)$.  
 Given that $\rank(\bar A)=n$ and  that $\sigma_{\min}(\bar A)$ is equal to the minimum of $\|\bar A x\|$ over the unit sphere, we have $\sigma_{\min}(\bar A)>0$. Moreover, the following hold for all $x\in\R^n$ with $\|x\|=1$:   
$$
\|Ax\|=\|\bar A x+(A-\bar A)x\|\geq\|\bar A x\|-\|A-\bar A\|\|x\|\geq \sigma_{\min}(\bar A)-\|A-\bar A\|>0,
$$
which implies that $\rank(A)=n$. Let $s \in \mathbb{R}^n$ be such that $s_i=1$ for all $i\in I^+$, $s_i=0$ for all $i\in I^0$ and $s_i=-1$ for all $i\in I^-$. Since   the linear system $A^\top y = s$ admits a solution $y \in \mathbb{R}^m$ due to $A^\top$ being of full row rank,  we get $(I^+, I^0, I^-)\in\mathcal{J}(A)$  by definition and hence $(I^+, I^0, I^-)\in\mathcal{J}^0(A)$ due to $A$ being of full column rank. Given that $(I^+, I^0, I^-)\in \mathcal{I}$ is arbitrarily fixed and that $\mathcal{J}^0(A)\subset \mathcal{J}(A)\subset \mathcal{I}$ holds by definition, we get $\mathcal{J}^0(A)=\mathcal{J}(A)=\mathcal{I}$.  

To show (iv), we note from (i) that $\mathcal{\bar I}\subset \mathcal{J}^0(A)$ for all $A$ near $\bar A$. Given that the inclusion  $\mathcal{J}^0(A)\subset\mathcal{J}(A)$ holds for all $A$ by definition,  it suffices  to show the inclusion  $\mathcal{J}(A)\subset \mathcal{\bar I}$ for all $A$ near $\bar A$. Suppose to the contrary  that there is some $A^\nu\to \bar A$ such that $\mathcal{J}(A^\nu)\backslash\mathcal{\bar I}\not=\emptyset$ for all $\nu$. 
  Given that $\mathcal{J}(A^\nu)\backslash\mathcal{\bar I}$ is a subset of the finite collection $\mathcal{I}$ for all $\nu$, by taking a subsequence if necessary, we assume that there is some $(I^+, I^0, I^-)\in \mathcal{I}$ such that   $(I^+, I^0, I^-)\in \mathcal{J}(A^\nu)\backslash\mathcal{\bar I}$ for all $\nu$. By the definition of $\mathcal{J}$ (cf. \eqref{def-partition-faces}), for every $\nu$, there is some $y^\nu\in \mathbb{Y}(A^\nu)$ such that 
    \begin{equation}\label{jicao}
   ( A_i^\nu)^\top y^\nu=1\;\forall i\in I^+,\;\;|( A_i^\nu)^\top y^\nu|<1\;\forall i\in I^0,\;\;( A_i^\nu)^\top y^\nu=-1\;\forall i\in I^-.
    \end{equation}
    We claim that $\{y^\nu\}$ is a bounded sequence, for otherwise we would have $\|y^\nu\|\to  \infty$ by taking a subsequence if necessary, which implies  that $(A^\nu)^\top\frac{y^\nu}{\|y^\nu\|}\to 0$ due to $\|(A^\nu)^\top\frac{y^\nu}{\|y^\nu\|}\|_\infty\leq \frac{1}{\|y^\nu\|}$ for all $\nu$ following  the definition of $\mathbb{Y}(A^\nu)$.  By taking a subsequence if necessary, we assume that $\frac{y^\nu}{\|y^\nu\|}\to \bar y$ for some unit vector  $\bar y$.  In light of $A^\nu\to \bar A$, we further get $\bar A^\top \bar y=0$, contradicting to the fact that $\bar A$ is of full row rank. Therefore, $\{y^\nu\}$ is a bounded sequence, and  by taking a subsequence if necessary, we assume that $y^\nu\to y^*$. Then it follows from \eqref{jicao} that 
        \[
   (\bar A_i)^\top y^*=1\;\forall i\in I^+,\;\;|(\bar A_i)^\top y^*|\leq 1\;\forall i\in I^0,\;\;(\bar A_i)^\top y^*=-1\;\forall i\in I^-.
    \]
This implies that  
            \begin{equation}\label{jicao000}
   (\bar A_i)^\top y^*=1\;\forall i\in I^+\cup I^{0+},\;\;|(\bar A_i)^\top y^*|< 1\;\forall i\in I^0\backslash (I^{0+}\cup I^{0-}),\;\;(\bar A_i)^\top y^*=-1\;\forall i\in I^-\cup I^{0-},
    \end{equation}
where $I^{0+}:=\{i\in I^0\mid (\bar A_i)^\top y^*=1\}$ and $I^{0-}:=\{i\in I^0\mid (\bar A_i)^\top y^*=-1\}$. So by the  definition of $\mathcal{J}$ (cf. \eqref{def-partition-faces}),  we get from \eqref{jicao000} that  $(I^+\cup I^{0+}, I^0\backslash (I^{0+}\cup I^{0-}), I^-\cup I^{0-})\in \mathcal{J}(\bar A)$. Given that $\mathcal{J}(\bar A)=\mathcal{J}^0(\bar A)$ by assumption, we have $(I^+\cup I^{0+}, I^0\backslash (I^{0+}\cup I^{0-}), I^-\cup I^{0-})\in \mathcal{J}^0(\bar A)$.  Then we get from (ii) that $(I^+, I^0, I^-)\in \mathcal{J}^0(\bar A)$.  Given that $\mathcal{J}^0(\bar A)=\mathcal{\bar I}$ by assumption, we obtain that  
 $(I^+, I^0, I^-)\in \mathcal{\bar I}$, a contradiction to   our previous assumption that $(I^+, I^0, I^-)\in \mathcal{J}(A^\nu)\backslash\mathcal{\bar I}$ for all $\nu$.  
 This indicates that $\mathcal{J}^0(A)=\mathcal{J}(A)=\mathcal{\bar I}$ for all $A$ near $\bar A$ and hence that 
 \[
 \Gamma:=\{\gamma>0\mid \mathcal{J}^0(A)=\mathcal{J}(A)=\mathcal{\bar I}\;\;\forall A\in \B(\bar A, \gamma)\}\not=\emptyset. 
 \]
This implies that $\kappa(\bar A):=\sup \Gamma>0$. 
 
  As  $\bar A$ is of full row rank and  $\mathbb{Y}(\bar A)$ is by definition the solution set of the system of linear inequalities $\bar A_i^\top y\leq 1$ and $-\bar A_i^\top y\leq 1$ for all $i\in I$,  $\mathbb{Y}(\bar A)$ is clearly  
 a non-empty  bounded polyhedral set  having at least one extreme point   (cf. \cite[Corollary 18.5.1]{CvxAn}).  Let $\tilde y\in \R^m$ be such an extreme point. Noting that $\mathcal{A}^+(\tilde y, \bar A)$ and $\mathcal{A}^-(\tilde y, \bar A)$ are the active index sets of $\tilde y$ with respect to the above system of inequalities   defining  $\mathbb{Y}(\bar A)$,  the rank of the corresponding columns $\bar A_i$ with $i\in \mathcal{A}^+(\tilde y, \bar A)\cup \mathcal{A}^-(\tilde y, \bar A)$ is $m$. Let  $(\tilde I^+, \tilde I^0, \tilde I^-):=(\mathcal{A}^+(\tilde y, \bar A), \mathcal{A}^0(\tilde y, \bar A), \mathcal{A}^-(\tilde y, \bar A))$. By Lemma \ref{lem-xxys} (i), we have $(\tilde I^+, \tilde I^0, \tilde I^-)\in \mathcal{J}(\bar A)$. Due to $\mathcal{J}(\bar A)=\mathcal{J}^0(\bar A)$ by assumption,   we have  $(\tilde I^+, \tilde I^0, \tilde I^-)\in \mathcal{J}^0(\bar A)$ and hence  $\rank (\bar A_{\tilde I^+\cup \tilde I^-})=|\tilde I^+\cup \tilde I^-|$  by definition. This gives  $|\tilde I^+\cup \tilde I^-|=m$.  Let $A\in \R^{m\times n}$ be such that $\|A-\bar A\|<\kappa (\bar A)$. By the definition of $\kappa(\bar A)$, we get $(\tilde I^+, \tilde I^0, \tilde I^-)\in\mathcal{J}^0(A)$ and hence $\rank ( A_{\tilde I^+\cup \tilde I^-})=|\tilde I^+\cup \tilde I^-|=m$.  Therefore, $\rank(A)=m$. The proof is completed.   
\end{proof}

The following lemma establishes the outer semicontinuity of $\mathcal{K}$  under three  conditions.  
 \begin{lemma}\label{lem-active-outer-conti}
  $\mathcal{K}$ is outer semicontinuous at some $(\bar\lambda, \bar b,\bar A)\in \R_+ \times \R^m\times \R^{m\times n}$, i.e., 
\[
\mathcal{K}(\lambda, b, A)\subset\mathcal{K}(\bar\lambda, \bar b,\bar A)\quad\forall (\lambda, b, A)\;\; \mbox{near}\;\; (\bar \lambda, \bar b, \bar A),
\]
 if one of the following  conditions holds:
\begin{description}
    \item[(i)] $\bar\lambda>0$ and  $\mathcal{K}(\bar \lambda, \bar b, \bar A)\subset \mathcal{J}^0(\bar A)$.
    \item[(ii)] $\rank (\bar A)=m$ and  $\mathcal{J}^0(\bar A)=\mathcal{J}(\bar A)$. 
    \item[(iii)] $\rank(\bar A)=n$.
\end{description}
\end{lemma}
\begin{proof}
To show the outer semicontinuity of $\mathcal{K}$ at  $(\bar\lambda, \bar b,\bar A)$, suppose to the contrary  that there  is  some  $(\lambda^\nu,b^\nu,A^\nu)\to(\bar\lambda,\bar b,\bar A)$ such that $\mathcal{K}(\lambda^\nu,b^\nu, A^\nu)\backslash \mathcal{K}(\bar\lambda,\bar b, \bar A)\not=\emptyset$ for all $\nu$.  For every $\nu$, let $(I_\nu^+,I_\nu^0,I_\nu^-)\in \mathcal{K}(\lambda^\nu,b^\nu, A^\nu)\backslash \mathcal{K}(\bar\lambda,\bar b, \bar A)$.   As $\mathcal{K}(\lambda^\nu,b^\nu, A^\nu)\backslash \mathcal{K}(\bar\lambda,\bar b, \bar A)$ is included in the finite collection $\mathcal{I}$ for all $\nu$,  by taking a subsequence if necessary, we can assume that  $(I_\nu^+,I_\nu^0,I_\nu^-)\equiv (I^+,I^0,I^-)$ for all $\nu$.  That is, we have $(I^+,I^0,I^-)\in \mathcal{K}(\lambda^\nu,b^\nu, A^\nu)$ for all $\nu$, but
\begin{equation}\label{zhongjiemubiao}
    (I^+,I^0,I^-)\not\in \mathcal{K}(\bar\lambda,\bar b, \bar A).
\end{equation}

First we assume  condition (i). Given that $\bar\lambda>0$, we can assume that $\lambda^\nu>0$ for all $\nu$. By the definition of $\mathcal{K}$ (cf. \eqref{def-active}), we have $(I^+, I^0, I^-)\in \mathcal{J}(A^\nu)$ and  $(\lambda^\nu, b^\nu)\in P_{1,2}(\mathbb{G}(A^\nu, (I^+, I^0, I^-)))$ for all $\nu$. In view of the definition of $P_{1,2}$, we can associate with  each $\nu$  some $x^\nu\in \R^n$ such that 
$(\lambda^\nu, b^\nu, x^\nu)\in \mathbb{G}(A^\nu, (I^+, I^0, I^-))$. Then for all $\nu$,  we get from the definition of $\mathbb{G}$  (cf. \eqref{def-sf-new}) that  
\begin{equation}\label{sos}
\begin{array}{lcl}
    x^{\nu}_i\geq 0, &\hskip-0.2cm (A_i^\nu)^\top (b^{\nu}-A^\nu x^{\nu})=\lambda^{\nu} &\hskip-0.2cm \forall i\in I^+,\\
     x_i^{\nu}=0, &\hskip-0.2cm  |(A_i^\nu)^\top (b^{\nu}-A^\nu x^{\nu})|\leq  \lambda^{\nu}&\hskip-0.2cm \forall i\in I^0,\\
     x_i^{\nu}\leq  0, &\hskip-0.2cm (A^\nu_i)^\top (b^{\nu}-A^\nu x^{\nu})=-\lambda^{\nu} &\hskip-0.2cm \forall i\in I^-.
     \end{array}
\end{equation}
This, together with the fact that   $\lambda^\nu>0$ and  the formula \eqref{norm-1-subdiff} for the subdifferential set of $\|\cdot\|_1$, implies that  $0\in(A^\nu)^\top(A^\nu x^\nu-b^\nu)+\lambda^\nu\partial\|\cdot\|_1(x^\nu)$ and hence by the optimality condition that $x^\nu$ is an optimal solution to the Lasso problem \eqref{lasso-problem-lambda-b1}. So for all $\nu$, we get   
$$
  \frac{1}{2\lambda^\nu}\|A^\nu x^\nu-b^\nu\|^2+\|x^\nu\|_1\leq  \frac{1}{2\lambda^\nu}\|A^\nu 0-b^\nu\|^2+\|0\|_1=\frac{\|b^\nu\|^2}{2\lambda^\nu}.
$$
This, together with the fact that $(\lambda^\nu, b^\nu)\to (\bar \lambda, \bar b)$ with $\bar \lambda>0$,  implies that $x^\nu$ is a bounded sequence. By taking a subsequence if necessary, we assume that $x^\nu\to \bar x$. By letting $\nu\to \infty$ and taking into account that 
$(\lambda^\nu, b^\nu,A^\nu, x^\nu)\to (\bar \lambda, \bar b,\bar A, \bar x)$, we get from \eqref{sos} that 
  \begin{equation}\label{haya}
      \begin{array}{lcl}
    \bar x_i\geq 0, &\hskip-0.2cm \bar A_i^\top (\bar b-\bar A \bar x)=\bar\lambda &\hskip-0.2cm \forall i\in I^+,\\
     \bar x_i=0, &\hskip-0.2cm  |\bar A_i^\top (\bar b-\bar A \bar x)|\leq \bar\lambda&\hskip-0.2cm \forall i\in I^0,\\
    \bar x_i\leq  0, &\hskip-0.2cm \bar A_i^\top (\bar b-\bar A \bar x)=-\bar\lambda &\hskip-0.2cm \forall i\in I^-.
     \end{array}
  \end{equation}
 Then by the definitions   of $\mathbb{G}$   and  $P_{1,2}$, we have  $(\bar \lambda, \bar b, \bar x)\in \mathbb{G}(\bar A, (I^+,I^0, I^-))$  and  \begin{equation}\label{buqiyan}
     (\bar \lambda, \bar b)\in P_{1,2}(\mathbb{G}(\bar A, (I^+,I^0, I^-))).
 \end{equation}
In terms of $I^{0+}:=\{i\in I^0\mid \bar A_i^\top (\bar b-\bar A \bar x)=\bar\lambda \}$, $I^{0-}:=\{i\in I^0\mid \bar A_i^\top (\bar b-\bar A \bar x)=-\bar\lambda \}$  and  $\bar y:=(\bar b-\bar A\bar x)/\bar \lambda$, we get from \eqref{haya} that 
\[
      \begin{array}{lcl}
    \bar x_i\geq 0, &\hskip-0.2cm \bar A_i^\top (\bar b-\bar A \bar x)=\bar\lambda &\hskip-0.2cm \forall i\in \bar I^+:=I^+\cup I^{0+},\\
     \bar x_i=0, &\hskip-0.2cm  |\bar A_i^\top (\bar b-\bar A \bar x)|< \bar\lambda&\hskip-0.2cm \forall i\in \bar I^0:=I^{0}\backslash (I^{0+}\cup I^{0-}) ,\\
    \bar x_i\leq  0, &\hskip-0.2cm \bar A_i^\top (\bar b-\bar A \bar x)=-\bar\lambda &\hskip-0.2cm \forall i\in \bar I^-:= I^-\cup I^{0-},
     \end{array}
\]
  and hence that 
 \[
 \bar A_i^\top \bar y=1\;\forall i\in \bar I^+,\;\;|\bar A_i^\top \bar y|<1\;\forall i\in  \bar I^{0},\;\;\bar A_i^\top \bar y=-1\;\forall i\in \bar I^-.
 \]
 Then  we get by definition that  $(\bar I^+, \bar I^0, \bar I^-) \in \mathcal{J}(\bar A)$ and  $(\bar\lambda, \bar b, \bar x)\in \mathbb{G}(\bar A,(\bar I^+, \bar I^0, \bar I^-))$, which implies  that $(\bar\lambda, \bar b)\in P_{1,2}(\mathbb{G}(\bar A,(\bar I^+, \bar I^0, \bar I^-)))$ and hence that $(\bar I^+, \bar I^0, \bar I^-)\in \mathcal{K}(\bar\lambda, \bar b,\bar A)$.  It then follows from the inclusion in   condition (i) that  $(\bar I^+, \bar I^0, \bar I^-)\in \mathcal{J}^0(\bar A)$.  Given  that $I^+\subset \bar I^+$ and $I^-\subset \bar I^-$, we get  from Lemma \ref{lem-J0A} (ii) that $(I^+, I^0, I^-)\in \mathcal{J}^0(\bar A)\subset \mathcal{J}(\bar A)$.  Then we get from \eqref{buqiyan} and the definition of $\mathcal{K}$ (cf. \eqref{def-active})  that 
  $(I^+, I^0, I^-)\in  \mathcal{K}(\bar \lambda, \bar b, \bar A)$, a contradiction to \eqref{zhongjiemubiao}. 
  
  Next we assume   condition (ii).  
 Given that $A^\nu\to \bar A$, we apply Lemma \ref{lem-J0A} (iv) to assert that $\mathcal{J}^0(A^\nu)=\mathcal{J}(A^\nu)=\mathcal{J}^0(\bar A)=\mathcal{J}(\bar A)$ for all $\nu$ large enough.  
 Since $(I^+,I^0,I^-)\in \mathcal{K}(\lambda^\nu,b^\nu, A^\nu)$ for all $\nu$, we get by the definition of $\mathcal{K}$ (cf. \eqref{def-active})  that $(I^+,I^0,I^-)\in \mathcal{J}(A^\nu)$ and
 $(\lambda^\nu, b^\nu)\in P_{1,2}(\mathbb{G}(A^\nu, (I^+,I^0,I^-)))$ for all $\nu$. In particular, we have 
 $(I^+,I^0,I^-)\in \mathcal{J}(\bar A)\cap \mathcal{J}^0(A^\nu)$ for all $\nu$ large enough.     
 Following the definition of $P_{1,2}$,  every $(\lambda^\nu, b^\nu)$ corresponds to some $x^\nu\in \R^n$ such that $(\lambda^\nu, b^\nu, x^\nu)\in \mathbb{G}(A^\nu,(I^+,I^0,I^-))$ or explicitly   
 \begin{equation}\label{sos’}
\begin{array}{lcl}
\lambda^\nu\geq 0, &&\\
    x^{\nu}_i\geq 0, &\hskip-0.2cm (A_i^\nu)^\top (b^{\nu}-A^\nu x^{\nu})=\lambda^{\nu} &\hskip-0.2cm \forall i\in I^+,\\
     x_i^{\nu}=0, &\hskip-0.2cm  |(A_i^\nu)^\top (b^{\nu}-A^\nu x^{\nu})|\leq  \lambda^{\nu}&\hskip-0.2cm \forall i\in I^0,\\
     x_i^{\nu}\leq  0, &\hskip-0.2cm (A^\nu_i)^\top (b^{\nu}-A^\nu x^{\nu})=-\lambda^{\nu} &\hskip-0.2cm \forall i\in I^-.
     \end{array}
\end{equation}
This, together with the fact that the column vectors $A^\nu_i$ with $i\in I^+\cup I^-$ are linearly independent  due to $(I^+,I^0,I^-)\in  \mathcal{J}^0(A^\nu)$ for all $\nu$ large enough, implies that  $x^\nu$ is uniquely determined by the linear equations in \eqref{sos’}  and hence   converges   to some $\bar x$.  By letting $\nu\to \infty$, we get from  \eqref{sos’}    that   $(\bar \lambda, \bar b, \bar x)\in \mathbb{G}(\bar A, (I^+,I^0, I^-))$  (cf.\eqref{def-sf-new})  and hence that  $(\bar \lambda, \bar b)\in P_{1,2}(\mathbb{G}(\bar A, (I^+,I^0, I^-)))$.   
 Given that $(I^+,I^0,I^-)\in \mathcal{J}(\bar A)$, we get from  the definition  of $\mathcal{K}$ (cf.\eqref{def-active}) again that  $(I^+, I^0, I^-)\in  \mathcal{K}(\bar \lambda, \bar b, \bar A)$, a contradiction to \eqref{zhongjiemubiao}. 
 
Finally, we invoke   condition \textnormal{(iii)}. By retracing the same steps as in the proof of \textnormal{(ii)}, we derive a contradiction to \eqref{zhongjiemubiao} by considering the following facts: (1) The matrix $A^\nu$ has full column rank for all sufficiently large $\nu$; (2) $\mathcal{J}^0(\bar A)=\mathcal{J}(\bar A) = \mathcal{I}$ (cf. Lemma \ref{lem-J0A} (iii)).  This completes the proof.
\end{proof}

\section{Lipschitz continuity when $A$ is fixed}\label{sec-lip-A-fixed}

 In this  section, for a fixed $A\in \R^{m\times n}$,  we will first apply Lemma \ref{lem-Lip-con-poly-theo-ori} to  establish the Lipschitz continuity of  $S(\cdot,\cdot, A)$ on its domain. As a consequence of this Lipschitz continuity of $S(\cdot,\cdot,A)$, we will obtain an explicit representation of $S(\cdot,\cdot, A)$ and then take advantage of this representation to establish the Lipschitz continuity of $S$ when $A$ varies locally in Section 4.

 Recall that $\dom S(\cdot, \cdot, A) = \R_+ \times \R^m$ is a closed and convex set.
Therefore, to apply Lemma \ref{lem-Lip-con-poly-theo-ori} to the solution multifunction mapping $S(\cdot, \cdot, A)$, it suffices to  verify  that $S(\cdot, \cdot, A)$ is a polyhedral multifunction and that condition (ii) in  Lemma \ref{lem-Lip-con-poly-theo-ori} is satisfied.

First, we verify that $S(\cdot, \cdot, A)$  is a polyhedral multifunction. 
To achieve this, we utilize the index partitions $(I^+, I^0, I^-)$ in $\mathcal{J}(A)$ to divide the graph  $\gph S(\cdot, \cdot, A)$ into  finitely many polyhedral cones $\mathbb{G}(A, (I^+, I^0, I^-))$, where $\mathcal{J}$ and $\mathbb{G}(\cdot, (I^+, I^0, I^-))$ are defined by \eqref{def-partition-faces} and 
\eqref{def-sf-new} respectively.

\begin{proposition}\label{prop-gphS-lambda-b} 
For any $A\in \R^{m\times n}$,    we have 
\begin{eqnarray}
 \gph S(\cdot, \cdot, A) 
 &=&\left\{(\lambda,b, x) \mid \lambda \geq 0, 0\in \rge A^\top + \partial\|\cdot\|_1(x),\;  A^\top (b-Ax)\in  \lambda\partial \|\cdot\|_1(x)\right\},
\label{lam-b-xs'}\\
 &=&\bigcup_{(I^+, I^0, I^-)\in \mathcal{J}(A)}\mathbb{G}(A, (I^+,I^0, I^-)).\label{key-decom}
\end{eqnarray}
\end{proposition}

Recall from Lemma \ref{lem-basic-proSF} that any two polyhedral cones $\mathbb{G}(A, (I^+,I^0, I^-))$ and $\mathbb{G}(A, (\bar I^+, \bar I^0, \bar I^-))$ with $(I^+,I^0, I^-)\not=(\bar I^+, \bar I^0, \bar I^-)$ can intersect only on their relative boundaries, which indicates that  the decomposition of $\gph S(\cdot, \cdot, A)$  by \eqref{key-decom}  can be considered as  a partition of $\gph S(\cdot, \cdot, A)$. 

\begin{proof}   We  divide the proof into  four steps as follows. 

{\bf Step 1.} We show \eqref{lam-b-xs'}. For any $x\in S(\lambda,b,A)$,  if $A^\top ( b-Ax)\in  \lambda\partial \|\cdot\|_1(x)$ holds with $\lambda>0$,  we get $0\in \rge A^\top +\partial\|\cdot\|_1(x)$, and alternatively if $A^\top (b-Ax)\in  \lambda\partial \|\cdot\|_1(x)$ with $\lambda=0$, we get 
$A^\top(Ax-b)=0$. This, together with  (\ref{zyx-lambda}), gives  
$$
S(\lambda,b, A)=\left\{x\in \R^n\mid    0\in \rge A^\top 
 + \partial\|\cdot\|_1(x),\;  A^\top (b-Ax)\in  \lambda\partial \|\cdot\|_1(x)\right\} 
$$
for all $(\lambda, b, A)\in \dom S$ and hence the equality  \eqref{lam-b-xs'}. 

In the following three steps we apply \eqref{lam-b-xs'} to establish \eqref{key-decom}.

{\bf Step 2.} We show that $\gph   S(\cdot, \cdot, A) $ is closed.   To that end, let $(\lambda^\nu,b^\nu,x^\nu)\to(\lambda,b,x)$ with $(\lambda^\nu,b^\nu,x^\nu)\in   \gph   S(\cdot, \cdot, A) $ for all $\nu$. According to  \eqref{lam-b-xs'}, we have
$$
\lambda^\nu\geq 0,\quad 0\in\operatorname{rge}~A^\top + \partial\|\cdot\|_1(x^\nu)~~~\mbox{and}~~~A^\top(b^\nu-Ax^\nu)\in\lambda^\nu\partial\|\cdot\|_1(x^\nu)\quad \forall \nu.
$$
Then we get $\lambda\geq 0$ and  there exist  $\{y^\nu,z^\nu \} \subset \partial\|\cdot\|_1(x^\nu) $  such that
$
y^\nu\in -\operatorname{rge}~A^\top$ and $A^\top(b^\nu-Ax^\nu)=\lambda^\nu z^\nu$, for all $\nu$.
As $\{y^\nu\}$ and $\{z^\nu\}$ are bounded,  by taking some subsequences if necessary, we assume that $y^\nu\to y$ and $z^\nu\to z$. Thus, 
$
y\in -\operatorname{rge}~A^\top$ and $A^\top(b-Ax)=\lambda z.
$
Moreover, it follows from the outer semicontinuity of the subdifferential (cf.\cite[Proposition 8.7]{Rock2009VaAn}) that  $\{y,z\}\subset\partial\|\cdot\|_1(x)$. Therefore,
$$
\lambda\geq 0,\quad 0\in\operatorname{rge}~A^\top+\partial\|\cdot\|_1(x)~~~\mbox{and}~~~A^\top(b-Ax)\in\lambda\partial\|\cdot\|_1(x).
$$
By \eqref{lam-b-xs'} again, we  have  $(\lambda,b,x)\in    \gph   S(\cdot, \cdot, A) $. This verifies the closedness of $ \gph   S(\cdot, \cdot, A) $.

{\bf Step 3.} We prove $\gph   S(\cdot, \cdot, A)\subset \bigcup_{(I^+, I^0, I^-)\in \mathcal{J}(A)}\mathbb{G}(A, (I^+,I^0, I^-))$. Fix any $(\lambda,b,x)\in\gph   S(\cdot, \cdot, A)$ and define 
$$
I^+:=\{i\in I\mid x_i>0\},\quad I^0:=\{i\in I\mid x_i=0\},\quad I^-:=\{i\in I\mid x_i<0\}.
$$
 Then by \eqref{lam-b-xs'}, we have $\lambda\geq 0$ and  $0\in \rge A^\top+\partial\|\cdot\|_1(x)$, the latter of which implies that $F:=\mathbb{F}(A, (I^+, I^0, I^-))\not=\emptyset$ and hence $F\in \mathcal{F}(A)$.  
For each $i\in\mathcal A^+(F,A)$,  we have $A_i^\top y=1$ for all $y\in F$, implying by the definition of $F$ that $i\notin I^-$ and hence that $x_i\geq 0$. Similarly, we get $x_i\leq0$ for all $i\in \mathcal A^-(F, A)$. For each  $i\in\mathcal A^0(F,A)$, 
 we have $|A_i^\top y|<1$ for some $y\in F$, implying by the definition of $F$ that $i\notin I^+\cup I^-$ and hence that $x_i=0$.
All above imply that $x\in \cl W(\bar I^+, \bar I^0, \bar I^-)$, where  
$$(\bar I^+, \bar I^0, \bar I^-):=(\mathcal A^+(F,A),\mathcal A^0(F,A),\mathcal A^-(F,A))\in \mathcal{J}(A).$$

If $\lambda=0$, we get from  \eqref{lam-b-xs'} that $A^\top(b-Ax)=0$. Then taking into account that $x\in \cl W(\bar I^+, \bar I^0, \bar I^-)$,  we get $(\lambda,b,x)\in \mathbb{G}(A, (\bar I^+,\bar I^0, \bar I^-))$. 

Alternatively  if $\lambda>0$,  we get from  \eqref{lam-b-xs'}  that   $ A^\top (b-Ax)/\lambda\in \partial \|\cdot\|_1(x)$.  In view of (\ref{norm-1-subdiff}) and the definition of $F$, we have $(b-Ax)/\lambda\in F$. 
As $F\in \mathcal{F}(A)$, we get from Corollary \ref{lem-one-to-one} that 
$
    F=\mathbb{F}(A, (\bar I^+, \bar I^0, \bar I^-)).
$
This implies that  $(b-Ax)/\lambda \in \mathbb{F}(A, (\bar I^+, \bar I^0, \bar I^-))$ and hence that 
$$
A_i^\top (b-Ax)=\lambda \; \forall i\in \bar I^+,\; |A_i^\top (b-Ax) |\leq  \lambda\; \forall i\in \bar I^0,\; A_i^\top (b-Ax) =-\lambda\;  \forall i\in \bar I^-.
$$
Taking into account that  $x\in \cl W(\bar I^+, \bar I^0, \bar I^-)$, we also get  $(\lambda,b,x)\in \mathbb{G}(A, (\bar I^+,\bar I^0, \bar I^-))$.

{\bf Step 4.} We prove   $\bigcup_{(I^+, I^0, I^-)\in \mathcal{J}(A)}\mathbb{G}(A, (I^+,I^0, I^-))\subset \gph   S(\cdot, \cdot, A)$. First fix any $(I^+,I^0, I^-)\in\mathcal J(A)$ and $(\lambda, b,x)\in \ri \mathbb{G}(A, (I^+,I^0, I^-))$. Then by (\ref{ri-S_F}) we have  $\lambda>0$ and 
$$
x_i>0,\ A_i^\top (b-Ax)=\lambda\ \forall i\in I^+,\; x_i=0, \ |A_i^\top (b-Ax)| < \lambda\ \forall i\in I^0,\;x_i<0,\ A_i^\top (b-Ax) =-\lambda\ \forall i\in I^-.
$$
It then follows from  (\ref{norm-1-subdiff}) that  $A^\top(b-Ax)\in\lambda\partial \|\cdot\|_1(x)$  and hence that  $0\in \rge A^\top+\partial \|\cdot\|_1(x)$. According to \eqref{lam-b-xs'}, we get  $(\lambda, b,x)\in \gph   S(\cdot, \cdot, A)$ and hence    $\ri \mathbb{G}(A, (I^+,I^0, I^-))\subset \gph   S(\cdot, \cdot, A)$.    As $\gph   S(\cdot, \cdot, A)$ is closed, we  have $  \mathbb{G}(A, (I^+,I^0, I^-))\subset \gph   S(\cdot, \cdot, A)$.   This completes the proof.
\end{proof}

As a consequence of  partitioning  $\gph S(\cdot, \cdot, A)$ into polyhedral cones  $\mathbb{G}( A, (I^+,I^0, I^-))$, where $(I^+,I^0,I^-)\in \mathcal{J}(A)$,  the following corollary domenstrates   that $\dom S(\cdot,\cdot, A)$ can be partitioned into polyhedral cones $P_{1,2}(\mathbb{G}( A, (I^+,I^0, I^-)))$,  any distinct  two of which intersect only on their boundaries.  Recall from Lemma \ref{lem-basic-pro} (i) that the polyhedral cones $P_{1,2}(\mathbb{G}( A, (I^+,I^0, I^-)))$ have nonempty interiors, where $(I^+,I^0, I^-)\in  \mathcal{J} (A)$.

\begin{corollary}\label{cor-domain-partition}
Let $A\in \R^{m\times n}$. We have  
\[
 \dom S(\cdot,\cdot, A)=\bigcup_{ (I^+,I^0, I^-)\in  \mathcal{J} (A)} P_{1,2}(\mathbb{G}( A, (I^+,I^0, I^-))),
   \]
where,   for any  distinct  $(I^+,I^0, I^-)$ and $(\bar I^+, \bar I^0, \bar I^-)$ in $\mathcal J (A)$, it holds that 
      $$\inte P_{1,2}(\mathbb{G}( A, (I^+,I^0, I^-)))\cap  P_{1,2}(\mathbb{G}(A, (\bar I^+, \bar I^0, \bar I^-)))=\emptyset.$$  
\end{corollary}
\begin{proof}  First   we  get  from  (\ref{key-decom}) the formula for $\dom S(\cdot, \cdot, A)$  by observing that $\dom S(\cdot, \cdot, A)=P_{1,2}(\gph   S(\cdot, \cdot, A))$.  Fix any  distinct  $(I^+,I^0, I^-)$ and $(\bar I^+, \bar I^0, \bar I^-)$ in $\mathcal J (A)$, and let $$F:=\mathbb{F}(A,(I^+,I^0,I^-))\quad \mbox{and}\quad \bar F:=\mathbb{F}(A, (\bar I^+, \bar I^0, \bar I^-)).$$  Then we get from Lemma  \ref{lem-one-to-one} that $F, \bar F\in \mathcal{F}(A)$ and $F\not=\bar F$. 
To show  that $\inte P_{1,2}(\mathbb{G}( A, (I^+,I^0, I^-)))$ and $P_{1,2}(\mathbb{G}(A, (\bar I^+,\bar I^0, \bar I^-)))$ have no common point, suppose by contradiction that    $(\lambda, b)$ is contained in both of them.  Note that both $P_{1,2}(\mathbb{G}( A, (I^+,I^0, I^-)))$ and $P_{1,2}(\mathbb{G}(A, (\bar I^+,\bar I^0, \bar I^-)))$ are convex.   
  Then we can find some $(\lambda', b')\in \inte P_{1,2}(\mathbb{G}(A, (\bar I^+,\bar I^0, \bar I^-)))$ and some $(\tilde{\lambda}, \tilde{b})$ close enough to $(\lambda, b)$ in the line segment $((\lambda, b), (\lambda', b')]$ such that 
  $$(\tilde{\lambda}, \tilde{b})\in \inte P_{1,2}(\mathbb{G}( A, (I^+,I^0, I^-)))\cap \inte P_{1,2}(\mathbb{G}(A, (\bar I^+,\bar I^0, \bar I^-))).$$   
  By  Lemma \ref{lem-basic-pro} (i),  we have $\tilde{\lambda}>0$ and   $\tilde{b}=Ax+\tilde{\lambda} y=Ax'+\tilde{\lambda} y'$ for some $x\in W(I^+,I^0, I^-)$, $x'\in W(\bar I^+, \bar I^0, \bar I^-)$, $y\in \ri F$ and $y'\in \ri \bar F$. By (\ref{ri-S_F}), we have $(\tilde{\lambda}, \tilde{b}, x)\in \ri \mathbb{G}(A, (I^+,I^0, I^-))$ and $(\tilde{\lambda}, \tilde{b}, x')\in \ri \mathbb{G}(A, (\bar I^+, \bar I^0, \bar I^-))$. In view of (\ref{key-decom}), we thus have $x, x'\in S(\tilde{\lambda}, \tilde{b}, A)$. Then we get from \cite[Lemma 4.1]{zhang2015yc} that $Ax=Ax'$ and hence $y=y'$ (i.e., the faces $F$ and $\bar F$ of $\mathbb{Y}(A)$  have a common point in their relative interiors).  By \cite[Theorem 18.1]{CvxAn}, we have $F=\bar F$, a contradiction.  This completes the proof.   
  \end{proof}
\begin{remark}\label{rem-kkk}
     Let $(\lambda, b, A)\in \dom S$.  By the definition of  $\mathcal{K}$ (cf. \eqref{def-active}),   Corollary \ref{cor-domain-partition} implies that $\dom \mathcal{K}=\dom S$, and that   $(\lambda, b)\in \inte P_{1,2}(\mathbb{G}( A, (I^+,I^0, I^-)))$ if and only if  $\mathcal{K}(\lambda,b,A)$ is a singleton. Moreover, when $\lambda > 0$, the union of all polyhedral cones $P_{1,2}(\mathbb{G}(A, (I^+, I^0, I^-)))$ with $(I^+, I^0, I^-) \in \mathcal{K}(\lambda, b, A)$ contains $(\lambda, b)$ in its interior.
\end{remark}

Next we verify condition (ii) in Lemma \ref{lem-Lip-con-poly-theo-ori}.

\begin{proposition}\label{convex} Let $A\in \R^{m\times n}$. For each  $(I^+, I^0, I^-)\in \mathcal{J}(A)$, the graph  of the restriction of  $S(\cdot, \cdot, A) $ on $P_{1,2}(\mathbb{G}( A, (I^+,I^0, I^-)))$ is convex. 
\end{proposition}

\begin{proof} 
Let $(I^+, I^0, I^-)\in \mathcal{J}(A)$. For any $(\hat\lambda,\hat b,\hat x),\;(\bar\lambda,\bar b,\bar x)\in \gph   S(\cdot, \cdot, A)|_{P_{1,2}(\mathbb{G}( A, (I^+,I^0, I^-)))}$ and $\theta\in (0,1)$.
 As $(\hat\lambda,\hat b),(\bar\lambda,\bar b)\in P_{1,2}(\mathbb{G}( A, (I^+,I^0, I^-)))$,   there exist $\hat x',\bar x'\in\R^n$ such that 
$
(\hat\lambda,\hat b,\hat x'),\;(\bar\lambda,\bar b,\bar x')\in \mathbb{G}(A, (I^+, I^0, I^-)).
$
Let
\begin{equation}\label{S-Lips2}
(\lambda_\theta,b_\theta,x_\theta') := \theta (\hat\lambda,\hat b,\hat x') + (1-\theta) (\bar\lambda,\bar b,\bar x') \mbox{ and } x_\theta:=\theta\hat x+(1-\theta)\bar x.
\end{equation}
It follows from the convexity of $\mathbb{G}(A, (I^+, I^0, I^-))$       that $(\lambda_\theta,b_\theta,x_\theta')\in \mathbb{G}(A, (I^+, I^0, I^-))$ and hence by (\ref{key-decom}) that
\begin{equation}\label{S-Lips3}
 x_\theta'\in S(\lambda_\theta, b_\theta,A).
\end{equation}
Moreover, as 
\[
 \mathbb{G}(A, (I^+, I^0, I^-))\subset \gph   S(\cdot, \cdot, A)|_{P_{1,2}(\mathbb{G}( A, (I^+,I^0, I^-)))}\subset \gph   S(\cdot, \cdot, A),
 \]
 we get $\hat x,\hat x'\in S(\hat \lambda,\hat b, A)$ and $\bar x,\bar x'\in S(\bar\lambda,\bar b, A)$. Then   it follows  that
\begin{equation}\label{S-Lips4}
A\hat x=A\hat x',~~\|\hat x\|_1=\|\hat x'\|_1~~\mbox{and}~~A\bar x=A\bar x',~~\|\bar x\|_1=\|\bar x'\|_1,
\end{equation}
which follow  from  \cite[Lemma 4.1]{zhang2015yc} for the cases of $\hat \lambda>0$ or $\bar \lambda>0$, and   from the definition for the cases of $\hat \lambda=0$ or $\bar \lambda=0$.

If $\hat\lambda=\bar\lambda=0$, then $\lambda_\theta=0$. Thus, $Ax_\theta=\theta A\hat{x}+(1-\theta)A\bar{x}=\theta \proj_{\rge A} \hat b+(1-\theta) \proj_{\rge A} \bar b=\proj_{\rge A}  b_\theta$ and
$$
\|x_\theta\|_1=\|\theta\hat x+(1-\theta)\bar x\|_1\leq\theta\|\hat x\|_1+(1-\theta)\|\bar x\|_1\hskip-0.1cm=\theta\|\hat x'\|_1+(1-\theta)\|\bar x'\|_1\hskip-0.1cm=\|x_\theta'\|_1,
$$
where the second last equality  follows from (\ref{S-Lips4}), the   last equality comes from \eqref{S-Lips2}  and the linearity of $\|\cdot\|_1$  on $P_3(\mathbb{G}( A, (I^+,I^0, I^-)))$ (see Lemma \ref{lem-basic-pro} (iii)).
In view of \eqref{S-Lips3},  we get  $x_\theta\in S(\lambda_\theta,b_\theta, A)$ by definition.

If $\hat\lambda>0$ or $\bar \lambda>0$, then $\lambda_\theta>0$. Thus
\begin{equation*}
\begin{array}{ll}
&\frac{1}{2\lambda_\theta}\|Ax_\theta-\proj{_{\rge A}b_\theta}\|^2+\|x_\theta\|_1\\
=&\frac{1}{2\lambda_\theta}\|\theta(A\hat x-\proj_{\rge A}\hat b)+(1-\theta)(A\bar x-\proj_{\rge A}\bar b)\|^2+\|\theta\hat x+(1-\theta)\bar x\|_1\\
\leq&\frac{1}{2\lambda_\theta}\|\theta(A\hat x-\proj_{\rge A}\hat b)+(1-\theta)(A\bar x-\proj{_{\rge A}}\bar b)\|^2+\theta\|\hat x\|_1+(1-\theta)\|\bar x\|_1\\
=&\frac{1}{2\lambda_\theta}\|\theta(A\hat x'-\proj_{\rge A}\hat b)+(1-\theta)(A\bar x'-\proj_{\rge A}\bar b)\|^2+\theta\|\hat x'\|_1+(1-\theta)\|\bar x'\|_1\\
=&\frac{1}{2\lambda_\theta}\|Ax_\theta'-\proj_{\rge A}b_\theta\|^2+\|x_\theta'\|_1,
\end{array}
\end{equation*}
where the first equality comes from \eqref{S-Lips2}, the second last equality comes from \eqref{S-Lips4}, and the   last equality comes from \eqref{S-Lips2} and the linearity of $\|\cdot\|_1$  on $P_3(\mathbb{G}( A, (I^+,I^0, I^-)))$ (see Lemma \ref{lem-basic-pro} (iii)). Then in light of \eqref{S-Lips3}, we  get  $x_\theta\in S(\lambda_\theta,b_\theta, A)$ by definition. 

In both cases, we have shown that  $x_\theta\in S(\lambda_\theta,b_\theta, A)$. Noting that due to $
(\hat\lambda,\hat b,\hat x), (\bar\lambda,\bar b,\bar x)\in \mathbb{G}(A, (I^+, I^0, I^-))$, we have $(\lambda_\theta,b_\theta) \in P_{1,2}(\mathbb{G}( A, (I^+,I^0, I^-)))$. Thus 
$$(\lambda_\theta, b_\theta, x_\theta) \in \gph   S(\cdot, \cdot, A)|_{P_{1,2}(\mathbb{G}( A, (I^+,I^0, I^-)))}.$$ This verifies the convexity of  $\gph  S|_{P_{1,2}(\mathbb{G}( A, (I^+,I^0, I^-)))}$ and thus completes the proof.
\end{proof}

Now noting that $\dom S(\cdot, \cdot, A)=\R_+\times \R^m$ is convex,  we are ready to obtain the global Lipschitz continuity  of $S(\cdot, \cdot, A)$ on its domain  by applying Lemma \ref{lem-Lip-con-poly-theo-ori} and Propositions \ref{prop-gphS-lambda-b}  and  \ref{convex} in a straightforward way.

{
\begin{theorem}\label{theo-S-Lips}
Let $A\in \R^{m\times n}$. Then  $S(\cdot, \cdot, A) $ is Lipschitz continuous on its domain.
\end{theorem}
}

Next, we present  several  properties of $S$ below, which are useful in the next section. In particular the explicit representation  of $S$ in \eqref{def-S-lambda-b-x}, a consequence of the Lipschitz continuity of $S$ with respect to  $\lambda$ and $b$, will play a key role in exploring the Lipschitz continuity of $S$ with respect to    $A$. 

\begin{corollary}\label{cor-big}
Let $A\in \R^{m\times n}$. The following hold:
\begin{description}
    \item[(i)] For all $(\lambda, b, A)\in \dom S$ and $(I^+, I^0, I^-)\in  \mathcal{K}(\lambda, b, A)$,   we have 
 \begin{equation}\label{def-S-lambda-b-x}
    \begin{array}{rl}
   S(\lambda,b, A)=&\left\{
    x \in \R^n
    \left|
    \begin{array}{lcl}
    x_i\geq 0, &\hskip-0.2cm A_i^\top (b-Ax)=\lambda &\hskip-0.2cm \forall i\in I^+\\
     x_i=0, &\hskip-0.2cm  | A_i^\top (b-Ax)|\leq \lambda&\hskip-0.2cm \forall i\in I^0\\
     x_i\leq 0, &\hskip-0.2cm A_i^\top (b-Ax)=-\lambda &\hskip-0.2cm \forall i\in I^-
     \end{array}
     \right.
     \right\}.
     \end{array}
    \end{equation}
  \item[(ii)] If $(I^+,I^0,I^-)\in \mathcal{J}^0(A)$, then $S(\cdot, \cdot, A)$ is  single-valued and  linear  on $P_{1,2}(\mathbb{G}(A,(I^+,I^0,I^-)))$ with      
     $$S(\lambda,b,A)=\{x\}\quad \forall (\lambda, b)\in P_{1,2}(\mathbb{G}(A,(I^+,I^0,I^-))),$$
     where $x$ is      the unique solution to the system of linear  equations 
      \begin{equation}\label{equ-single-valued}
     \left\{
     \begin{array}{rlll}
 A_i^\top  Ax&=&A_i^\top b-\lambda  & \forall i\in I^+,\\
     x_i&=&0 & \forall i\in I^0,\\
 A_i^\top  Ax&=&A_i^\top b+\lambda & \forall i\in I^-.
     \end{array}
     \right.
    \end{equation}
    If $(I^+,I^0,I^-) \in \mathcal{J}(A)\backslash\mathcal{J}^0(A)$,   $S(\cdot, \cdot, A)$ is not single-valued on  $\inte P_{1,2}(\mathbb{G}(A,(I^+,I^0,I^-)))$.  
 \item[(iii)] If $S$ is single-valued at $(\lambda, b, A)$ with $S(\lambda,b,A)=\{x\}$, then $(I^+(x),I^0(x),I^-(x))\in \mathcal{K}(\lambda, b, A)\cap \mathcal{J}^0(A)$, where $(I^+(x),I^0(x),I^-(x))$ is an index partition of $I$ defined by \eqref{def-index-partition-signs-x}. 
     \item[(iv)] For any  distinct index partitions $(I^+,I^0, I^-)$ and $(\bar I^+, \bar I^0, \bar I^-)$ in $\mathcal J (A)$ and any  $(\lambda, b, x), (\lambda', b',x')\in \gph   S(\cdot, \cdot, A)$ such that $(\lambda, b)\in \inte P_{1,2}(\mathbb{G}( A, (I^+,I^0, I^-)))$ and $(\lambda', b')\in \inte P_{1,2}(\mathbb{G}(A, (\bar I^+, \bar I^0, \bar I^-)))$, there is some $\theta\in (0, 1)$ such that
\begin{equation}\label{zxd-all-in}
\theta (\lambda, b, x)+(1-\theta)(\lambda', b',x')\not\in \gph   S(\cdot, \cdot, A).
\end{equation}
 \end{description}
 \end{corollary}
\begin{proof} 
First we show (i).  By the partition of $\dom S(\cdot, \cdot, A)$ presented in  Corollary \ref{cor-domain-partition}, we  have
\begin{equation}\label{bibi}
\gph   S(\cdot, \cdot, A)|_{\inte P_{1,2}(\mathbb{G}( A, (I^+,I^0, I^-)))}\subset \mathbb{G}(A, (I^+, I^0, I^-)).
\end{equation}
Fix any $(\lambda,b)\in \bdry P_{1,2}(\mathbb{G}( A, (I^+,I^0, I^-)))$ and  $x\in S(\lambda,b, A)$.  Since $\inte P_{1,2}(\mathbb{G}( A, (I^+,I^0, I^-)))\not=\emptyset$ (cf. Lemma \ref{lem-basic-pro} (i)), there exist  $(\lambda^\nu,b^\nu)\in \inte P_{1,2}(\mathbb{G}( A, (I^+,I^0, I^-)))$ such that $(\lambda^\nu,b^\nu)\to(\lambda,b)$. According to the Lipschitz continuity of $S(\cdot, \cdot, A)$ established in Theorem \ref{theo-S-Lips},  there exist $x^\nu\in\R^n$ such that $x^\nu \in S(\lambda^\nu,b^\nu, A)$ for all $\nu$ and $x^\nu\to x$. In view of (\ref{bibi}), we  have $(\lambda^\nu,b^\nu,x^\nu)\in \mathbb{G}(A, (I^+, I^0, I^-))$ for all $\nu$. As $\mathbb{G}(A, (I^+, I^0, I^-))$ is closed, we have  $(\lambda,b,x)\in \mathbb{G}(A, (I^+, I^0, I^-))$. This verifies $\gph   S(\cdot, \cdot, A)|_{\bdry P_{1,2}(\mathbb{G}( A, (I^+,I^0, I^-)))}\subset \mathbb{G}(A, (I^+, I^0, I^-))$, which, together with \eqref{bibi}, implies that  $\gph   S(\cdot, \cdot, A)|_{ P_{1,2}(\mathbb{G}( A, (I^+,I^0, I^-)))}\subset \mathbb{G}(A, (I^+, I^0, I^-))$.  In light of  (\ref{key-decom}), we get $\mathbb{G}(A, (I^+, I^0, I^-))\subset \gph   S(\cdot, \cdot, A)|_{P_{1,2}(\mathbb{G}( A, (I^+,I^0, I^-)))}$. This yields  the equality 
\[
\gph   S(\cdot, \cdot, A)|_{ P_{1,2}(\mathbb{G}( A, (I^+,I^0, I^-)))}=\mathbb{G}(A, (I^+, I^0, I^-)).
\]
By  the definitions of $\mathbb{G}(A, (I^+, I^0, I^-))$ and $\mathcal{K}(\lambda, b, A)$,  we get  (i) in a straightforward way.

 Next we  show (ii).   It follows from  (i) that   any $x\in S(\lambda,b,A)$ with $(\lambda, b)\in  P_{1,2}(\mathbb{G}( A, (I^+,I^0, I^-)))$ is a 
  solution to the system \eqref{equ-single-valued}  of linear  equations. 
  
If $(I^+,I^0,I^-)\in \mathcal{J}^0(A)$,  
 then the columns $A_i$ of $A$ with $i\in I^+\cup I^-$ are linearly independent by definition. Consequently,  any  solution $x$  to   \eqref{equ-single-valued} must be unique and linearly dependent on $(\lambda, b)$.  This implies that $S(\cdot, \cdot, A)$ is  single-valued and linear  on  $P_{1,2}(\mathbb{G}( A, (I^+,I^0, I^-)))$. 
 
 Alternatively,  let $(I^+,I^0,I^-)\in \mathcal{J}(A)\backslash \mathcal{J}^0(A)$ and fix any $(\lambda, b)\in \inte P_{1,2}(\mathbb{G}( A, (I^+,I^0, I^-)))$. Then by Lemma \ref{lem-basic-pro} (i),  there is some  $\bar x\in S(\lambda, b, A)$ such that $(\lambda,b, \bar x)\in \ri \mathbb{G}(A, (I^+, I^0, I^-))$. According to \eqref{ri-S_F} , we have $\lambda>0$ and $b=A\bar x+\lambda y$ for some  
$\bar x\in W(I^+,I^0,I^-)$ and $y\in \ri \mathbb F(A,(I^+,I^0,I^-))$. As  the vectors $A_i$ with $i\in I^+\cup I^-$ are linearly dependent due to $(I^+,I^0, I^-)\not\in \mathcal{J}^0(A)$,  there is  some $w\not=0$ such that $Aw=0$ and $w_i=0$ for all $i\in I^0$. Then $\bar x +t^*w\in W(I^+,I^0,I^-)$ for some $t^*\not=0$ small enough. It  follows again from \eqref{ri-S_F} that $(\lambda, A(\bar x +t^*w)+\lambda y, \bar x+t^*w)\in \ri \mathbb{G}(A, (I^+, I^0, I^-))$. 
In view of $Aw=0$ and $b=A\bar x+\lambda y$, we get   $(\lambda, b, \bar x+t^*w)\in \ri \mathbb{G}(A, (I^+, I^0, I^-))$. So we get from Proposition \ref{prop-gphS-lambda-b} that 
 $\bar x+t^*w\in S(\lambda, b, A)$. As $w\not=0$ and $t^*\not=0$, we have $\bar x+t^*w\not=\bar x$. Given that we have shown that  $\bar x\in S(\lambda, b, A)$ and $\bar x+t^*w\in S(\lambda, b, A)$,  $S$ is obviously  not single-valued at  $(\lambda, b, A)$.

To show (iii),  we first get from  \cite[Theorem 2.1]{zhang2015yc} (see also \cite[Proposition 3.2]{Gilbert17}) that    Condition 1.1 is satisfied at $(\lambda, b,A)$ for $x$. Then we have $(I^+(x),I^0(x),I^-(x))\in \mathcal{J}^0(A)$ by definition. Given that $x\in S(\lambda, b, A)$, we get  from \eqref{lam-b-xs'} in Proposition \ref{prop-gphS-lambda-b} that $A^\top(   b-A   x)\in  \lambda\partial\|\cdot\|_1(x)$.  By the formula  (\ref{norm-1-subdiff}) for the subdifferential set of $\|\cdot\|_1$, we get
$$
    \begin{array}{lcl}
        A_i^\top (b-Ax)=   \lambda\; \forall i\in I^+(x),\;\; |A_i^\top (b-Ax)|\leq\lambda\; \forall i\in I^0(x),\;\; A_i^\top (b-Ax)=-\lambda\;\forall i\in I^-(x).
     \end{array}
$$
This, together with the fact that  $\lambda\geq 0$, implies that  $(\lambda,b,x)\in \mathbb{G}(A, (I^+(x), I^0(x), I^-(x)))$ and hence $(\lambda, b)\in P_{1,2}(\mathbb{G}(A, (I^+(x), I^0(x), I^-(x))))$ by definition. It then follows from the definition  of $\mathcal{K}$  that $(I^+(x), I^0(x), I^-(x))\in \mathcal{K}(\lambda, b, A)$.

It remains to  show (iv).  Suppose by contradiction that
\[
(\lambda_\theta, b_\theta, x_\theta):=\theta (\lambda, b, x)+(1-\theta)(\lambda', b',x')\in \gph   S(\cdot, \cdot, A)\quad\forall \theta\in [0, 1].
\]
Then in view of $(\lambda, b)\in \inte P_{1,2}(\mathbb{G}( A, (I^+,I^0, I^-)))$, we have $(\lambda_\theta, b_\theta)\in \inte P_{1,2}(\mathbb{G}( A, (I^+,I^0, I^-)))$  for all $\theta$ close enough to 1 in $[0, 1]$. 
It then follows from (i) and Corollary \ref{cor-domain-partition} that  $(\lambda_\theta, b_\theta, x_\theta)\in \mathbb{G}(A, (I^+, I^0, I^-))$ for all $\theta$ close enough to 1 in $[0, 1]$, and hence  that $(\lambda_\theta, b_\theta, x_\theta)\in \aff \mathbb{G}(A, (I^+, I^0, I^-))$ for all $\theta\in [0, 1]$. By letting $\theta=0$, we have
  $(\lambda', b',x')\in \aff \mathbb{G}(A, (I^+, I^0, I^-))$. According to  (\ref{def-sf-aff}), we have
  \begin{equation}\label{affine-jieguo}
  A_i^\top \left(\frac{b'-Ax'}{\lambda'}\right)=1\quad \forall i\in I^+.
  \end{equation}
  As $(\lambda', b')\in\inte P_{1,2}(\mathbb{G}(A, (\bar I^+, \bar I^0, \bar I^-)))$, we get from Lemma \ref{lem-basic-pro} (i) the existence of some $\tilde{x}'\in W (\bar I^+, \bar I^0, \bar I^-)$ such that $\frac{b'-A\tilde{x}'}{\lambda'}\in \ri \mathbb F(A, (\bar I^+, \bar I^0, \bar I^-))$.
  It follows from  Lemma \ref{lem-basic-proSF} that $(\lambda', b', \tilde{x}')\in\ri \mathbb{G}(A, (\bar I^+, \bar I^0, \bar I^-))$. In view of the fact that $x', \tilde{x}'\in S(\lambda', b',A)$, we get from \cite[Lemma 4.1]{zhang2015yc} that $A\tilde{x}'=Ax'$. Thus, we have $\frac{b'-Ax'}{\lambda'}\in \ri \mathbb F(A, (\bar I^+, \bar I^0, \bar I^-))$. 
  Then by \eqref{xiangduineibu} and (\ref{affine-jieguo}),  we get  $I^+\subset \bar I^+$.
  In a parallel manner, we can show $\bar I^+\subset I^+$.   That is, we have $I^+=\bar I^+$. Similarly, we obtain $I^-= \bar I^-$. 
  So we have  $(I^+,I^0, I^-)=(\bar I^+, \bar I^0, \bar I^-)$, a contradiction. 
This completes the proof. \end{proof}


  Below  we apply Corollary \ref{cor-big} (iv) to show that the partition of  $\gph S(\cdot, \cdot, A)$ presented  in  Proposition \ref{prop-gphS-lambda-b}   is optimal and unique. 
\begin{proposition}[Optimal and Unique Partition]\label{prop-NS}
Let $A\in \R^{m\times n}$. 
 If  there are  closed and convex sets $C_l$ ($l=1,\cdots, L$) such that
\[
\gph   S(\cdot, \cdot, A)=\bigcup_{l=1}^L C_l,
\]
then $L\geq |\mathcal{J}(A)|$. If  in addition  $L=|\mathcal{J}(A)|$, then
  $$
 \{C_l\mid l=1,\cdots, L\}= \{\mathbb{G}(A,(I^+,I^0,I^-))\mid (I^+,I^0,I^-)\in \mathcal{J}(A)\}.
  $$ 
  \end{proposition}
   \begin{proof}  To show  $L\geq |\mathcal{J}(A)|$, suppose by contradiction that  $L<|\mathcal{J}(A)|$.  Following   (\ref{key-decom}) and Lemma \ref{lem-basic-proSF}, 
   there are some $C_l$ with $1\leq l\leq L$ and two distinct partitions  $(I^+, I^0, I^-)$ and $(\bar I^+, \bar I^0, \bar I^-)$ in $\mathcal{J}(A)$  such that $C_l\cap \ri \mathbb{G}(A, (I^+, I^0, I^-))\not=\emptyset$ and $C_l\cap \ri \mathbb{G}(A, (\bar I^+, \bar I^0, \bar I^-))\not=\emptyset$. 
   Let $(\bar\lambda, \bar b, \bar x)\in C_l\cap \ri \mathbb{G}(A, (I^+, I^0, I^-))$ and  $(\tilde\lambda , \tilde b, \tilde x)\in C_l\cap  \ri \mathbb{G}(A, (\bar I^+, \bar I^0, \bar I^-))$ be two distinct triples and  $(\lambda_\theta, b_\theta, x_\theta):=(1-\theta) (\bar\lambda, \bar b, \bar x)+\theta(\tilde\lambda, \tilde b, \tilde x)$ for all $\theta \in [0,1]$. Then it follows from Lemma \ref{lem-basic-pro} (ii) that $(\bar\lambda, \bar b)\in \inte P_{1,2}(\mathbb{G}( A, (I^+,I^0, I^-)))$  and  $(\tilde\lambda, \tilde b)\in \inte P_{1,2}(\mathbb{G}(A, (\bar I^+, \bar I^0,  \bar I^-)))$.
As $C_l$ is convex, we have $$(\lambda_\theta, b_\theta, x_\theta)\in C_l\subset \gph   S(\cdot, \cdot, A), \quad 
\forall \theta\in [0, 1],$$ contradicting to Corollary  \ref{cor-big} (iv). This indicates that  $L\geq |\mathcal{J}(A)|$.

Now we consider the case that $L=|\mathcal{J}(A)|$. For the sake of simplicity of notation, let $\mathcal{J}(A):=\{(I^+_1, I^0_1, I^-_1),\cdots, (I^+_L, I^0_L, I^-_L)\}$.  According to the previous argument, each $C_l$ can contain at most one $\ri \mathbb{G}(A, (I^+_l, I^0_l, I^-_l))$. Without loss of generality  we can assume that  
$ \ri \mathbb{G}(A, (I^+_l, I^0_l, I^-_l))\subset C_l$ for all $l=1,\dots, L$. As  each $C_l$ is  closed by assumption, we get 
\begin{equation} \label{unincl}
\mathbb{G}(A, (I^+_l, I^0_l, I^-_l)) \subset C_l,\quad \forall l=1,\dots, L.
\end{equation}
 To show that the above inclusion is actually an equality for completing the proof, suppose by contradiction that   there is 
$(\lambda^*, b^*, x^*)\in C_{k}\backslash \mathbb{G}(A, (I^+_{k}, I^0_{k}, I^-_{k}))$ for some $k\in \{1,\dots, L\}$.  Let  $(\lambda', b')\in \inte P_{1,2}(\mathbb{G}(A, (I^+_{k}, I^0_{k}, I^-_{k})))$ and  
$x'\in S(\lambda', b', A)$. Then we get from Corollary \ref{cor-domain-partition} and Corollary \ref{cor-big} (i)  that  $(\lambda^*, b^*)\not\in P_{1,2}(\mathbb{G}(A, (I^+_{k}, I^0_{k}, I^-_{k})))$ 
and 
$(\lambda', b', x')\in \mathbb{G}(A, (I^+_{k}, I^0_{k}, I^-_{k}))$,  implying that $(\lambda', b', x')\in C_{k}$. By definition we have 
$(\lambda^*, b^*)\in P_{1,2}(C_k)$ and hence 
$P_{1,2}(\mathbb{G}(A, (I^+_{k}, I^0_{k}, I^-_{k}))) \subsetneq  P_{1,2}(C_{k})$. 
Noting that $P_{1,2}(\mathbb{G}(A, (I^+_{k}, I^0_{k}, I^-_{k})))$ is closed and convex 
with  $(\lambda', b')$ being an interior of  $P_{1,2}(\mathbb{G}(A, (I^+_{k}, I^0_{k}, I^-_{k})))$
(hence of $P_{1,2}(C_{k})$) and 
that $P_{1,2}(C_{k})$ is convex with  $$(\lambda^*, b^*)\in P_{1,2}(C_k)\backslash P_{1,2}(\mathbb{G}(A, (I^+_{k}, I^0_{k}, I^-_{k}))),$$ 
we obtain the existence of some $(\lambda, b)$ close enough to $(\lambda^*, b^*)$ in the line segment $((\lambda', b'), (\lambda^*, b^*))$ such that $(\lambda, b)$ is an exterior of $P_{1,2}(\mathbb{G}(A, (I^+_{k}, I^0_{k}, I^-_{k})))$ and an interior of $P_{1,2}(C_{k})$. Then there is an neighborhood $V$ of $(\lambda, b)$ such that 
\[
P_{1,2}(\mathbb{G}(A, (I^+_{k}, I^0_{k}, I^-_{k})))\cap V=\emptyset\quad\mbox{and}\quad V\subset P_{1,2}(C_{k})\subset\dom S(\cdot, \cdot, A). 
\]   
It then follows from the partition of $\dom S(\cdot, \cdot, A)$ given in 
Lemma \ref{lem-basic-pro} (ii) that there is some $j\in \{1,\cdots, L\}$ such that  $j\not=k$ and  
    $\inte  P_{1,2}(\mathbb{G}(A, (I^+_{j}, I^0_{j}, I^-_{j}))) \cap V\not=\emptyset$.  
    
    Let $(\lambda, b)\in    \inte  P_{1,2}(\mathbb{G}(A, (I^+_{j}, I^0_{j}, I^-_{j})))\cap V$. Then we get $(\lambda, b)\in  P_{1,2}(C_{k})$  and hence  $(\lambda, b, x)\in C_{k}$ for some $x\in \R^n$.   As $C_{k}$ is convex and we have already shown $(\lambda', b', x')\in C_{k}$, we thus have  
$
[(\lambda, b, x),\,(\lambda', b', x')]\subset C_{k}\subset \gph   S(\cdot, \cdot, A),
$
contradicting  to  Corollary \ref{cor-big} (iv)  due to  $(\lambda, b)\in \inte P_{1,2}(\mathbb{G}(A, (I^+_{j}, I^0_{j}, I^-_{j})))$ and $(\lambda', b')\in \inte  P_{1,2}(\mathbb{G}(A, (I^+_{k}, I^0_{k}, I^-_{k})))$. Thus \eqref{unincl} holds with an equality.
\end{proof}

Below we apply \eqref{def-S-lambda-b-x} and \eqref{equ-single-valued} to give an analytic expression of the solution multifunction for an example of \eqref{ext-lasso-problem-lambda}.

\begin{example}\label{exam-nice}
Let $m=2$, $n=3$ and $I=\{1,2, 3\}$. Consider the following matrix of full row rank:
\[
A=
\left[
\begin{array}{ccc}
  1 & 0 & 2 \\
  0 & 2 & -2
\end{array}
\right],
\]
which has appeared  in both \cite[Example 4.18]{Berk2023bh} and \cite{zhang2015yc}. Clearly, we have   $$\mathbb{Y}(A) =\co\{V_1, V_2, V_3, V_4\},$$ 
where $V_1:=(1, 1/2)^\top$, $V_2:=(0, -1/2)^\top$, $V_3:=(-1, -1/2)^\top$ and $V_4:=(0, 1/2)^\top$ are extreme points of $\mathbb{Y}(A)$.  Clearly, the parallelogram $\mathbb{Y}(A)$ has 9 nonempty faces $F_i$ ($i=1,\cdots, 9$). According to  Lemmas \ref{lem-xxys} and \ref{lem-one-to-one}, each $F_i$ corresponds one-to-one to an index partition $(I_i^+, I_i^0, I_i^-)$ in $\mathcal{J}(A)$ with 
$$
(I_i^+, I_i^0, I_i^-)=(\mathcal{A}^+(y_i, F_i),\mathcal{A}^0(y_i, F_i),\mathcal{A}^-(y_i, F_i))\quad \forall y_i\in \ri F_i.
$$
It  follows from  Corollary \ref{cor-domain-partition} that 
$\dom S(\cdot, \cdot, A)$ can be partitioned into 9 polyhedral cones  
$$D_i:=P_{1,2}(\mathbb{G}(A, (I_i^+, I_i^0, I_i^-)))\quad \forall i=1,\dots, 9,$$
whose extreme directions can be calculated, according to Lemma \ref{lem-basic-pro} (ii), in terms of columns $A_j$ of $A$ and extreme points $V_k$ of $\mathbb{Y}(A)$. We record all  extreme directions of the polyhedral cones  $D_i$ as columns $U_l$  ($l=1,\cdots,8$) of the matrix
\[
U:=\left[
\begin{array}{cccccccc}
  0 & 0 & 0 & 0 & 1 & 1 & 1 & 1 \\
  A_2 & A_3 & -A_2 & -A_3& V_1 & V_2 & V_3 & V_4  \\
\end{array}
\right].
\]

\begin{table}[h]
\begin{tabular}{lllll}
\hline
$i$  & $(I_i^+, I_i^0, I_i^-)$ & Ext($F_i$)    &Ext-Dir($D_i$) & $S|_{D_i}(\cdot, \cdot, A)$ \\[0.2cm]
  \hline
1  & $(\{1,2,3\},\emptyset, \emptyset)$  & $V_1$ & $U_1, U_2, U_5$ &  see (\ref{duojie1})   \\[0.2cm]
\hline
2 & $(\emptyset, \emptyset, \{1,2,3\})$ & $V_3$ & $U_3, U_4, U_7$  &  see (\ref{duojie2})    \\[0.2cm]
\hline
3 & $(\{2\},\{1,3\},\emptyset)$ & $V_1, V_4$ & $U_1, U_5, U_{8}$ & $\{(0,\,-\frac{\lambda}{4}+\frac{b_2}{2},\,0)^\top\}$\\[0.2cm]
\hline
4 & $(\emptyset,\{1,3\},\{2\})$ & $V_3, V_2$  & $U_3, U_6, U_{7}$ & $\{(0,\,\frac{\lambda}{4}+\frac{b_2}{2},\,0)^\top\}$\\[0.2cm]
\hline
5 & $(\{3\},\{1,2\},\emptyset)$ &$V_1, V_2$ &  $U_2, U_5, U_{6}$ & $\{(0,\,0,\, \frac{b_1}{4}-\frac{b_2}{4}-\frac{\lambda}{8})^\top\}$\\[0.2cm]
  \hline
6 & $(\emptyset,\{1,2\},\{3\})$ & $V_3, V_4$ & $U_4, U_7, U_{8}$ & $\{(0,\,0,\, \frac{b_1}{4}-\frac{b_2}{4}+\frac{\lambda}{8})^\top\}$\\[0.2cm]
\hline
7& $(\{2\},\{1\},\{3\})$  & $V_4$ & $U_1, U_4, U_{8}$ & $\{(0,\,\frac{b_1}{2}+\frac{b_2}{2}-\frac{\lambda}{4},\,\frac{b_1}{2})^\top\}$\\[0.2cm]
  \hline
8 &$(\{3\},\{1\},\{2\})$   & $V_2$ & $U_2, U_3, U_{6}$ & $\{(0,\,\frac{b_1}{2}+\frac{b_2}{2}+\frac{\lambda}{4},\,\frac{b_1}{2})^\top\}$\\[0.2cm]
\hline
9 & $(\emptyset,\{1,2,3\},\emptyset)$  & $V_1, V_2, V_3, V_4$  & $U_5, U_6, U_{7}, U_{8}$ & $\{(0,\,0,\,0)^\top\}$\\[0.2cm]
\hline
\end{tabular}
\centering
\caption{The index partitions $(I_i^+, I_i^0, I_i^-)$ in $\mathcal{J}(A)$, the extreme points of the corresponding faces $F_i$ of $\mathbb{Y}(A)$, the extreme directions $U_l$ of the polyhedral cones  $D_i$, and the explicit representations  of the restriction  of $S(\cdot,\cdot, A)$ on $D_i$}
      \label{table-analyticalofs} 
\end{table}

By virtue of \eqref{def-S-lambda-b-x} and \eqref{equ-single-valued}, we can get the analytic expression of $S(\cdot, \cdot, A)$ on all polyhedral cones $D_i$, which  in particular gives 
\begin{equation}\label{duojie1}
      S(\cdot, \cdot, A)|_{D_1}=\left\{\left(b_1-\lambda-2t,\, \frac{b_2}{2}-\frac{\lambda}{4}+t, \,t\right)^\top\ \left | \ \left(\frac{\lambda}{4}-\frac{b_2}{2}\right)_+\leq t\leq \frac{b_1-\lambda}{2} \right. \right\} 
 \end{equation}
 and
\begin{equation}\label{duojie2}
S(\cdot, \cdot, A)|_{D_2}=\left\{\left(b_1+\lambda-2t,\, \frac{b_2}{2}+\frac{\lambda}{4}+t, \,t\right)^\top\ \left | \ \frac{b_1+\lambda}{2}\leq t\leq -\left(\frac{\lambda}{4}+\frac{b_2}{2} \right. \right)_+\right\}.
 \end{equation}
 It is evident that  $S(\cdot, \cdot, A)$ is single-valued and linear on $D_i$ with $i=3,\cdots, 9$  and single-valued on the boundaries  of $D_1$ and $D_2$, but not single-valued on the interiors of $D_1$ and $D_2$.   In Table \ref{table-analyticalofs}, we summarize the details of all index  partitions  $(I_i^+, I_i^0, I_i^-)$ in $\mathcal{J}(A)$, all extreme points  of $F_i$ (labeled Ext($F_i$)),   all extreme directions of $D_i$ (labeled Ext-Dir($D_i$)),  and all  analytic expressions of $S(\cdot, \cdot, A)$ on $D_i$ (labeled $S(\cdot, \cdot, A)|_{D_i}$).    

In Figure \ref{dom-parti}, we plot the truncated partitions of $\dom S(\cdot, \cdot, A)$ with the unit ball.
 
\begin{figure}[ht]
\vskip -3cm  
\hskip1.5cm
\includegraphics[width=10cm]{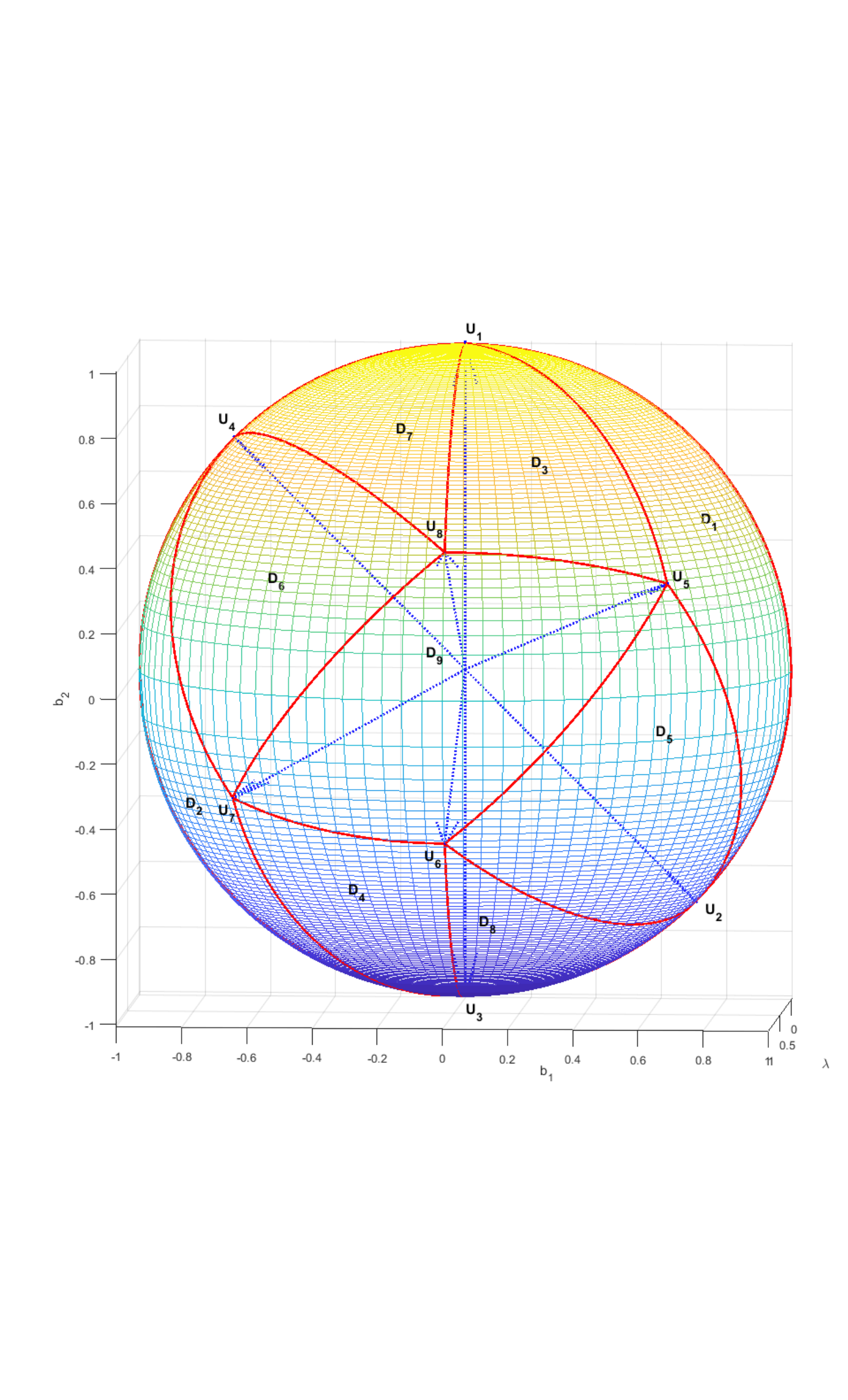}\\
\vskip-3.5cm
\caption{Tilted half unit ball with $\lambda$-axis facing the readers:  Truncations of the polyhedral cones $D_i$ ($i=1,\cdots, 9$) with the unit ball, where the $U_l$'s ($l=1,\cdots, 8$) denote their unit extreme directions}\label{dom-parti}
\end{figure}

 \end{example}


 \section{Lipschitz continuity when $A$ varies locally}\label{sec-lip-a-varies}

In this section, we will investigate the Lipschitz continuity of $S$ on the triple variables $(\lambda, b, A)$. 
To begin with, we  utilize the   multifunction $\mathcal{K}$  defined by \eqref{def-active} and  the explicit representation of $S$ presented in Corollary \ref{cor-big}  to give characterizations of Conditions 1.2 and 1.3 that have been employed to deduce the  property of single-valuedness and local Lipschitz continuity  for $S$.  

Note that the implications in \eqref{impls} for Conditions 1.1-1.3 are shown to hold for $\lambda>0$ in \cite{Berk2023bh}. We observe from the definitions that these implications also hold in the case of  $\lambda=0$. { However, we note that Condition 1.2 at $(\lambda,b, A)$ with $\lambda=0$ if and only if $A$ is of full column rank. Similarly,  Condition 1.3 at $(\lambda,b, A)$ with $\lambda=0$ if and only if $A$ is of full column rank and there is some $x\in S(\lambda,b,A)$ such that  $x_i\not=0$ for all $i$. Therefore, to address   Conditions 1.2 and 1.3 for some $(\lambda, b, A)\in \dom S$, we typically assume that $\lambda>0$ when  $A$ is not  explicitly stated to be of full column rank.}
 
{Recall from \cite{zhang2015yc} that $S$ is single-valued at some $(\bar \lambda,  \bar b,\bar A)\in \dom S$ if and only if Condition 1.1 holds at $(\bar \lambda,  \bar b,\bar A)$ under the full row rank assumption on $\bar A$. 
It was noticed in \cite{Bello-Cruz2022,Fadili2024,Gilbert17,Liu-Y2021} that this assumption is not necessary. }

Now we discuss Condition 1.2. 
Starting from  any $(I^+,I^0,I^-)\in \mathcal{K}(\bar \lambda, \bar b, \bar A)$ with  $\bar\lambda>0$, we show  in the following a method to construct a sequence of matrices   $A^\nu\to \bar A$  such that   $\rank (A^\nu_{I^+\cup I^-})=\rank (\bar A_{I^+\cup I^-})$   and 
$\mathcal{K}(\bar \lambda, \bar b,  A^\nu)=\{ (I^+,I^0,I^-) \}$ for all $\nu$. In our argument, the matrices $A^\nu$ are constructed by rotating and scaling $\bar A$.

\begin{lemma}\label{lem-boundary-to-int}
Let $(I^+,I^0,I^-)\in \mathcal{K}(\bar \lambda, \bar b, \bar A)$ with  $\bar\lambda>0$. There exists a sequence $A^\nu\to\bar A$ such that $\rank (A^\nu_{I^+\cup I^-})=\rank (\bar A_{I^+\cup I^-})$ and 
$\mathcal{K}(\bar \lambda, \bar b,   A^\nu)=\{(I^+,I^0,I^-)\}$ for all $\nu$.
\end{lemma}
\begin{proof}
Given that $(I^+,I^0,I^-)\in \mathcal{K}(\bar \lambda, \bar b, \bar A)$, we get from the definition of $\mathcal{K}$ (cf. \eqref{def-active}) that  $(I^+,I^0,I^-)\in\mathcal{J}(\bar A)$ and there is some $\bar x\in \R^n$ such that 
$(\bar\lambda, \bar b, \bar x)\in \mathbb{G}(\bar A, (I^+,I^0,I^-))$. Then by Proposition \ref{prop-gphS-lambda-b} and the definition of $\mathbb{G}(\cdot, (I^+,I^0,I^-))$ (cf. \eqref{def-sf-new}), we have  $\bar x\in S(\bar \lambda, \bar b, \bar A)\cap \cl W(I^+,I^0,I^-)$  and 
$$
\bar A_i^\top(\bar b-\bar A\bar x)=\bar\lambda\ \forall i\in I^+,\quad |\bar A_i^\top(\bar b-\bar A\bar x)|\leq \bar\lambda\ \forall i\in I^0,\quad \bar A_i^\top(\bar b-\bar A\bar x)=-\bar\lambda\ \forall i\in I^-.
$$
This implies by \eqref{def-FFF} that  $\bar y:=(\bar b-\bar A\bar x)/\bar\lambda$ belongs to the polyhedral face $\mathbb{F}(\bar A,(I^+,I^0,I^-))$ of $\mathbb{Y}(\bar A)$.  So there are some sequences  $\hat x^\nu\to\bar x$ and $\hat y^\nu\to\bar y$ such that $\hat x^\nu\in W(I^+,I^0,I^-)$ and $\hat y^\nu\in\ri \mathbb{F}(\bar A,(I^+,I^0,I^-))$ for all $\nu$. Let $\bar z:=\bar y+ \bar A(\bar x/\bar\lambda)= \bar b/\bar\lambda$ and $z^\nu:=\hat y^\nu+ \bar A(\hat x^\nu /\bar\lambda)$ for all $\nu$. Then we get $z^\nu\to \bar z$.

In the case of $\bar b=0$,  we have $\bar x=0$ due to $\bar x\in S(\bar\lambda, \bar b, \bar A)$. Then in light of $\bar \lambda>0$,  we get  $(I^+,I^0,I^-)=(\emptyset, I, \emptyset)$ and  $(\bar\lambda, \bar b)\in  \inte P_{1,2}(\mathbb{G} (\bar A, (I^+,  I^0, I^-)))$ (cf. Lemma \ref{lem-basic-pro} (i)). By Corollary \ref{cor-domain-partition} and the definition of $\mathcal{K}$, we have $\mathcal{K}(\bar \lambda, \bar b,   \bar A)=\{(I^+,I^0,I^-)\}$. 
In this case,  the result holds trivially with $A^\nu=\bar A$ for all $\nu$. So we only need to consider the case that $\bar b\not=0$ in what follows.

As $\bar z\not=0$ due to $\bar b\not=0$, we assume that $z^\nu\not=0$ for all $\nu$.   Clearly, we have  $\gamma^\nu:=\|z^\nu\|/\|\bar z\|\to 1$ and
\begin{equation}\label{jiaomaya111'}
   x^\nu:= \frac{\hat x^\nu}{(\gamma^\nu)^2}\in W(I^+, I^0, I^-)\quad \forall \nu.
\end{equation}
Given that $z^\nu/\|z^\nu\|\to \bar z/\|\bar z\|$,  there is a sequence of rotation matrices $M^\nu$ that  converges  to the identity matrix such that  $M^\nu (z^\nu/\|z^\nu\|)=\bar z/\|\bar z\|$ for all $\nu$ or equivalently  
\begin{equation}\label{jiaomaya222'}
 \bar z=\frac{M^\nu z^\nu}{\gamma^\nu}\quad \forall \nu.
\end{equation}
Let $A^\nu:=\gamma^\nu M^\nu\bar A$ and $y^\nu:=M^\nu\hat y^\nu/\gamma^\nu$ for all $\nu$. Clearly, we have $A^\nu\to \bar A$ and $\rank (A^\nu_{I^+\cup I^-})=\rank (\bar A_{I^+\cup I^-})$ for all $\nu$. Given  that $\hat y^\nu\in\ri \mathbb{F}(\bar A,(I^+,I^0,I^-))$ for all $\nu$, we get from \eqref{xiangduineibu} that the following hold for all $\nu$: $\bar A_i^\top\hat y^\nu=1$ for all $i\in I^+$, $|\bar A_i^\top\hat y^\nu|<1$ for all $i\in I^0$ and $\bar A_i^\top\hat y^\nu=-1$ for all $i\in I^-$. 
Since $(A_i^\nu)^\top y^\nu=\gamma^\nu\bar A_i^\top(M^\nu)^\top M^\nu\hat y^\nu/\gamma^\nu=\bar A_i^\top\hat y^\nu$ for all $i$ and $\nu$,  we get from the definition of $\mathcal{J}$ (cf.\eqref{def-partition-faces}) and  \eqref{xiangduineibu} that  $(I^+,I^0,I^-)\in\mathcal J(A^\nu)$ and $y^\nu\in\ri\mathbb F(A^\nu,(I^+,I^0,I^-))$ for all $\nu$. Moreover, we deduce from \eqref{jiaomaya222'} and the equalities that  $\bar z=\bar b/\bar\lambda$, $z^\nu=\hat y^\nu+ \bar A(\hat x^\nu /\bar\lambda)$, $y^\nu=M^\nu\hat y^\nu/\gamma^\nu$, $x^\nu=\hat x^\nu/(\gamma^\nu)^2$ and $A^\nu=\gamma^\nu M^\nu\bar A$ that
$$
\bar b=\bar\lambda \bar z=\frac{\bar\lambda M^\nu z^\nu}{\gamma^\nu}=\frac{\bar\lambda M^\nu\hat y^\nu}{\gamma^\nu}+\frac{M^\nu\bar A \hat x^\nu}{\gamma^\nu}=\bar\lambda y^\nu+A^\nu x^\nu\quad\forall \nu.
$$
This, combined with \eqref{jiaomaya111'} and the fact that $y^\nu\in\ri\mathbb F(A^\nu,(I^+,I^0,I^-))$ for all $\nu$, implies by Lemma \ref{lem-basic-pro}(i) that $(\bar\lambda,\bar b)\in\inte P_{1,2}(\mathbb G(A^\nu,(I^+,I^0,I^-)))$ for all $\nu$. By Corollary \ref{cor-domain-partition}, we get  $\mathcal{K}(\bar \lambda, \bar b,   A^\nu)=\{(I^+,I^0,I^-)\}$ for all $\nu$. This completes the proof. 
\end{proof}

We now apply Lemma \ref{lem-boundary-to-int} to obtain the following result, which is of independent interest itself, and says that the Lipschitz continuity of $S$ cannot be guaranteed around $(\bar\lambda, \bar b, \bar A)\in \dom S$ 
whenever  $S(\cdot,\cdot, \bar A)$ is not single-valued  around  $(\bar \lambda, \bar b)$ (cf. Corollary \ref{cor-big} (ii)).

\begin{lemma}\label{lem-no-lip}
Assume that $\mathcal{K}(\bar \lambda, \bar b, \bar A)\backslash\mathcal{J}^0(\bar A)\not=\emptyset$.
The following hold:
\begin{description}
    \item[(i)] If  $\bar\lambda>0$,  $S(\bar \lambda, \bar b, \cdot)$ is not Lipschitz continuous around $\bar A$.    
    \item[(ii)] If   $\bar \lambda=0$,  $S$ is not Lipschitz continuous around $(\bar \lambda, \bar b, \bar A)$ relative to $\dom S$. 
\end{description}
\end{lemma}
\begin{proof}  Let $(I^+, I^0, I^-)\in \mathcal{K}(\bar \lambda, \bar b, \bar A)\backslash\mathcal{J}^0(\bar A)$. Then by definition we get   $(I^+, I^0, I^-)\in \mathcal{J}(\bar A)$ and $(\bar\lambda, \bar b)\in  P_{1,2}(\mathbb{G} (\bar A, (I^+, I^0, I^-)))$.     Given that $(I^+,I^0,I^-)\not\in \mathcal{J}^0(\bar A)$, we get from the definition of $\mathcal{J}^0$ (cf. \eqref{def-J0}) that  $I^+\cup I^-\not=\emptyset$ 
and that the columns $\bar A_i$ with $i\in I^+\cup I^-$ are linearly   dependence. We claim that $I^+\cup I^-$ cannot be a singleton, for otherwise the equality $I^+\cup I^-=\{k_0\}$ would imply that $\bar A_{k_0}=0$ and hence that $|\bar A_{k_0}^\top y^*|=0$ for all $y^*\in\R^m$, a contradiction to the fact that $(I^+, I^0, I^-)\in \mathcal{J}(\bar A)$. 
So   $I^+\cup I^-$ must contain at least two indices  and there is at least one index $i_0\in I^+\cup I^-$ such that the column $\bar A_{i_0}$ can be expressed as a linear combination of the other  ones. 
That is, there exists some $\mu_i\in \R$ for each  $i\in (I^+\cup I^-)\backslash \{i_0\}$ such that 
\begin{equation}\label{lem-nece-Ai0}
\bar A_{i_0}=\sum\limits_{i\in (I^+\cup I^-)\backslash \{i_0\}}\mu_i\bar A_i.
\end{equation}
By definition of $\mathcal{J}$ (cf.\eqref{def-partition-faces}), there is some $y\in \R^m$ such that 
\[
\bar A^\top_i y=1\quad \forall i\in I^+\quad \mbox{and}\quad \bar A^\top_i y=-1\quad \forall i\in I^-.
\]
This, together with \eqref{lem-nece-Ai0},  implies that 
\begin{equation}\label{guoduguodu}
    \left|\sum\limits_{i\in I^+\backslash \{i_0\}}\mu_i-\sum\limits_{i\in I^-\backslash \{i_0\}}\mu_i\right|=1.
\end{equation}

To show (i), suppose to the contrary that  there are some $\tau>0$ and a neighborhood $\mathcal{V}$ of $\bar A$ such that 
    \begin{equation}\label{zhenxindaistu}
        S(\bar\lambda,\bar b, A)\subset S(\bar\lambda, \bar b, A')+\tau\|A-A'\|\B\quad \forall A, A'\in \mathcal{V}. 
    \end{equation}
    Consider a sequence $A^\nu\to \bar A$  such that $A^\nu\in \mathcal{V}$ with $A^\nu_{i_0}:=(1-\frac{1}{\nu})\bar A_{i_0}$ and $A^\nu_i:=\bar A_i$ for all $i\not=i_0$ and all $\nu$.  We claim that $(I^+, I^0, I^-)\not\in \mathcal{J}(A^\nu)$ for all $\nu$, for otherwise there would exists by definition  some $\tilde y$ such that
\[
(1-\frac{1}{\nu})|\bar A_{i_0}^\top \tilde y|=1,\quad \bar A^\top_i \tilde y=1\quad \forall i\in I^+\backslash\{i_0\}\quad \mbox{and}\quad \bar A^\top_i \tilde y=-1\quad \forall i\in I^-\backslash\{i_0\},
\]
which is impossible due to \eqref{lem-nece-Ai0} and \eqref{guoduguodu}.
    
First we consider the case that $(\bar\lambda, \bar b)\in \inte P_{1,2}(\mathbb{G} (\bar A, (I^+, I^0, I^-)))$. In this case,  we get from Proposition \ref{prop-gphS-lambda-b},  Lemma \ref{lem-basic-pro} (i) and   Lemma \ref{lem-basic-proSF}  that $\bar \lambda>0$ and there is some $\bar x\in S(\bar \lambda,  \bar b, \bar A)$ such that 
\begin{equation}\label{lem-nece-ri-S000}
 \begin{array}{lcl}
    \bar x_i> 0, &\hskip-0.2cm \bar A_i^\top (\bar b-\bar A \bar x)=\bar\lambda &\hskip-0.2cm \forall i\in I^+,\\
     \bar x_i=0, &\hskip-0.2cm  |\bar A_i^\top ( \bar b-\bar A \bar x)|< \bar\lambda&\hskip-0.2cm \forall i\in I^0,\\
    \bar x_i<  0, &\hskip-0.2cm \bar A_i^\top (\bar b-\bar A \bar x)=-\bar\lambda &\hskip-0.2cm \forall i\in I^-.
     \end{array}
\end{equation}
It then follows from the   inclusion \eqref{zhenxindaistu} and the fact that $A^\nu \in \mathcal{V}$ for all $\nu$ that there are some   $x^\nu\in S(\bar \lambda,\bar b,A^\nu)$ and $u^\nu\in \B$ such that $\bar x=x^\nu+\tau\|A^\nu-\bar A\| u^\nu$ for all $\nu$. Thus, we have  $x^\nu=\bar x-\tau\|A^\nu-\bar A\| u^\nu\to \bar x$ due to $A^\nu\to \bar A$ and $\{u^\nu\}$ being bounded.   In light of  $x^\nu\in S(\bar \lambda, \bar b,A^\nu)$ ,  we get from  \eqref{lam-b-xs'} in Proposition \ref{prop-gphS-lambda-b} that  $(A^\nu)^\top(\bar b-A^\nu x^\nu)\in \bar \lambda\partial\|\cdot\|_1(x^\nu)$) for all $\nu$.     
 Then for all $\nu$ sufficiently large, we get  from \eqref{lem-nece-ri-S000} that  $ x_i^\nu>0$ for all $i\in I^+$ and $x_i^\nu<0$ for all $i\in I^-$, and then from the formula \eqref{norm-1-subdiff} for the subdifferential set of $\|\cdot\|_1$ that 
 \begin{equation}\label{tmckeai}
 (A_i^\nu)^\top ( \bar b-A^\nu  x^\nu)=\bar \lambda\;\; \forall i\in I^+\quad \mbox{and}\quad
 (A_i^\nu)^\top ( \bar b-A^\nu x^\nu)=-\bar \lambda\;\;\forall i\in I^-.
\end{equation}
 Due to $A^\nu \to \bar A$ and $x^\nu\to \bar x$, we get from \eqref{lem-nece-ri-S000} that for all $\nu$ sufficiently large, 
 \[
 |(A_i^\nu)^\top ( \bar b-A^\nu  x^\nu)|<\bar \lambda\;\;\forall i\in I^0,
 \]
 which, together with \eqref{tmckeai}, implies by definition  that $(I^+, I^0, I^-)\in \mathcal{J}(A^\nu)$ for all $\nu$ sufficiently large,  
a contradiction to our previous claim  that $(I^+, I^0, I^-)\not\in \mathcal{J}(A^\nu)$ for all $\nu$. 

Next we consider the case that $(\bar\lambda, \bar b)\in \bdry P_{1,2}(\mathbb{G} (\bar A, (I^+, I^0, I^-)))$.
In this case, due to  $(I^+, I^0, I^-)\}\in \mathcal{K}(\bar \lambda, \bar b, \bar A)$ and $\bar\lambda>0$,  we get from Lemma \ref{lem-boundary-to-int} that there is some $A^\nu\to \bar A$ such that $\rank (A^\nu_{I^+\cup I^-})=\rank (\bar A_{I^+\cup I^-})$ and 
$\mathcal{K}(\bar \lambda, \bar b, A^\nu)=\{(I^+, I^0, I^-)\}$ for all $\nu$. 
By  Remark \ref{rem-kkk}, 
we get $(\bar\lambda, \bar b)\in  \inte P_{1,2}(\mathbb{G} ( A^\nu, ( I^+,  I^0,  I^-)))$. Noting that  $\rank (A^\nu_{I^+\cup I^-})=\rank (\bar A_{I^+\cup I^-})$ and  the columns $\bar A_i$ of $\bar A$  with $i\in I^+\cup I^-$  are linearly  dependent due to $(I^+,I^0,I^-)\not\in \mathcal{J}^0(\bar A)$,   we get  $( I^+,  I^0,  I^-)\in \mathcal{J}(A^\nu)\backslash\mathcal{J}^0(A^\nu)$ for all $\nu$.  So following the same arguments as in the first case, we assert that $S(\bar\lambda,\bar b,\cdot)$ is not locally Lipschitz continuous around $A^\nu$ for all $\nu$. 
However, we have $A^\nu\in \mathcal{V}$ for all $\nu$ sufficiently large due to $A^\nu\to \bar A$. This implies by  
 \eqref{zhenxindaistu} that $S(\bar\lambda,\bar b,\cdot)$ is locally Lipschitz continuous around $A^\nu$ for  all $\nu$ sufficiently large, a contradiction. 
The contradictions obtained in the above two cases indicate  that $S(\bar \lambda, \bar b, \cdot)$ is indeed not Lipschitz continuous 
around $\bar A$ whenever $\bar\lambda>0$. 

To show (ii), suppose to the contrary that there are some $\tau>0$ and some    neighborhoods $\mathcal{U}$ of $(\bar\lambda, \bar b)$ and  $\mathcal{V}$ of $\bar A$   such that the following holds for all $((\lambda, b), A), ((\lambda', b'), A')\in  \mathcal{U}\times \mathcal{V} \cap \dom S$: 
        \begin{equation}\label{kekuozhan}
                S(\lambda,b, A)\subset S(\lambda',b', A')+\tau\|(\lambda,b, A)-(\lambda',b', A')\|\B.
        \end{equation}
Noting that we have $(\bar\lambda, \bar b)\in \bdry P_{1,2}(\mathbb{G} (\bar A, (I^+, I^0, I^-)))$ when $\bar\lambda=0$, we get from  Lemma \ref{lem-basic-pro} (i)   that $\inte  P_{1,2}(\mathbb{G} (\bar A, (I^+, I^0, I^-)))\cap \mathcal{U}\not=\emptyset$. 
Let $(\lambda^*, b^*)\in \inte  P_{1,2}(\mathbb{G} (\bar A, (I^+, I^0, I^-)))\cap \mathcal{U}$. By Lemma \ref{lem-basic-pro} (i) again, we have $\lambda^*>0$. 
Then we get from (i) that  $S(\lambda^*,b^*, \cdot)$ is not Lipschitz continuous around $\bar A$. However, it follows from \eqref{kekuozhan} that $S(\lambda^*,b^*, A)\subset S(\lambda^*, b^*, A')+\tau\| A-A'\|\B$ for all $A,A'\in \mathcal{V}$, which indicates that $S(\lambda^*,b^*, \cdot)$ is Lipschitz continuous around $\bar A$, a contradiction. 
\end{proof}

{
Let $(\bar\lambda, \bar b,\bar A)\in \dom S$ with $\bar\lambda>0$. The equivalences among Condition 1.2, the single-valuedness and Lipschitz continuity of $S$ around $(\bar\lambda, \bar b,\bar A)$, the single-valuedness of $S$ around $(\bar\lambda, \bar b,\bar A)$, and others have been established in \cite[Theorem 3.8]{Cui2024} (see also \cite[Theorem 3.13]{Nghia2024}) for some more general models by employing Robinson's strong regularity on the dual problem (see also \cite[Remark 4.16 (a)]{Berk2023bh} for the Lasso.) Below by making use of the explicit representation of $S$ given in Corollary \ref{cor-big}, we show the equivalences among Condition 1.2, the Lipschitz continuity of $S$ around $(\bar\lambda, \bar b,\bar A)$ and the piecewise linearity of $S(\cdot,\cdot,\bar A)$ around $(\bar\lambda, \bar b)$.} 

{
\begin{proposition}\label{prop-local-lip-condition4.2}
Let $(\bar\lambda, \bar b,\bar A)\in \dom S$ with $\bar\lambda>0$. The following are equivalent: 
\begin{description}
\item[(i)]   Condition 1.2  holds at  $(\bar\lambda, \bar b,\bar A)$. 
\item[(ii)] $S$ is Lipschitz continuous around $(\bar\lambda, \bar b,\bar A)$.
\item[(iii)] $S(\bar\lambda, \bar b, \cdot)$ is Lipschitz continuous around $\bar A$.
\item[(iv)] $S(\cdot, \cdot, \bar A)$ is single-valued and  linear  on  each polyhedral cone $P_{1,2}(\mathbb{G}(\bar A, (I^+,I^0, I^-)))$ containing $(\bar \lambda, \bar b)$, where  $(I^+, I^0, I^-)\in \mathcal{J}(\bar A)$. 
\end{description}
\end{proposition}}
\begin{proof}  Note that Condition 1.2 holds  at $(\lambda, b, \bar A)$ for some $\bar x\in S(\lambda, b, \bar A)$  if and only if  \cite[condition (45)]{Cui2024} holds  with $g:=\|\cdot\|_1$ for some  $\bar z:=\frac{1}{\bar \lambda}\bar A^\top(b-\bar A\bar x)$, i.e., 
\[  
\ker \bar A\cap {\rm par}\;\partial \|\cdot\|_\infty (\bar z)=\{0\},
\]
where  ${\rm par}\;\partial \|\cdot\|_\infty (\bar z)$ denotes the parallel subspace of the subdifferential set $\partial \|\cdot\|_\infty (\bar z)$, which coincides with $\{x\in\R^n\mid x_i=0\;\forall i\;\mbox{with}\; |\bar z_i|<1\}$. 
 In light of this equivalence, the implication from (i) to (ii) can be found in  \cite[Theorem 3.5]{Cui2024}, see also  \cite[Remark 4.16 (a)]{Berk2023bh}.   As the implication from (ii) to (iii) is trivial, it suffices to show the   implications from (iii) to (iv), and from (iv) to (i). 


To  show the implication from (iii) to (iv), suppose to the contrary that (iv) does not hold, i.e., there is some $(I^+, I^0, I^-)\in \mathcal{J}(\bar A)$ such that  $(\bar \lambda, \bar b)\in P_{1,2}(\mathbb{G}(\bar A, (I^+,I^0, I^-)))$  and    $S(\cdot, \cdot, \bar A)$ is not  linear on  $P_{1,2}(\mathbb{G}(\bar A, (I^+,I^0, I^-)))$. By Corollary \ref{cor-big} (ii), we know that $(I^+, I^0, I^-)\not\in \mathcal{J}^0(\bar A)$. Given that $\bar\lambda>0$, we get from Lemma \ref{lem-no-lip} that $S(\bar\lambda, \bar b, \cdot)$ is not Lipschitz continuous around $\bar A$, a contradiction to (iii).  

To get the implication from (iv) to (i), let $\bar x\in S(\bar \lambda, \bar b,\bar A)$.  Given that $\bar \lambda>0$ and that  
  $\bar A^\top(\bar b-\bar A\bar x)\in \bar \lambda\partial\|\cdot\|_1(\bar x)$ due to \eqref{lam-b-xs'} in Proposition \ref{prop-gphS-lambda-b}, we get from the formula \eqref{norm-1-subdiff} for $\partial\|\cdot\|_1$  that  $\|\bar A^\top(\bar b-\bar A\bar x)\|_\infty\leq \bar \lambda$ and that  
  $\bar x_i\geq 0$ for all $i\in I^+$, $\bar x_i=0$ for all $i\in I^0$, and $\bar x_i\leq 0$ for all $i\in I^-$,  where $I^+:=\{i\in I\mid \bar A_i^\top(\bar b-\bar A\bar x)=\bar \lambda\}$, $I^0:=\{i\in I\mid |\bar A_i^\top(\bar b-\bar A\bar x)|<\bar \lambda\}$ and $I^-:=\{i\in I\mid \bar A_i^\top(\bar b-\bar A\bar x)=-\bar \lambda\}$. 
  Moreover,  we get $(I^+,I^0,I^-)\in\mathcal J(\bar A)$ by definition.    It then follows from  the definitions of $\mathbb G(\cdot,(I^+,I^0,I^-))$ (cf. \eqref{def-sf-new})  and $P_{1,2}$ that  $(\bar \lambda,\bar b,\bar x)\in\mathbb G(\bar A,(I^+,I^0,I^-))$ and   $(\bar \lambda,\bar b)\in P_{1,2}(\mathbb G(\bar A,(I^+,I^0,I^-)))$. By (iv) and  Corollary \ref{cor-big} (ii), we get  $(I^+,I^0,I^-)\in \mathcal J^0(\bar A)$, which implies by definition  that the columns $\bar A_i$ with $i\in I^+\cup I^-$ are linearly independent. Thus, Condition 1.2  holds at  $(\bar\lambda, \bar b,\bar A)$ for $\bar x\in S(\bar \lambda, \bar b,\bar A)$ by definition.   This  completes the proof.
\end{proof}

It was pointed out by \cite[Remark 4.16 (b)]{Berk2023bh} that if Condition 1.3  holds at  some $(\bar\lambda, \bar b,\bar A)\in \dom S$ with $\bar\lambda>0$, then $S$ is single-valued and Lipschitz continuous around $(\bar\lambda, \bar b,\bar A)$. However, under the same condition,  we are able to show that $S$ is  actually   $\mathcal{C}^\infty$ around $(\bar\lambda, \bar b,\bar A)$. This is done by utilizing  our explicit representation of $S$ presented in Corollary \ref{cor-big} (ii) to show that $S$ coincides locally with a $\mathcal{C}^\infty$ implicit function determined by a system of linear equations. 

{
\begin{proposition}\label{prop-local-lip-condition4.3}
Let $(\bar\lambda, \bar b,\bar A)\in \dom S$ with $\bar\lambda>0$. The following are equivalent:  
\begin{description}
    \item[(i)] Condition 1.3  holds at  $(\bar\lambda, \bar b,\bar A)$. 
    \item[(ii)]  $S$ is   $\mathcal{C}^\infty$ around $(\bar\lambda, \bar b,\bar A)$.
    \item[(iii)] $S(\cdot, \cdot, \bar  A)$ is  single-valued and linear  on the unique polyhedral cone  $P_{1,2}(\mathbb{G}(\bar  A, (I^+,I^0, I^-)))$ 
    having $(\bar \lambda, \bar  b)$ in its interior. 
\end{description} 
\end{proposition}}
\begin{proof}  
To show    the implication from (i) to (ii), assume that Condition 1.3 holds at $(\bar \lambda, \bar b, \bar A)$ for some $\bar x\in S(\bar \lambda, \bar b, \bar A)$.  Given that $\bar x\in S(\bar \lambda, \bar b, \bar A)$, we get from \eqref{lam-b-xs'} in Proposition \ref{prop-gphS-lambda-b} that $\bar A^\top(\bar b-\bar A\bar x)\in  \bar \lambda\partial\|\cdot\|_1(\bar x)$. It then  follows from  the formula (\ref{norm-1-subdiff}) for  $\partial\|\cdot\|_1$ that $\bar A_i^\top (\bar b-\bar A\bar x)=\bar \lambda$ for all $i\in I^+(\bar x)$ and $\bar A_i^\top (\bar b-\bar A\bar x)=-\bar \lambda$ for all $i\in I^-(\bar x)$. By  Condition 1.3  at $(\bar \lambda, \bar b, \bar A)$ for $\bar x$,  we have $|\bar A_i^\top (\bar b-\bar A\bar x)|< \bar \lambda$  for all $i\in I^0(\bar x)$.
This implies by definition that  $(I^+(\bar x), I^0(\bar x), I^-(\bar x))\in \mathcal{J}(\bar A)$ and hence that $(I^+(\bar x), I^0(\bar x), I^-(\bar x))\in \mathcal{J}^0(\bar A)$ by Condition 1.3  at $(\bar \lambda, \bar b, \bar A)$ for some $\bar x$.  Moreover, we get   from  Lemmas \ref{lem-basic-proSF} and \ref{lem-basic-pro}  that $(\bar \lambda, \bar b, \bar x)\in \ri \mathbb{G}(\bar A, (I^+(\bar x), I^0(\bar x), I^-(\bar x)))$ and hence that $(\bar \lambda, \bar b)\in\inte P_{1,2}(\mathbb{G}(\bar A, (I^+(\bar x), I^0(\bar x), I^-(\bar x))))$.  This implies by  Remark  \ref{rem-kkk} that  $\mathcal{K}(\bar\lambda, \bar b, \bar A)=\{(I^+(\bar x),I^0(\bar x),I^-(\bar x))\}$.   Then by   Lemmas \ref{lem-J0A} (i) and \ref{lem-active-outer-conti} (i), 
there is some open neighborhood $\mathcal{W}$ of $(\bar\lambda,  \bar b, \bar A)$ such that
   \[
      \mathcal{K}(\lambda, b, A)=\{(I^+(\bar x),I^0(\bar x),I^-(\bar x))\}\subset \mathcal{J}^0(A)\quad \forall (\lambda, b, A)\in\mathcal{W}.
   \]
This, together with Corollary \ref{cor-big} (ii), implies that $S(\lambda, b, A)=\{x\}$ for $(\lambda, b, A)\in \mathcal{W}$,  
    where $x$ is the unique solution to the system of linear equations
          \begin{equation}\label{equ-around-single-valued}
     \left\{
     \begin{array}{rlll}
 A_i^\top  Ax&=&A_i^\top b-\lambda  & \forall i\in I^+(\bar x),\\
     x_i&=&0 & \forall i\in I^0(\bar x),\\
 A_i^\top  Ax&=&A_i^\top b+\lambda & \forall i\in I^-(\bar x).
     \end{array}
     \right.
    \end{equation}  
    That is, on the open neighborhood $\mathcal{W}$ of $(\bar\lambda,   \bar b,  \bar A)$, $S$ is single-valued and    coincides with the implicit function determined by the linear  system \eqref{equ-around-single-valued}. Therefore, by Theorem 3 in Chapter 8 of the reference  \cite{Magnus1999}, $S$ is $\mathcal{C}^\infty$  on $\mathcal{W}$.   This verifies (ii).

To show the implication from (ii) to (iii),  we note from (ii)  that $S$ is Lipschitz continuous around $(\bar\lambda, \bar b, \bar A)$. Then by Proposition \ref{prop-local-lip-condition4.2}, $S(\cdot, \cdot, \bar A)$ is single-valued and  linear  on  each polyhedral cone $P_{1,2}(\mathbb{G}(\bar A, (I^+,I^0, I^-)))$ containing $(\bar \lambda, \bar b)$, where  $(I^+, I^0, I^-)\in \mathcal{J}(\bar A)$. 
Therefore, to show (iii), it suffices to show the  aforementioned polyhedral cone is unique and contains $(\bar \lambda, \bar  b)$ in its interior. 
  Suppose to the contrary this is not the case. Then by Remark  \ref{rem-kkk},  
  there exist   at least two distinct index partitions $(I^+, I^0, I^-)$ and $(\bar I^+, \bar I^0, \bar I^-)$ such that 
\begin{equation}\label{ggzj}
    (\bar \lambda, \bar b)\in \bdry P_{1,2}(\mathbb{G}(\bar A,(I^+, I^0, I^-)))\cap \bdry P_{1,2}(\mathbb{G}(\bar A,(\bar I^+, \bar I^0, \bar I^-))).
\end{equation}
By Corollary \ref{cor-big} (iv), any line segment joining $(\lambda', b', S(\lambda',b', \bar A))$ 
and $(\lambda'', b'', S(\lambda'',b'', \bar A))$ with  $(\lambda', b')\in \inte P_{1,2}(\mathbb{G}(\bar A,(I^+, I^0, I^-)))$ and $(\lambda'', b'')\in \inte P_{1,2}(\mathbb{G}(\bar A,(\bar I^+, \bar I^0, \bar I^-)))$  does not lie entirely on the graph of $S(\cdot, \cdot, \bar A)$. Combining this with the fact that  $S(\cdot, \cdot, \bar A)$  is a piecewise linear function on the union of $ P_{1,2}(\mathbb{G}(\bar A,(I^+, I^0, I^-)))$ and $P_{1,2}(\mathbb{G}(\bar A,(\bar I^+, \bar I^0, \bar I^-)))$, we confirm  that  $S(\cdot, \cdot, \bar A)$ is not differentiable at $(\bar \lambda, \bar b)$, contradicting to (ii).

 To show the implication from (iii) to (i),  assume that  (iii) holds with some $(I^+, I^0, I^-)\in \mathcal{J}(\bar A)$. Given that $S(\cdot, \cdot, \bar  A)$ is a linear function on   $P_{1,2}(\mathbb{G}(\bar  A, (I^+,I^0, I^-)))$, we get from  Corollary \ref{cor-big}  (ii) that $(I^+, I^0, I^-)\in \mathcal{J}^0(\bar A)$. Then by definition,   the columns $\bar A_i$  with $i\in I^+\cup I^-$ are linearly independent. In view of  
$(\bar \lambda,\bar b)\in\inte P_{1,2}(\mathbb{G}(\bar A, (I^+,I^0, I^-)))$, we  get from Lemma \ref{lem-basic-proSF}, Lemma \ref{lem-basic-pro} (i) and Proposition \ref{prop-gphS-lambda-b} that there exists some $\bar x\in S(\bar\lambda, \bar b, \bar A)$ such that $I^+=I^+(\bar x)$, $ I^-=I^-(\bar x)$, and   $|A_i^\top (\bar b-\bar A\bar x)|< \bar \lambda$ for all $i\in I^0(\bar x)$.   This verifies  Condition 1.3   at $(  \bar \lambda,  \bar b, \bar A)$ for $\bar x$ by definition.  This completes the proof.
\end{proof} 

\begin{remark}\label{lem-vip-linear-piecewise-linear}
Let $(\lambda,b,A)\in \dom S$ with $\lambda>0$. Under the assumption that Condition 1.2 holds at  $(\lambda,b,A)$,  Proposition \ref{prop-local-lip-condition4.2} (iv), together with Remark \ref{rem-kkk},  implies that $S(\cdot, \cdot, A)$ is single-valued and piecewise linear  around $(\lambda,  b)$. Under Condition 1.3 being held at  $(\lambda,b,A)$, Proposition \ref{prop-local-lip-condition4.3} (iii)  yields that   $S(\cdot, \cdot, A)$ is single-valued and linear around $(\lambda,  b)$.   Notably,  Berk et al. \cite{Berk2023bh} established the directional differentiability and continuous differentiability of $S(\cdot, \cdot, A)$  under Conditions 1.2 and 1.3,  respectively. 
 It is worth noting that when Condition 1.2 and 1.3 hold, respectively,  at $(\lambda, b, A)$,  the piecewise linearity and  linearity of $S(\cdot, \cdot, A)$ around $(\lambda, b)$  directly imply the directional differentiability and the continuous differentiability of $S(\cdot, \cdot, A)$  around $(\lambda, b)$, respectively. 
\end{remark}

The following example shows that $S(\lambda, b, A)$ is in general not  Lipschitz continuous when $A$ varies along a line. Therefore when the Lipschitz continuity of $S$ is studied with $A$ being taken as a   parameter as well, we need only to look for the local change of $A$.
 
\begin{example} (Cont. of Example \ref{exam-nice})
Consider a line $\mathcal{L}:=\{(\bar \lambda, \bar b, M(t))\mid t\in \R\}$ that lies in $\dom S=\R_+\times \R^m\times \R^{m\times n}$, where    $\bar\lambda=2$, $\bar b=(4, 1)^\top$ and  $M:\R\rightarrow \R^{m\times n}$ is a matrix-valued   mapping  defined by 
\[
t\mapsto A+t\left[
\begin{array}{ccc}
  1 & 0 & 0 \\
  0 & 0 & 0
\end{array}
\right],
\] 
where $A$ is given as in Example \ref{exam-nice}. 
Following the same process as in Example \ref{exam-nice}, we can calculate $P_{1,2}(\mathbb{G}(M(t), (I^+,I^0, I^-)))$   for all $t\in \R$ and $(I^+,I^0, I^-)\in \mathcal{J}(M(t))$. Then in light of  the definition of  $\mathcal{K}$ (cf. \eqref{def-active}), Lemma \ref{lem-basic-pro} (i)-(ii),  and  Corollary \ref{cor-big} (i)-(ii), we are able to calculate the values of $\mathcal{K}$ and $S$ at $(\bar\lambda, \bar b, M(t))$ for all $t\in \R$. In Table \ref{table-along-line},  we summarize the details of  the mappings    $\mathcal{K}$ and $S$    when $t$ varies.  
By the expressions of $S(\bar\lambda, \bar b, M(t))$ for all $t\in \R$, we observe that $S$ displays non-polyhedral structure  when $t$  changes linearly on $(-\infty, -2)\cup (0,  \infty)$ and  that $S$ is not Lipschitz continuous on   $\mathcal{L}$. 
\begin{table}[ht]
\begin{tabular}{cll}
\hline
$t$ & $\mathcal{K}(\bar\lambda, \bar b, M(t))$ & $S(\bar\lambda, \bar b, M(t))$  \\ 
  \hline
 $(-\infty, -2)$& $\left\{(\{2\}, \{3\}, \{1\}),\;(\emptyset,\{2,3\}, \{1\})\right\}$ & $\left\{\left(\frac{2(3+2t)}{(1+t)^2},0,0\right)^\top \right\} $ \\[0.2cm]
    \hline
       $-2$ & $\{ (\{2,3\}, \emptyset, \{1\})\}$ & $\left\{(2\alpha-2,\alpha,\alpha)^\top\mid \alpha\in [0, 1] \right\} $\\[0.2cm]
 \hline
   $(-2, 0)$ & $\{ (\{2,3\}, \{1\}, \emptyset)\}$  & $ \{(0, 1, 1)^\top\}  $ \\[0.2cm]
  \hline
    $0$ &$\{ (\{1, 2,3\}, \emptyset, \emptyset)\}$ & $\left\{(2-2\alpha,\alpha,\alpha)^\top\mid  \alpha\in [0, 1]\right\} $ \\[0.2cm]
  \hline
   $(0,  \infty)$ & $\left\{(\{1, 2\}, \{3\}, \emptyset),\;(\{1\},\{2,3\}, \emptyset)\right\}$ & $\left\{\left(\frac{2(1+2t)}{(1+t)^2},0,0\right)^\top \right\}$ \\[0.2cm]
  \hline
  \end{tabular}
            \centering
\caption{The values that the multifunctions  $\mathcal{K}$ and $S$ take along the line $\mathcal{L}$
}
      \label{table-along-line}
\end{table}
\end{example}

Recall that we have shown in Theorem \ref{theo-S-Lips} that for any $\bar A\in \R^{m\times n}$,   $S$ is  Lipschitz continuous  on  $\R_+\times \R^m\times \{\bar A\}$.   By restricting $(\lambda, b)$  to a nonempty,  compact and convex subset $\mathcal{U}$   of $\R_{++}\times \R^m$, we give a characterization for the Lipschitz continuity of $S$ on $\mathcal{U}\times \mathcal{V}$ for some neighborhood $\mathcal{V}$ of $\bar A$ by virtue of Proposition \ref{prop-local-lip-condition4.2}.

 \begin{theorem}\label{theo-lip-global-equi}
   Let $\bar A\in \R^{m\times n}$  and  $\mathcal{U}$ be a nonempty,  compact and convex subset  of $\R_{++}\times \R^m$.
  Then  $S$ is    Lipschitz continuous on $\mathcal{U}\times \mathcal{V}$ for some neighborhood $\mathcal{V}$ of $\bar A$   if and only if Condition 1.2 holds at $(\lambda, b, \bar A)$ for all $(\lambda, b)\in \mathcal{U}$. 
\end{theorem}
\begin{proof}  We first show the `if' part. Given that Condition 1.2 holds at $(\lambda, b, \bar A)$ for all $(\lambda, b)\in \mathcal{U}$, we get from Proposition \ref{prop-local-lip-condition4.2}  that  for all $(\lambda, b)\in \mathcal{U}$,   $S(\cdot, \cdot, \bar A)$ is single-valued and linear on each polyhedral cone $P_{1,2}(\mathbb{G}(\bar A, (I^+,I^0, I^-)))$ containing $(\lambda, b)$, where $(I^+, I^0, I^-)\in \mathcal{J}(\bar A)$.  By Corollary \ref{cor-big} (ii)  and the definition of $\mathcal{K}$, we get $\mathcal{K}( \lambda,   b,  \bar A)\subset \mathcal{J}^0(\bar A)$ for all $(\lambda, b)\in \mathcal{U}$.  

We claim   that $\mathcal{K}(\lambda, b,   A)\subset \mathcal{J}^0(A)$ for all $(\lambda, b, A)\in \mathcal{U}\times \B(\bar A, \alpha)$ for all $\alpha>0$ small enough. Suppose to the contrary that there is some sequence $\{(\lambda^\nu, b^\nu,A^\nu)\}\subset \mathcal{U}\times \R^{m\times n}$ such that $A^\nu\to \bar A$ and $\mathcal{K}(\lambda^\nu, b^\nu,   A^\nu)\not\subset \mathcal{J}^0(A^\nu)$ for all $\nu$.  As $\mathcal{U}$ is compact by assumption,  by taking a subsequence if necessary, we may assume that $(\lambda^\nu, b^\nu)\to (\bar \lambda, \bar b)\in \mathcal{U}$, which implies that   $\bar \lambda>0$ due to $\mathcal{U}\subset \R_{++}\times \R^m$ by assumption. 
Given that $\mathcal{K}(\bar \lambda, \bar b,  \bar A)\subset \mathcal{J}^0(\bar A)$ by the previous argument,    we get from Lemma \ref{lem-J0A} (i) and Lemma \ref{lem-active-outer-conti} (i) that   $\mathcal{K}(\lambda^\nu, b^\nu,   A^\nu)\subset \mathcal{J}^0(A^\nu)$ for all $\nu$ sufficiently large, a contradiction to the assumption that $\mathcal{K}(\lambda^\nu, b^\nu,   A^\nu)\not\subset \mathcal{J}^0(A^\nu)$ for all $\nu$. 
 This indicates that   $\mathcal{K}(\lambda, b,   A)\subset \mathcal{J}^0(A)$ for all $(\lambda, b, A)\in \mathcal{U}\times \B(\bar A, \alpha)$ for all $\alpha>0$ small enough.  
   By  Corollary \ref{cor-big} (ii), the definition of $\mathcal{K}$  and Proposition \ref{prop-local-lip-condition4.2} again, we assert  that $S$ is single-valued and Lipschitz continuous around  each   $(\lambda, b, A) \in \mathcal{U}\times \B(\bar A, \alpha)$ for all $\alpha>0$ small enough.  Given that  $\mathcal{U}$ is compact and  convex by assumption, we deduce from \cite[Theorem 9.2]{Rock2009VaAn} that $S$ is Lipschitz continuous on  $\mathcal{U}\times \B(\bar A, \alpha)$ for all $\alpha>0$ small enough.

    To show the `only if' part, assume that $S$ is Lipschitz continuous on $\mathcal{U}\times \mathcal{V}$ for some neighborhood $\mathcal{V}$. This implies that $S(\lambda,b,\cdot)$ is Lipschitz continuous around $\bar A$ for all $(\lambda, b)\in \mathcal{U}$.
    Given that $\lambda>0$ for all   $(\lambda, b)\in \mathcal{U}$ by assumption, we get from  Proposition \ref{prop-local-lip-condition4.2}    that   Condition 1.2 holds at $(\lambda, b, \bar A)$ for all $(\lambda, b)\in \mathcal{U}$.  This completes the proof. 
\end{proof}

When the matrix $\bar A$ has full row rank satisfying certain   condition or full row column rank, and 
$\mathcal{U}$ is chosen from $\mathbb{R}_+\times \mathbb{R}^m$ (allowing $\lambda=0$), we show that $S$ is Lipschitz continuous on $\mathcal{U} \times \mathcal{V}$, where the radius of the neighborhood $\mathcal{V}$ of $\bar{A}$ can be explicitly quantified. 
Again, the explicit representation  obtained for $S$ in Corollary \ref{cor-big} plays a key role in our argument.

\begin{theorem}\label{theo-lip-glob-local-full-rank}
 Let $\bar A\in \R^{m\times n}$  and  $\mathcal{U}$ be a nonempty,  compact and convex subset  of $\R_+\times \R^m$. Assume that one of the following   conditions holds:
\begin{description}
    \item[(i)]$\rank(\bar A)=m$, $\mathcal{J}^0(\bar A)=\mathcal{J}(\bar A)$ and $0<\alpha<\kappa(\bar A)$, where $\kappa(\bar A)$ is defined in \eqref{hmz-banjing}. 
    \item[(ii)]$\rank(\bar A)=n$ and $0<\alpha<\sigma_{\min}(\bar A)$.
\end{description} 
Then $S$ is single-valued and Lipschitz continuous on $\mathcal{U}\times \B(\bar A, \alpha)$.
\end{theorem}
\begin{proof} 
First assume   condition (i). Let $\mathcal{\bar I}\subset \mathcal{I}$ be such that $\mathcal{J}^0(\bar A)=\mathcal{J}(\bar A)=\mathcal{\bar I}$. 
By Corollary \ref{cor-domain-partition}, Corollary \ref{cor-big} (ii) and   Lemma \ref{lem-J0A} (iv), the single-valuedness of $S$ on $\mathcal{U}\times \B(\bar A, \alpha)$ follows immediately. 
It remains to prove   the Lipschitz continuity of $S$ on $\mathcal{U}\times \B(\bar A, \alpha)$.
For every $(I^+,I^0,I^-)\in \mathcal{\bar I}$ and $A\in \B(\bar A, \alpha)$, it follows from  Lemma \ref{lem-J0A}  (iv) that 
$(I^+,I^0,I^-)\in \mathcal{J}^0(A)$ and hence that $A_{I^+\cup I^-}$ is of full column rank. Thus, for every $(I^+,I^0,I^-)\in \mathcal{\bar I}$, the system of linear equations 
       \begin{equation}\label{xiaoxiaosasa}
               \left\{
     \begin{array}{rlll}
 A_i^\top  Ax&=&A_i^\top b-\lambda  & \forall i\in I^+\\
     x_i&=&0 & \forall i\in I^0\\
 A_i^\top  Ax&=&A_i^\top b+\lambda & \forall i\in I^-
     \end{array}
     \right.
       \end{equation}
defines an implicit function, which maps each $(\lambda, b, A)\in \R\times \R^m\times \B(\bar A, \alpha)$  to a unique solution $x$ to \eqref{xiaoxiaosasa}. Every such defined implicit function is, according to Theorem 3 in Chapter 8 of the reference  \cite{Magnus1999},  $\mathcal{C}^\infty$ and thus Lipschitz continuous on the compact and convex subset  $\mathcal{U}\times \B(\bar A, \alpha)$ with  constant $\tau_{(I^+, I^0, I^-)}$, where  $(I^+,I^0,I^-)\in \mathcal{\bar I}$.
Let 
$$\tau:=\max_{(I^+,I^0,I^-)\in \mathcal{\bar I} }\tau_{(I^+, I^0, I^-)}.$$
Then for all $(\lambda', b', A'), (\lambda'', b'', A'')\in \mathcal{U}\times \B(\bar A, \alpha)$ such that   $\mathcal{K}(\lambda', b',A')\cap \mathcal{K}(\lambda'', b'',A'')\not=\emptyset$,  the following  inequality holds:  
\begin{equation}\label{haoyongno3}
    \|S(\lambda', b', A')-S(\lambda'', b'', A'')\|\leq \tau\|(\lambda', b', A')-(\lambda'', b'', A'')\|. 
\end{equation}

  Fix any $(\lambda^0, b^0, A^0)$  and  $(\lambda^1, b^1, A^1)$ in  $\mathcal{U}\times \B(\bar A, \alpha)$, and  define
$$
(\lambda^\theta, b^\theta, A^\theta):=(1-\theta)(\lambda^0, b^0, A^0)+\theta (\lambda^1, b^1, A^1)\quad \forall \theta\in [0, 1]. 
$$
Given that $\mathcal{U}\times \B(\bar A, \alpha)$ is convex by assumption,  we  have $ (\lambda^\theta, b^\theta, A^\theta) \in \mathcal{U}\times \B(\bar A, \alpha)$  for all $\theta\in[0,1]$.  Noting that  $A^\theta\in \B(\bar A, \alpha)$ and $\alpha\in (0, \kappa(\bar A))$, we get from  Lemma \ref{lem-J0A} (iv) that $\rank (A^\theta)=m$  and $\mathcal{J}^0(A^\theta)=\mathcal{J}(A^\theta)$ for all $\theta\in [0, 1]$. 
It then follows from Lemma \ref{lem-active-outer-conti} (ii)  that there is some $\delta_\theta>0$ such that 
 \begin{equation}\label{xiaojubu}
      \mathcal{K}(\lambda,b, A)\subset\mathcal{K}(\lambda^\theta, b^\theta, A^\theta)  \quad  \forall
      (\lambda, b, A)\in \mathbb{B}((\lambda^\theta, b^\theta, A^\theta),\delta_\theta). 
 \end{equation}
Given that the line segment $\{(\lambda^\theta, b^\theta, A^\theta)~|~\theta\in[0,1]\}$   can be covered by   the  balls 
$\inte \mathbb{B}((\lambda^\theta, b^\theta, A^\theta),\delta_\theta)$ with $\theta\in [0, 1]$, we get from the Heine–Borel theorem  (see \cite[Theorem 2.41]{rudin1953}) that  it can also be covered by finitely many such balls $\inte \mathbb{B}((\lambda^{\theta_j}, b^{\theta_j}, A^{\theta_j}),\delta_{\theta_j})$, where $1\leq j\leq r$ and $\theta_j\in [0,1]$.  Without loss of generality, we may assume that $0=\theta_1<\theta_2<\cdots<\theta_r=1$ and   
$$
\inte\mathbb{B}((\lambda^{\theta_j}, b^{\theta_j}, A^{\theta_j}),\delta_{\theta_j})\cap \inte\mathbb{B}( (\lambda^{\theta_{j+1}}, b^{\theta_{j+1}}, A^{\theta_{j+1}}),\delta_{\theta_{j+1}})\not=\emptyset\quad\forall j\in \Gamma :=\{1,\cdots,r-1\}.
$$
Then for each $j\in \Gamma $,  there is some  $\tilde\theta_j\in(\theta_j,\theta_{j+1})$ such that 
$$
(\lambda^{\tilde\theta_j}, b^{\tilde\theta_j}, A^{\tilde\theta_j})\in\inte\mathbb{B}((\lambda^{\theta_j}, b^{\theta_j}, A^{\theta_j}),\delta_{\theta_j})\cap \inte\mathbb{B}( (\lambda^{\theta_{j+1}}, b^{\theta_{j+1}}, A^{\theta_{j+1}}),\delta_{\theta_{j+1}}),
$$
which, together with   \eqref{xiaojubu}, implies  that $\mathcal{K}(\lambda^{\tilde\theta_j},b^{\tilde\theta_j},A^{\tilde\theta_j})\subset \mathcal{K}(\lambda^{\theta_j},b^{\theta_j},A^{\theta_j})\cap \mathcal{K}(\lambda^{\theta_{j+1}},b^{\theta_{j+1}},A^{\theta_{j+1}})$. 
Therefore, for each $j\in \Gamma $,  we get from Corollary \ref{cor-domain-partition} that $\mathcal{K}(\lambda^{\tilde\theta_j},b^{\tilde\theta_j},A^{\tilde\theta_j})\not=\emptyset$ and hence 
\[
  \mathcal{K}(\lambda^{\theta_j},b^{\theta_j},A^{\theta_j})\cap \mathcal{K}(\lambda^{\theta_{j+1}},b^{\theta_{j+1}},A^{\theta_{j+1}})\not=\emptyset.   
\]
In light of $(\lambda^{\theta_j},b^{\theta_j},A^{\theta_j}), (\lambda^{\theta_{j+1}},b^{\theta_{j+1}},A^{\theta_{j+1}})\in \mathcal{U}\times \B(\bar A, \alpha)$, we get from \eqref{haoyongno3} that 
\[
\|S(\lambda^{\theta_j},b^{\theta_j},A^{\theta_j})-S(\lambda^{\theta_{j+1}},b^{\theta_{j+1}},A^{\theta_{j+1}})\|\leq \tau\|(\lambda^{\theta_j},b^{\theta_j},A^{\theta_j})-(\lambda^{\theta_{j+1}},b^{\theta_{j+1}},A^{\theta_{j+1}})\|\quad \forall j\in \Gamma.
\]
This further implies  that 
\[
\begin{aligned}
\|S(\lambda^0,b^0,A^0)-S(\lambda^1,b^1,A^1)\|
\leq&\displaystyle\sum_{i\in \Gamma}\|S(\lambda^{\theta_j},b^{\theta_j},A^{\theta_j})-S(\lambda^{\theta_{j+1}},b^{\theta_{j+1}},A^{\theta_{j+1}})\| \\
\leq&\tau\sum_{i\in \Gamma} \|(\lambda^{\theta_j},b^{\theta_j},A^{\theta_j})-(\lambda^{\theta_{j+1}},b^{\theta_{j+1}},A^{\theta_{j+1}})\|\\
=&\tau\displaystyle \|(\lambda^0,b^0,A^0)-(\lambda^1,b^1,A^1)\|.
\end{aligned}
\]
 Given that $(\lambda^0, b^0, A^0), (\lambda^1, b^1, A^1)\in \mathcal{U}\times \B(\bar A, \alpha)$ are given arbitrarily, we get 
 the Lipschitz continuity of  $S$  on $\mathcal{U}\times \B(\bar A, \alpha)$.  This verifies  the result under condition (i). 

Now assume    condition (ii). Following analogous arguments to those used under condition (i), 
the single-valuedness and Lipschitz continuity of $S$ on $\mathcal{U}\times \B(\bar A, \alpha)$  can be shown   via Lemma \ref{lem-J0A} (iii) and Lemma \ref{lem-active-outer-conti} (iii).  The details of the proof are omitted for the sake of simplicity. 
\end{proof}


\section{Conclusions}

In this paper, we obtained the Lipschitz continuity of the solution multifunction of both the Lasso and the Basis Pursuit problem in a unified way by using the polyhedral theory. We also provided full characterizations such as single-valuedness, piecewise linearity and linearity for the solution multifunctions of the Lasso in terms of some known conditions in the literature. We believe that the results in this paper on global Lipschitz continuity may have some important impacts in learning
global convergence rates of different methods for solving the Lasso (e.g., \cite[Theorem 3.2]{Bello-Cruz2022}) and the BP problem and in global sensitivity analysis.

\bibliographystyle{plain}
\bibliography{References}
\end{document}